\documentclass[10pt,letterpaper]{article}

\newcommand*{\mypart}[1]{\addtocounter{part}{1}\setcounter{page}{1}%
	\setcounter{figure}{0}\setcounter{table}{0}\setcounter{section}{0}%
	\setcounter{footnote}{0}\setcounter{theorem}{0}\setcounter{remark}{0}%
	\setcounter{equation}{0}\thispagestyle{empty}	
	\part*{\centering\normalsize\MakeUppercase{#1: Part \thepart}}%
		\begin{center}\footnotesize TIMOTHY F.~HAVEL\footnote{%
			Retired research staff, MIT; biography given at the end of Part II.}\end{center}} 
\newcommand*{\mytitle}{An Extension of Heron's Formula to Tetrahedra,\\[3pt]%
						and the Projective Nature of Its Zeros}

\usepackage{titlesec}
\titleformat{\section}[runin]%
	{\normalfont\normalsize\bfseries}{\indent\thesection.}{1em}{}[.]
\usepackage[runin]{abstract}
\abslabeldelim{.}

\usepackage[format=plain,justification=centerlast,%
				textfont={small,it},labelfont={small,rm},%
				labelsep=period]{caption}

\usepackage{xr} 
\usepackage{chapterbib}
\sectionbib{\bigskip\centering}{section} 

\newcommand*{\fit}[1][0.05em]{\hspace*{#1}} 
\newcommand*{\dit}[1][1ex]{\rule[-#1]{0pt}{1ex}} 
\newcommand*{\hit}[1][2ex]{\rule{0pt}{#1}} 
\newcommand*{\myref}[2][0]{%
	\ifcase#1 \ref{#2}\or I\fit:\ref{I:#2}\or II\fit:\ref{II:#2}\or %
		III\fit:\ref{III:#2}\or IV:\ref{IV:#2}\fi}

\usepackage{amsmath} 
\usepackage{amssymb} 
\usepackage{amsthm}  

\allowdisplaybreaks 

\newtheoremstyle{pmetheorem} 
	{}{} 
	{\itshape}{\parindent} 
	{\scshape} 
	{.}{0.5em}{} 
\theoremstyle{pmetheorem}
\newtheorem{theorem}{Theorem}
\newtheorem{lemma}[theorem]{Lemma}
\newtheorem{corollary}[theorem]{Corollary}
\newtheorem{proposition}[theorem]{Proposition}
\newtheorem{conjecture}[theorem]{Conjecture}
\newtheorem{definition}[theorem]{Definition}
\theoremstyle{remark} 
\newtheorem{remark}{Remark} 
\usepackage{xpatch} 
\xpatchcmd{\proof}{\hskip\labelsep}{\hskip4\labelsep}{}{}
\makeatletter\apptocmd\th@remark{\let\thm@indent\indent}{}{}\makeatother

\usepackage[helvratio=0.833]{newtxtext} 
\usepackage[nosymbolsc,bigdelims]{newtxmath} 

\usepackage{csquotes}	
\usepackage[a]{esvect}	
\usepackage{graphicx}	
\usepackage{moresize}	
\usepackage{xurl}	

\newcommand*{\TDEF}[1]{\emph{#1}} 
\newcommand*{\EMPH}[1]{\underline{\smash{#1}}} 
\newcommand*{\msf}{\mathsf} \newcommand*{\mbf}{\mathbf} 
\newcommand*{\mrm}{\mathrm} \newcommand*{\mbb}{\mathbb} 
\newcommand*{\mathbi}[1]{\ensuremath\boldsymbol{#1}} 
\newcommand*{\mbi}{\ensuremath\mathbi} 
\newcommand{\Dprod}{
	\mathchoice{\mathbin{\raisebox{0.18ex}{$\scriptstyle\bullet$}}} 
					 {\mathbin{\raisebox{0.18ex}{$\scriptstyle\bullet$}}}
					 {\mathbin{\raisebox{0.09ex}{$\scriptscriptstyle\bullet$}}}
					 {\mathbin{\raisebox{0.09ex}{$\scriptscriptstyle\bullet$}}}}
\newcommand*{\Xprod}{\boldsymbol\times} 
\newcommand*{\Vplus}{\boldsymbol+} 
\newcommand*{\Vdiff}{\boldsymbol-} 
\newcommand*{\Tau}{\ensuremath\text{\small$\boldsymbol{\mathcal T}$}} 
\newcommand*{\Del}{\deltaup\hspace{0.03em}} 
\newcommand*{\CMD}{\ensuremath\scalebox{1.1}[1.1]{$\Delta$}} 
\newcommand*{\OL}{\overline}  
\newcommand*{\fnsize}{\footnotesize} 
\newlength{\upbgdl}
\newcommand*{\rfrac}[2]{\ensuremath\raisebox{0.5\upbgdl}{\Large$\scriptstyle\frac{#1}{#2}$}}
\newcommand*{\sfrac}[2]{\ensuremath\raisebox{0.5\upbgdl}{\small$\displaystyle\frac{#1}{#2}$}}

\newcommand*{\GEO}[1]{\ensuremath\OL{\msf{#1}}} 
\newcommand{\tmv}{ 
	\mathchoice{\displaystyle\vert} 
					 {\mathbin{\fit[-0.07em]\rule{0.67pt}{1.5ex}\fit[-0.07em]}}
					 {\scriptstyle\vert} 
					 {\scriptscriptstyle\vert}} 
\newcommand*{\Lvert}{
	\mathchoice{\raisebox{0.5\upbgdl}{$\displaystyle\big\lvert{}$}}
					 {\raisebox{\upbgdl}{$\lvert{}$}} 
					 {\raisebox{0.5\upbgdl}{$\scriptstyle\lvert{}$}}
					 {\raisebox{0.25\upbgdl}{$\scriptstyle\lvert{}$}}}
\newcommand*{\Rvert}{
	\mathchoice{\raisebox{0.5\upbgdl}{$\displaystyle{\fit}\big\rvert$}}
					 {\raisebox{\upbgdl}{${\fit}\rvert$}} 
					 {\raisebox{0.5\upbgdl}{$\scriptstyle{\fit}\rvert$}}
					 {\raisebox{0.25\upbgdl}{$\scriptstyle{\fit}\rvert$}}}
\newcommand*{\MAG}[1]{\ensuremath\left.\fit[-0.08em]\Lvert\GEO{#1}\Rvert\fit[-0.08em]\right.}
\newcommand*{\VEC}[1]{\ensuremath\vv{\msf{#1}}} 
\newcommand*{\VCP}[2]{\ensuremath\VEC{#1}\Xprod\VEC{#2}} 
\newcommand*{\LVert}{
	\mathchoice{\raisebox{0.5\upbgdl}{$\displaystyle\big\lVert{}$}}
					 {\raisebox{\upbgdl}{$\lVert{}$}} 
					 {\raisebox{0.5\upbgdl}{$\textstyle\lVert{}$}}
					 {\raisebox{0.25\upbgdl}{$\scriptstyle\lVert{}$}}}
\newcommand*{\RVert}{
	\mathchoice{\raisebox{0.5\upbgdl}{$\displaystyle{\fit}\big\rVert$}}
					 {\raisebox{\upbgdl}{${\fit}\rVert$}} 
					 {\raisebox{0.5\upbgdl}{$\textstyle{\fit}\rVert$}}
					 {\raisebox{0.25\upbgdl}{$\scriptstyle{\fit}\rVert$}}}
\newcommand*{\NMV}[1]{\ensuremath\left.\fit[-0.08em]\LVert\VEC{#1}\RVert\fit[-0.08em]\right.}
\newcommand*{\NCP}[2]{\ensuremath\left.\fit[-0.08em]\LVert\VCP{#1}{#2}\RVert\fit[-0.08em]\right.}

\begin{document}

\mypart{\mytitle} 

\begin{abstract}
This is the first part of a series of four papers in The $\Pi\mrm{ME}$ Journal.
The focus of these papers is on a natural extension of Heron's 2000 year old formula for the area of a triangle to the volume of a tetrahedron.
While it is already quite remarkable that such an extension should have been overlooked for so long, that turns out to be just the first of a long, perhaps even infinite, sequence of surprises to be unveiled by the new formula.
This paper gives a high-level overview of these results, and then uses a novel fusion of classical vector algebra with Euclidean distance geometry to show that the areal vectors, or normal vectors scaled by the areas, of the four faces of a tetrahedron are intimately related to those of its three \TDEF{medial parallelograms}.
These relations justify the heterodox point-of-view that the medial parallelograms should be regarded as \EMPH{interior faces}.
Of particular note in subsequent parts of the series is the fact (as first shown by B.~D.~S. McConnell) that the areas of the interior and exterior faces together determine a non-degenerate tetrahedron up to isometry, and that they jointly satisfy a quadratic identity together with a system of 18 linear inequalities.
These inequalities, which extend the good old triangle inequality from lengths to areas, will henceforth be known as the \TDEF{tetrahedron inequalities}.
\end{abstract}

\section{Introduction} \label{sec:intro} 
Heron's formula for the squared area of a triangle is one of the oldest and most celebrated equations in classical Euclidean geometry \cite{Dunham:1990,Liberti:2016}.
It has been extended to $n$-dimensional simplices for all positive integers $n$ via Cayley-Menger determinants, which similarly give the $n$-simplices' squared hyper-volumes as homogeneous polynomials in their edge lengths \cite{Alfakih:2018, Blumenthal:1953,  Bowers:2017, Crippen:1988, Liberti:2017, Menger:1931}, but with one important difference:
For $n = 2$ the three-point determinant can be written as a product of four sums of the signed edge lengths, which is generally what is meant by \textquote{Heron's formula,} whereas for $n > 2$ the Cayley-Menger determinants are polynomials in the \EMPH{squared} edge lengths that do not factorize.
As a result, the combinatorial geometry of their zeros is far less transparent than it is with Heron's formula \cite{Dress:1986,Dress:1991,Easthope:1989}, where one can see at a glance that there are exactly three ways in which a triangle can have an area of zero, depending on which one of its vertices lies on the edge spanned by the other two.
Indeed three of the factors  in Heron's formula are simply the deviations of the three triangle inequalities among the edge lengths from saturation, meaning from holding as equalities, while the fourth is a non-degeneracy condition that vanishes if \& only if all three vertices coincide.

This series of four papers in The $\Pi\text{ME}$ Journal presents a rather different, but geometrically natural, extension of Heron's formula to tetrahedra.
This extension gives the fourth power of the volume as a polynomial in six simple rational functions of seven areal magnitudes that are canonically associated with each and every tetrahedron.
Four of these magnitudes are the areas of the usual four faces of the tetrahedron, while the remaining three are the areas of its \TDEF{medial parallelograms\/} (as defined in Section \myref{sec:vector} below).
Accordingly, the latter will be referred to herein as \textquote{interior faces.}
As will be shown in Section \myref{sec:arealgram}, the exterior and interior areas together determine a non-degenerate tetrahedron up to isometry.
The denominator of all six rational functions is just the exterior surface area of the tetrahedron, while each numerator factorizes into a product of two linear factors, one of which is a non-degeneracy condition and the other of which is the deviation from saturation of an areal generalization of the triangle inequality.

The significance of this extension lies not in providing yet-another means of calculating the volume of a tetrahedron \textit{per se\/}, but in the rather surprising nature of the geometric insights it yields into all the ways in which a tetrahedron can become \textquote{flat.}
Almost all of the formula's zeros, in fact, correspond to collinear tetrahedra with vertices at infinite distances from one another, although the ratios of those distances remain generically well defined.
The interpretation of these unconventional Euclidean configurations, and what they may have to tell us about the physical space in which we live, are questions of a kind generally seen as too obvious to even think about, and this paper ineluctably challenges that assumption.
Readers who doubt that such questions could be interesting are invited to consider the following innocent example:
\begin{displayquote}
\centering
\textit{How can the normal vectors of the usual four faces of a tetrahedron\\be coplanar but not collinear?}
\end{displayquote}

To make it subsequently clear that this extension is indeed geometrically natural, let us briefly revisit Heron's formula and its connection to the in-circle of a triangle $\GEO{ABC}$.
Hence let $a = \MAG{BC}$, $b = \MAG{AC}$, $c = \MAG{AB}$ be the lengths of the edges of $\GEO{ABC}$ opposite its vertices $\GEO{A}$, $\GEO{B}$, $\GEO{C}$, respectively, and let $s \coloneq \frac12(a + b + c)$ be its semi-perimeter.
Then the deviations of the three triangle inequalities from saturation are
\begin{equation}
u ~\coloneq~ \tfrac12 ( -a + b + c ) ~,\quad v ~\coloneq~ \tfrac12 (a - b + c) ~,\quad w ~\coloneq~ \tfrac12 (a + b - c) ~,
\end{equation}
where the factor of $1/2$ was introduced so that $a = v + w$, $b = u + w$, $c = u + v$.
These deviations have been called the \TDEF{Heron parameters\/} of a triangle \cite{Beardon:2015} (as well as \textquote{Gromov products} \cite{Cantarella:2019}), and clearly determine it up to isometry.
The Heron parameters are however not constrained by the triangle inequality, in that any $u, v, w \ge 0$ will yield distances that satisfy all three triangle inequalities among them.
Together with $s = u+v+w$, they also enable the squared area of the triangle to be expressed simply as
\begin{equation}
\MAG{ABC}^2 ~=~ s\, uvw ~=~ \tfrac12\, (u+v+w)\; \mrm{det\!} \begin{bmatrix} ~0&u&v~\\[2pt] ~u&0&w~\\[2pt] ~v&w&0~ \end{bmatrix} ~.
\end{equation}
Although this compact version of Heron's formula is well known, the product $uvw$ therein has not previously been viewed as a determinant.
Nevertheless, an analogous $4\times4$ determinant will be found in its extension to tetrahedra.

\begin{figure}
\input{Fig_I-1.tex}
\label{fig:heron}
\end{figure}

As illustrated in Fig.~\myref{fig:heron}, the Heron parameters are geometrically the distances from the vertices of the triangle to the \TDEF{in-touch points\/} at which its in-circle \textquote{touches} its edges.
They are also equal to the distances from the vertices to the \TDEF{ex-touch points\/} at which the triangle's ex-circles touch its edges, as well as the lines spanned by those edges.
Analogously, the aforementioned rational functions in our extension are the areas of the three triangles into which each exterior face of a tetrahedron is divided by its in-touch point.
There are twelve such triangles but, just as occurs with the Heron parameters of a triangle, these \textquote{contact triangles} will be found to occur in congruent pairs, giving rise to only six independent areas.
The \TDEF{natural parameters\/} of the tetrahedron will be defined as the values of these areas.
Additional parameters will be defined that are similarly related to the areas of the triangles into which the exterior faces are divided by their ex-touch points, and which are again rational functions of the seven facial areas.

Like the tetrahedron itself, all these parameters are uniquely determined by the areas of the exterior and interior faces together.
Expressing this geometric fact in algebraic terms will require us to take a bit of a detour through some very basic, though not very widely taught, vector geometry, to which we now turn.

\section{Areal relations from elementary vector algebra} \label{sec:vector}
\begin{figure}
\input{Fig_I-2.tex}
\label{fig:medoct}
\end{figure}
The nearly trivial relations among the inter-vertex vectors of a tetrahedron $\GEO{ABCD}$,
\begin{equation}
\VEC{AB} ~=~ \VEC{AC} \,\Vplus\, \VEC{CB} ~=~ \VEC{AC} \,\Vdiff\, \VEC{BC} ~=~ \VEC{CB} \,\Vdiff\, \VEC{CA}  ~=~ \VEC{DB} \,\Vdiff\, \VEC{DA} ~,  \label{eq:vecid1}
\end{equation}
are the basis for much of what follows.
An immediate consequence is that the cross product of the vectors between any two distinct pairs of vertices can be expanded in various ways, e.g.
\begin{equation} \begin{aligned}
\VCP{AB}{CD} ~=~ {}&{} \VEC{AB} \Xprod \big( \VEC{AD} \,\Vdiff\, \VEC{AC} \big) ~=~ \VCP{AB}{AD} \,\Vdiff\, \VCP{AB}{AC} \\
=~ {}&{} \VCP{AC}{CD} \,\Vdiff\, \VCP{BC}{CD} ~=~ \VCP{AC}{AD} \,\Vdiff\, \VCP{BC}{BD} ~.
\end{aligned} \label{eq:vecid2} \end{equation}
Up to sign, the cross products on the right are of course twice the areal vectors of the faces ($2$-faces, or facets) of the tetrahedron $\GEO{ABCD}$, which are just the (conventionally) outwards pointing normal vectors of those faces scaled by their areas, but what is the left-hand side?
It can be viewed as four times the cross product of the vector between the midpoints of the edges  $\GEO{AC}$ \& $\GEO{BC}$ and the vector between the midpoints of $\GEO{AC}$ \& $\GEO{AD}$:
\begin{equation}
\Big( \tfrac12 \big(\fit \GEO{B} \Vplus \GEO{C} \fit\big) \Vdiff \tfrac12 \big(\fit \GEO{A} \Vplus \GEO{C} \fit\big) \Big) \,\Xprod\, \Big( \tfrac12 \big(\fit \GEO{A} \Vplus \GEO{D} \fit\big) \Vdiff \tfrac12 \big(\fit \GEO{A} \Vplus \GEO{C} \fit\big) \Big) ~=~ \rfrac14\, \VCP{AB}{CD} ~.
\end{equation}
This is easily seen to be the same (up to sign) as the cross product of the vectors from the midpoint of any one of the four edges in the Hamiltonian cycle $\GEO{A} \mapsto \GEO{C} \mapsto \GEO{B} \mapsto \GEO{D} \mapsto \GEO{A}$ to the midpoints of the other two of those edges sharing a vertex with the first.
Thus \smash{$\VCP{AB}{CD}$} is four times an areal vector of the \TDEF{medial parallelogram\/} spanned by the midpoints of those four edges.
Similar interpretations also hold for the cross products \smash{$\VCP{AC}{BD}$} \& \smash{$\VCP{AD}{BC}$}.
This is further clarified and expanded upon in Fig.~\myref{fig:medoct}.

In the following, the areas of the (exterior) faces will be denoted by 
\begin{equation}
\MAG{ABC} ~=~ \rfrac12 \NCP{AB}{AC} ~=~ \rfrac12 \NCP{AB}{BC} ~=~ \rfrac12 \NCP{AC}{BC}
\end{equation}
etc., and the areas of the medial parallelograms (aka \TDEF{interior faces}) by
\begin{equation} \begin{aligned}
\MAG{AB|CD} ~=~ \rfrac14 \NCP{AB}{CD} ~,\quad \MAG{AC|BD} ~=~ \rfrac14 \NCP{AC}{BD} & \\
\text{and}\quad \MAG{AD|BC} ~=~ \rfrac14 \NCP{AD}{BC} & ~.
\end{aligned} \end{equation}
Then our first (new?) result is:
\begin{proposition} \label{thm:tetraineqs}
The areas of the interior and exterior faces of a tetrahedron $\GEO{ABCD}$ satisfy a system of $18$ linear inequalities, each of which involves one interior and two exterior faces.
These may logically be grouped into six triples, with two triples for each interior face, a typical example of which is:
\begin{subequations}
\label{eq:aineq0}
\begin{align}
4 \MAG{AB|CD} \,=\, \NCP{AB}{CD} \,\le\, {}&{} \NCP{AB}{AC} + \NCP{AB}{AD} \,=\, 2 \MAG{ABC} + 2 \MAG{ABD} \label{eq:aineq1} \\
2 \MAG{ABC} \,=\, \NCP{AB}{AC} \,\le\, {}&{} \NCP{AB}{AD} + \NCP{AB}{CD} \,=\, 2 \MAG{ABD} + 4 \MAG{AB|CD} \label{eq:aineq2} \\
2 \MAG{ABD} \,=\, \NCP{AB}{AD} \,\le\, {}&{} \NCP{AB}{AC}+ \NCP{AB}{CD} \,=\, 2 \MAG{ABC} + 4 \MAG{AB|CD} \label{eq:aineq3}
\end{align}
\end{subequations}
\end{proposition}
\begin{proof}
Equation (\myref{eq:aineq1}) follows immediately from the standard triangle inequality for vectors, $\mathbf v_1 = \mathbf v_2 \Vplus \mathbf v_3$ $\implies$ $\| \mathbf v_1 \| \le \| \mathbf v_2 \| + \| \mathbf v_3 \|$, applied to  the identity given by the second equality in Eq.~(\myref{eq:vecid2}), while Eqs.~(\myref{eq:aineq2}) \& (\myref{eq:aineq3}) follow from the two equations obtained by swapping terms between its left- \& right-hand sides.
The remaining five triples of inequalities are obtained simply by permuting the labels $\msf A, \msf B, \msf C$ \& $\msf D$.
\end{proof}	
\noindent Note that these are inequalities amongst the areas of the parallelograms spanned by the inter-vertex and inter-midpoint vectors, \underline{not} their lengths.
For this reason, although they are technically \textquote{triangle inequalities,} it seems more appropriate to call them \TDEF{tetrahedron inequalities}.
The deviations of the tetrahedron inequalities from saturation, expressed in terms of the facial areas, will henceforth be denoted by
\begin{subequations} \label{eq:taus}
\begin{align}
\Tau_{\!1}[\msf a, \msf b] ~\coloneq~ & 2 \MAG{abc} \,+\, 2 \MAG{abd} \,-\, 4 \MAG{ab|cd} ~,  \\
\Tau_{\!2}[\msf a, \msf b] ~\coloneq~ & 4 \MAG{ab|cd} \,+\, 2 \MAG{abd} \,-\, 2 \MAG{abc} ~, \\
\Tau_{\!3}[\msf a, \msf b] ~\coloneq~ & 4 \MAG{ab|cd} \,+\, 2 \MAG{abc} \,-\, 2 \MAG{abd} ~, \\
\intertext{and the corresponding non-degeneracy condition by}
\Tau_{\!0}[\msf a, \msf b] ~\coloneq~ & 2 \MAG{abc} \,+\, 2 \MAG{abd} \,+\, 4 \MAG{ab|cd} ~,
\end{align}
\end{subequations}
where $\{ \msf a, \msf b, \msf c, \msf d \} = \{ \msf A, \msf B, \msf C, \msf D \}$ and $\msf c < \msf d$ in alphabetic order.
In the Euclidean plane, any single $\Tau_{\!0}[\msf a, \msf b] = 0$ if \& only if $\GEO{a} = \GEO{b}$ is not on the line spanned by $\GEO{c} \ne \GEO{d}$.
\begin{remark}
Many additional, albeit weaker, linear inequalities among the seven facial areas can be derived by adding these deviations together, along with a great many more lower bounds on the areas analogous to the usual inverse triangle inequality.
An untypically well-known example is the upper bound on the area of any one exterior face given by the sum of the other three, e.g.
\begin{equation}
\MAG{ABC} ~\le~ \MAG{ABD} \,+\, \MAG{ACD} \,+\, \MAG{BCD} ~,
\end{equation}
along with the three others obtained by permuting the vertex labels \cite{Klamkin:1970, Richardson:1902}.
These four inequalities are known to be necessary and sufficient for the existence of a tetrahedron exhibiting the given exterior areas \cite{Hajji:2022, Mitrinovic:1989}.
\end{remark}

The following identity is usually attributed to Hermann Minkowski \cite{Minkowski:1903}.
\begin{lemma}[Minkowski's Identity]
The outwards (or inwards) pointing areal vectors of the exterior faces of a tetrahedron $\GEO{ABCD}$ (times $2$) satisfy
\begin{equation}
\VCP{AB}{AC} \,\Vdiff\, \VCP{AB}{AD} \,\Vplus\, \VCP{AC}{AD} \,\Vdiff\, \VCP{BC}{BD} ~=~ \mbf 0 ~. \label{eq:minkow}
\end{equation}
\end{lemma}
\begin{proof}
$\VEC{BC}\Xprod\VEC{BD} \,=\, \big( \VEC{BA} \Vplus \VEC{AC} \big) \Xprod \big( \VEC{BA} \Vplus \VEC{AD} \big) \,=\, \VEC{AB}\Xprod\VEC{AC} \,\Vdiff\, \VEC{AB}\Xprod\VEC{AD} \,\Vplus\, \VEC{AC}\Xprod\VEC{AD}$.
\end{proof}
\noindent This may be extended to the areal vectors of the seven faces together as follows.
\begin{proposition} \label{thm:arealvec}
The areal vectors of the exterior faces (times $4$) are equal to the following signed sums of the areal vectors of the interior faces (also times $4$):
\begin{subequations} \label{eq:havel0} \begin{align}
\VCP{AB}{CD} \,\Vplus\, \VCP{AC}{BD} \,\Vplus\, \VCP{AD}{BC} ~=~ & 2\, \VCP{AB}{AD} \label{eq:havel1} \\[2pt]
\Vdiff\, \VCP{AB}{CD} \,\Vplus\, \VCP{AC}{BD} \,\Vplus\, \VCP{AD}{BC} ~=~ & 2\, \VCP{AB}{AC} \label{eq:havel2} \\[2pt]
\VCP{AB}{CD} \,\Vplus\, \VCP{AC}{BD} \,\Vdiff\, \VCP{AD}{BC} ~=~ & 2\, \VCP{AC}{AD} \label{eq:havel3} \\[2pt]
\Vdiff\, \VCP{AB}{CD} \,\Vplus\, \VCP{AC}{BD} \,\Vdiff\, \VCP{AD}{BC} ~=~ & 2\, \VCP{BC}{BD} \label{eq:havel4}
\end{align} \end{subequations}
\end{proposition}
\begin{proof}
One can prove Eq.~(\myref{eq:havel1}) simply as follows:
\begin{multline*}
\VCP{AB}{CD} \,\Vplus\, \VCP{AC}{BD} \,\Vplus\, \VCP{AD}{BC} ~= \\[2pt]
\VEC{AB} \Xprod \big( \VEC{AD} \Vdiff \VEC{AC} \big) \,\Vplus\, \VEC{AC} \Xprod \big( \VEC{AD} \Vdiff \VEC{AB} \big) \,\Vplus\, \VEC{AD} \Xprod \big( \VEC{AC} \Vdiff \VEC{AB} \big) ~=~ 2\, \VCP{AB}{AD}
\end{multline*}
The proofs of the remaining identities are similar save for Eq.~(\myref{eq:havel4}), where Minkowski's identity (\myref{eq:minkow}) is also needed.
\end{proof}
\begin{remark} \label{rem:int-face-obs}
Applying the triangle inequality for vectors to these relations shows that the area of each exterior face is bounded above by the sum of the interior areas.
They also show that the areal vectors of the three interior faces (however oriented) determine those of the exterior faces which, by another well-known theorem of Minkowski \cite{Minkowski:1897}, determine the tetrahedron uniquely up to translation.
Finally, they show that the tetrahedron is equi-facial (also termed \textquote{equi-areal} \cite{McMullen:2000} or, in the older literature \cite{Court:1935}, \textquote{isosceles}) if \& only if the areal vectors of the interior faces are mutually orthogonal.
\end{remark}

We now turn to the trigonometric relations among the areal vectors.
While these formulae can only be ascribed to folklore \cite{Crane:2020,Koranyi:2006}, they are not given explicitly in other\-wise comprehensive surveys of tetrahedral geometry from the early 20$^\text{th}$ century \cite{Court:1935, Richardson:1902}.
\begin{lemma}[The Areal Law of Cosines] \label{thm:aloc}
Given a tetrahedron $\GEO{ABCD}$, the areal vectors of its interior and exterior faces satisfy
\begin{align} \label{eq:cosines}
\fit[-0.5em] \big( \VCP{AB}{AC} \big) \Dprod \big( \VCP{AB}{AD} \big) ~=~ {}&{} \NCP{AB}{AC} \NCP{AB}{AD} \cos( \varphi_\msf{AB} ) \nonumber \\[2pt]
=~ {}&{} \sfrac12\, \Big( \NCP{AB}{AC}^2 +\, \NCP{AB}{AD}^2 -\,  \NCP{AB}{CD}^2 \Big) \\[2pt]
=~ {}&{} \tfrac14\, \big( \Tau_0[\msf A, \msf B]\, \Tau_1[\msf A, \msf B] \,-\, \Tau_2[\msf A, \msf B]\, \Tau_3[\msf A, \msf B]\, \big) , \nonumber
\end{align}
where \textquote{$\,\Dprod$} is the vector dot product and $\varphi_\msf{AB}$ is the internal dihedral angle between $\GEO{ABC}$\linebreak[2] \& $\GEO{ABD}$, along with the five analogous relations obtained by permuting vertex labels.
\end{lemma}
\begin{proof}
The first line of Eq.~(\myref{eq:cosines}) is just the geometric definition of the dot product of the cross products and of the dihedral angle itself.
The second line is obtained by dotting each side of the second equality in Eq.~(\myref{eq:vecid2}) with itself, i.e.
\begin{equation*}
\NCP{AB}{CD}^2 ~=~ \NCP{AB}{AC}^2 +\, \NCP{AB}{AD}^2 -\, 2\, \big( \VCP{AB}{AC} \big) \Dprod \big( \VCP{AB}{AD} \big) ~,
\end{equation*}
followed by rearrangement.
The third line follows from the definitions in Eq.~(\myref{eq:taus}).
\end{proof}
\begin{lemma}[The Areal Law of Sines] \label{thm:alos}
Given a tetrahedron $\GEO{ABCD}$, its edge lengths, the areas of its exterior faces, and its volume $\MAG{ABCD} = \raisebox{0.5\upbgdl}{$\big|$}\fit \VEC{AB} \Dprod \big( \VCP{AC}{AD} \big) \fit\raisebox{0.5\upbgdl}{$\big|$} \!\bigm/\! 6$ satisfy
\begin{align} \label{eq:sines}
\fit[-0.5em] \NMV{AB} \raisebox{\upbgdl}{\fnsize$\Big|$}\fit \VEC{AB} \Dprod \big( \VCP{AC}{AD} \big) \raisebox{\upbgdl}{\fnsize$\Big|$} ~ {}&{} =\, \NCP{AB}{AC} \NCP{AB}{AD} \sin( \varphi_\msf{AB} ) \\
{}&{}=\, \sqrt{\, \NCP{AB}{AC}^2 \NCP{AB}{AD}^2 -\, \Big(\fit[-0.05em] \big( \VCP{AB}{AC} \big) \Dprod \big( \VCP{AB}{AD} \big) \fit[-0.05em]\Big)^{\!2} } \nonumber \\[2pt]
{}&{}=~ \tfrac12 \sqrt{\, \Tau_0[\msf A, \msf B]\, \Tau_1[\msf A, \msf B]\, \Tau_2[\msf A, \msf B]\, \Tau_3[\msf A, \msf B]\, } ~, \nonumber
\end{align}
where $\varphi_\msf{AB}$ is the dihedral angle as above, along with the five other relations obtained by permuting the vertex labels.
\end{lemma}
\begin{proof}
The standard vector algebra identity $(\mbf p \Xprod \mbf q ) \Xprod ( \mbf p \Xprod \mbf r) = \big( \mbf p \Dprod (\mbf q \Xprod \mbf r) \big)\, \mbf p$ implies
\begin{equation}
\raisebox{\upbgdl}{\fnsize$\Big\|$} \big( \VCP{AB}{AC} \big) \Xprod \big( \VCP{AB}{AD} \big) \raisebox{\upbgdl}{\fnsize$\Big\|$}^{\,2} =~ \NMV{AB}^2 \! \Big( \VEC{AB} \Dprod \big( \VCP{AC}{AD} \big) \Big)^{\!2} ~. \label{eq:sines0}
\end{equation}
Together with $\| \mbf p \Xprod \mbf q\fit \| = \| \mbf p\fit \| \| \mbf q\fit \| \sin(${\large$\angle$}$(\mbf p, \mbf q) )$, this gives the first line of Eq.~(\myref{eq:sines}).
The second line follows from Lagrange's identity $\| \mbf p \Xprod \mbf q \|^2 \fit=\fit \| \mbf p \|^2 \, \| \mbf q \|^2 \,-\, (\mbf p \fit\Dprod\fit \mbf q)^2$ applied to the left-hand side of Eq.~(\myref{eq:sines0}).
The last line of Eq.~(\myref{eq:sines}) then follows from the last line of Eq.~(\myref{eq:cosines}) together with the easily proven complementary identity $4 \NCP{AB}{AC} \! \NCP{AB}{AD} =\, \Tau_0[\msf A, \msf B]\, \Tau_1[\msf A, \msf B] \,+\, \Tau_2[\msf A, \msf B]\, \Tau_3[\msf A, \msf B]$.
\end{proof}
\noindent This areal law of sines can readily be shown to imply the spherical law of sines.

\begin{remark}
It is also possible to derive a law of cosines for the dot product of the areal vectors of an interior and an exterior face, e.g.
\begin{equation} \label{eq:cosine2}
\big( \VCP{AB}{AC} \big) \Dprod \big( \VCP{AB}{CD} \big) ~=~ \sfrac12\, \Big( \NCP{AB}{AD}^2 -\, \NCP{AB}{CD}^2 -\,  \NCP{AB}{AC}^2 \Big) ~,
\end{equation}
and for the areal vectors of two interior faces, e.g.
\begin{equation} \label{eq:cosine3} \begin{aligned}
\big( \VCP{AB}{CD} \big) \Dprod \big( \VCP{AC}{BD} \big) ~= & \\
\sfrac12\, \Big( \NCP{AB}{AD}^2 +\, {}&{} \NCP{AC}{AD}^2 -\, \NCP{AB}{AC}^2 -\, \NCP{BC}{BD}^2 \Big) ~.
\end{aligned} \end{equation}
The extant formula for the volume as $2/3$ the area of an interior face $\MAG{AB \tmv CD}$ times the perpendicular distance between the lines of $\GEO{AB}$ \& $\GEO{CD}$ (see e.g.~Ex.~12 on pg.~91 of Ref.~\cite{Court:1935}) can also be seen as a kind of areal law of sines, as can the lovely formula
\begin{equation}
\big( \VCP{AB}{CD} \big) \Dprod \Big( \big( \VCP{AC}{BD} \big) \Xprod \big( \VCP{AD}{BC} \big) \Big) ~=~ 2\, \Big( \VEC{AB} \fit\Dprod \big( \VCP{AC}{AD} \big) \Big)^{\!2} ~.
\end{equation}
\end{remark}

Another formula which has also been called the law of cosines for a tetrahedron \cite{Hampton:2014, Murray-Lasso:2002} (and hence our addition of the qualifier \textquote{areal} above) is:
\begin{lemma}
Given a tetrahedron $\GEO{ABCD}$, the areas of and dihedral angles between its exterior faces satisfy:
\begin{align} \label{eq:loctet}
&{} \fit[-0.5em] \NCP{BC}{BD}^2 \,=~ \NCP{AB}{AC}^2 +\, \NCP{AB}{AD}^2 +\, \NCP{AC}{AD}^2 \\[2pt]
&{} \nonumber \fit[2em]
\begin{aligned} -\; 2 \NCP{AB}{AC}\! \NCP{AB}{AD} \cos\!\big( \varphi_\msf{AB} \big) -\, 2\, \NCP{AB}{AC}\! \NCP{AC}{AD} \cos\!\big( \varphi_\msf{AC} \big) & \\[2pt]
-\: 2 \NCP{AB}{AD}\! \NCP{AC}{AD} \cos\!\big( \varphi_\msf{AD} \big) & \end{aligned}
\end{align}
\end{lemma}
\begin{proof}
Simply solve Eq.~(\myref{eq:minkow}) for $\VCP{BC}{BD}$, then dot each side with itself and apply Lemma \myref{thm:aloc} taking account of $\big( \VCP{AB}{AC} \big) \Dprod \big( \VCP{AC}{AD} \big) = -\, \big( \VCP{AC}{AB} \big) \Dprod \big( \VCP{AC}{AD} \big)$.
\end{proof}
This leads to the algebraic identity that connects the interior and exterior areas.
\begin{proposition}[Yetter's Identity] \label{thm:yid}
Given a tetrahedron $\GEO{ABCD}$, the areas of its interior and exterior faces satisfy
\begin{align} &
\begin{aligned} \NCP{AB}{AC}^2 +\, \NCP{AB}{AD}^2 +\, \NCP{AC}{AD}^2 +\, \NCP{BC}{BD}^2 & \\[2pt]
=~ \NCP{AB}{CD}^2 +\, \NCP{AC}{BD}^2 +\, \NCP{AD}{BC}^2 & \end{aligned} \label{eq:yetter} \\[4pt]
\Longleftrightarrow\quad & \Xi\Big( 2\! \MAG{ABC}\!,\fit 2\! \MAG{ABD}\!,\fit 2\! \MAG{ACD}\!,\fit 2\! \MAG{BCD}\!,\fit 4\! \MAG{AB|CD}\!,\fit 4\! \MAG{AC|BD}\!,\fit 4\! \MAG{AD|BC} \!\Big) =~ 0 ~, \nonumber
\end{align}
where the quadratic polynomial $\Xi(a,b,c,d,e,f,g) \coloneq a^2 + b^2 + c^2 + d^2 - e^2 - f^2 - g^2$\linebreak[2] may be abbreviated as \textquote{$\Xi$} when the indeterminates representing twice the exterior and four times the interior facial areas are clear from context.
\end{proposition}
\begin{proof}
This follows from the second equality of Eq.~(\myref{eq:cosines}) in Lemma \myref{thm:aloc} applied to Eq.~(\myref{eq:loctet}).
Alternatively, the first equality of Eq.~(\myref{eq:cosines}) shows that Eq.~(\myref{eq:loctet}) is equivalent to
\begin{multline*}
\raisebox{0.5\upbgdl}{$\big\|\fit$} \VCP{BC}{BD} \raisebox{0.5\upbgdl}{$\fit\big\|$}^2 ~=~ 
\raisebox{0.5\upbgdl}{$\big\|\fit$} \VCP{AB}{AC} \,\Vdiff\, \VCP{AB}{AD} \raisebox{0.5\upbgdl}{$\fit\big\|$}^2 +\, \raisebox{0.5\upbgdl}{$\big\|\fit$} \VCP{AB}{AC} \,\Vplus\, \VCP{AC}{AD} \raisebox{0.5\upbgdl}{$\fit\big\|$}^2 \\[1pt] 
+\, \raisebox{0.5\upbgdl}{$\big\|\fit$} \VCP{AB}{AD} \,\Vdiff\, \VCP{AC}{AD} \raisebox{0.5\upbgdl}{$\fit\big\|$}^2 -\, \NCP{AB}{AC}^2 -\, \NCP{AB}{AD}^2 -\, \NCP{AC}{AD}^2 .
\end{multline*}
But by Eq.~(\myref{eq:vecid2}) and its permutations, it is easily seen that the sum and differences of the cross products inside the norms in this equation are four times the areal vectors of the interior faces, whence Eq.~(\myref{eq:yetter}) follows.
\end{proof}
\noindent This identity was given as Ex.~17 on pg.~294 of Altshiller-Court's 1935 text \cite{Court:1935}.
More recently, it has been extended by David N.~Yetter to a family of identities connecting the \textquote{hyper-areas} of the facets and medial sections of $n$-simplices for all $n \ge 3$ \cite{Yetter:2010} (hence its attri\-bution to him), although only the $n = 3$ case above will be used in what follows.

\section{The exterior areal Gram matrices} \label{sec:arealgram}
The (exterior) \TDEF{areal Gram matrix\/} at any vertex of a tetrahedron $\GEO{ABCD}$, say $\GEO{A}$, plays a central role in what follows (extensions to $n$-dimensional spaces of constant curvature may be found in Refs.~\cite{Abrosimov:2014,Kokkendorff:2006}).
This is the $3\times3$ symmetric matrix $\mbf G_\msf A$ consisting of the dot products of twice the outwards-pointing areal vectors of the three exterior faces meeting at $\GEO{A}$, and as such is positive definite for any non-degenerate tetrahedron.
Using the areal law of cosines (Eq.~(\myref{eq:cosines})), it may be expressed as a matrix of linear polynomials in indeterminates representing
the \EMPH{squared} facial areas $F_\msf{ABC} \leftrightarrow 4 \smash{\MAG{ABC}}^2$, $\ldots\,$, $F_\msf{AD|BC} \leftrightarrow 16 \smash{\MAG{AD\tmv BC}}^2$, namely:
\begin{subequations} \label{eq:gramatA} \begin{align}
\fit[-0.5em] \mbf G_F[\msf A] ~\coloneq~ \left[ \begin{smallmatrix} \hit[1.5ex]
F_\msf{ABC} & \frac12 \big( F_\msf{AB|CD} \,-\, F_\msf{ABC} \,-\, F_\msf{ABD} \big) & \frac12 \big( F_\msf{AC|BD} \,-\, F_\msf{ABC} \,-\, F_\msf{ACD} \big) \\[1pt]
\frac12 \big( F_\msf{AB|CD} \,-\, F_\msf{ABC} \,-\, F_\msf{ABD} \big) & F_\msf{ABD} & \frac12 \big( F_\msf{AD|BC} \,-\, F_\msf{ABD} \,-\, F_\msf{ACD} \big) \\[1pt]
\frac12 \big( F_\msf{AC|BD} \,-\, F_\msf{ABC} \,-\, F_\msf{ACD} \big) & \frac12 \big( F_\msf{AD|BC} \,-\, F_\msf{ABD} \,-\, F_\msf{ACD} \big) & F_\msf{ACD} \dit[1ex] \end{smallmatrix} \right] & \label{eq:gramatAa} \\ 
\fit[-0.5em] \longleftrightarrow\quad \mbf G_\msf A ~\coloneq~ \left[ \begin{smallmatrix} \hit[1.5ex]
{\LVert \vv{\msf{AB}} \fit\Xprod\fit \vv{\msf{AC}} \RVert}^2 & -\, \big( \vv{\msf{AB}} \fit\Xprod\fit \vv{\msf{AC}} \big) \,\Dprod\, \big( \vv{\msf{AB}} \fit\Xprod\fit \vv{\msf{AD}} \big) & -\,\big( \vv{\msf{AC}} \fit\Xprod\fit \vv{\msf{AB}} \big) \,\Dprod\, \big( \vv{\msf{AC}} \fit\Xprod\fit \vv{\msf{AD}} \big) \\[2pt]
-\, \big( \vv{\msf{AB}} \fit\Xprod\fit \vv{\msf{AC}} \big) \,\Dprod\, \big( \vv{\msf{AB}} \fit\Xprod\fit \vv{\msf{AD}} \big) & {\LVert \vv{\msf{AB}} \fit\Xprod\fit \vv{\msf{AD}} \RVert}^2 & -\, \big( \vv{\msf{AB}} \fit\Xprod\fit \vv{\msf{AD}} \big) \,\Dprod\, \big( \vv{\msf{AC}} \fit\Xprod\fit \vv{\msf{AD}} \big) \\[2pt]
-\, \big( \vv{\msf{AC}} \fit\Xprod\fit \vv{\msf{AB}} \big) \,\Dprod\, \big( \vv{\msf{AC}} \fit\Xprod\fit \vv{\msf{AD}} \big) & -\, \big( \vv{\msf{AB}} \fit\Xprod\fit \vv{\msf{AD}} \big) \,\Dprod\, \big( \vv{\msf{AC}} \fit\Xprod\fit \vv{\msf{AD}} \big) & {\LVert \vv{\msf{AC}} \fit\Xprod\fit \vv{\msf{AD}} \RVert}^2 \dit[1ex] \end{smallmatrix} \right] & \label{eq:gramatAb}
\end{align} \end{subequations}
Note that the negative signs before the dot products in the entries adjacent to the diagonal are due to the way these cross products are signed in Minkowski's identity (\myref{eq:minkow}), while the negative sign in the corner entries is due to the swap of the vectors in the cross product of the dot product's first factor that is needed to apply Eq.~(\myref{eq:cosines}) directly.

The determinant of this matrix will be denoted by $\Gamma_{\!F} [\msf A] \coloneq \mrm{det}\big( \mbf G_F[\msf A] \big)$ $\leftrightarrow$ $\mrm{det} \big( \mbf G_\msf A \big)$,
and referred to as the \TDEF{Gramian at\/} \smash{$\GEO{A}$}. 
Expansion shows that $\Gamma_{\!F} [\msf A]$ is a homogeneous cubic polynomial containing $30$ terms, as are the Gramians at the other three vertices $\Gamma_{\!F}[\msf B]$, $\Gamma_{\!F}[\msf C]$ \& $\Gamma_{\!F}[\msf D]$.
They are related as follows.
\begin{lemma} \label{thm:gramequiv}
As polynomials in the indeterminates $F_\msf{ABC}$, $\ldots\,$, $F_\msf{AD|BC}$ $\in \mbb R$, the four Gramians satisfy
\begin{equation}
\Gamma_{\!F}[\msf A] ~\equiv~ \Gamma_{\!F}[\msf B] ~\equiv~ \Gamma_{\!F}[\msf C] ~\equiv~ \Gamma_{\!F}[\msf D] \mod \breve\Xi_F ~,
\end{equation}
where $\breve\Xi_F$ is the linear polynomial in the indeterminates $F$ corresponding to the quadratic polynomial $\Xi$ of Eq.~(\myref{eq:yetter}).
Thus if the indeterminates satisfy Yetter's identity $\breve\Xi_{F\!} = 0$, the Gramians are all equal.\hspace*{-2pt} 
\end{lemma}
\begin{proof}
Using computer algebra, it is easily shown that
\begin{multline*}
\Gamma_{\!F}[\msf B] \,-\ \Gamma_{\!F}[\msf A] ~=~ \sfrac14\, \Big( \big( F_\msf{ACD} - F_\msf{BCD} \big) \big( 2\, F_\msf{ABC} + 2\, F_\msf{ABD} - F_\msf{AB|CD} \big) \\[-4pt]
+\, \big( F_\msf{ABC} - F_\msf{ABD} \big) \big( F_\msf{AC|BD} - F_\msf{AD|BC} \big) \Big)\, \breve\Xi_F \,.
\end{multline*}
Similar results are obtained for $\Gamma_{\!F}[\msf C] \,-\ \Gamma_{\!F}[\msf A]$ and $\Gamma_{\!F}[\msf D] \,-\ \Gamma_{\!F}[\msf A]$.
\end{proof}
\noindent The $2\times2$ principal minors of $\mbf G_\msf A$ are also of interest, and will be denoted by
\begin{equation}
\Gamma_{\!F}[\msf A; \msf B] ~\coloneq~ \mrm{det}\! \left[ \begin{smallmatrix} \hit[1.5ex] F_\msf{ABC} & \frac12\, \big( F_\msf{AB|CD} \,-\, F_\msf{ABC} \,-\, F_\msf{ABD} \big) \\[1pt] \frac12\, \big( F_\msf{AB|CD} \,-\, F_\msf{ABC} \,-\, F_\msf{ABD} \big) & F_\msf{ABD} \dit[1ex] \end{smallmatrix} \right] ~=~ \Gamma_{\!F}[\msf B; \msf A] ~,
\end{equation}
with analogous definitions for $\Gamma_{\!F}[\msf A; \msf C] = \Gamma_{\!F}[\msf C; \msf A]$, $\Gamma_{\!F}[\msf A; \msf D] = \Gamma_{\!F}[\msf D; \msf A]$ and the other $2\times2$ principal minors of the four Gram matrices.
\begin{lemma} \label{thm:2x2minors}
Given any indeterminates $f_\msf{ABC}$, $\ldots\,$, $f_\msf{AD|BC}$ $\in \mbb R$, and letting $F_\msf{abc} \coloneq f_\msf{abc}^2$ \& $F_\msf{ab|cd} \coloneq f_\msf{ab|cd}^2\,$ for all $\{ \msf a,\, \ldots,\, \msf d \} = \{ \msf A,\, \ldots,\, \msf D \}$, we have
\begin{equation} \label{eq:2x2minors}
\Gamma_{\!F}[\msf a; \msf b] ~=~ \tfrac14\, \Tau_{\!0:f}[\msf a,\msf b]\,\Tau_{\!1:f}[\msf a,\msf b]\,\Tau_{\!2:f}[\msf a,\msf b]\,\Tau_{\!3:f}[\msf a,\msf b] ~,
\end{equation}
where $\Tau_{\!1:f}[\msf a, \msf b] \coloneq f_\msf{abc} + f_\msf{abd} - f_\msf{ab|cd}$, $\Tau_{\!2:f}[\msf a, \msf b] \coloneq f_\msf{ab|cd} + f_\msf{abd} - f_\msf{abc}$, $\Tau_{\!3:f}[\msf a, \msf b] \coloneq f_\msf{ab|cd} + f_\msf{abc} - f_\msf{abd}$ ($\msf c < \msf d$) are the linear polynomials corresponding to the deviations of tetrahedron inequalities from saturation in Eq.~(\myref{eq:taus}), and $\Tau_{\!0:f}[\msf a, \msf b] \coloneq f_\msf{abc} + f_\msf{abd} + f_\msf{ab|cd}$ are those of the associated non\-degeneracy factors.
\end{lemma}
\begin{proof}
For Euclidean areas this follows directly from Eq.~(\myref{eq:sines}).
For general indeterminates one need only replace each $F$ in the Gram matrices (as in Eq.~(\myref{eq:gramatA})) by $f^2$, then factorize the $2\times2$ principal minors to get the right-hand side of Eq.~(\myref{eq:2x2minors}).
\end{proof}
\noindent Observe that the four factors in these formulae are analogous to those in Heron's formula,\pagebreak[2] so these minors can likewise be written as three-point Cayley-Menger determinants, albeit in indeterminates representing one interior and two exterior areas rather than inter-vertex distances (as in Eq.~(\myref{eq:cmd3}) below).

In order to see what these polynomials are geometrically, recall that the squared areas of the exterior faces of a tetrahedron $\GEO{ABCD}$ may be written as three-point Cayley-Menger determinants in its squared edge lengths $D_\msf{AB} \leftrightarrow \smash{\MAG{AB}}^2$ etc.~\cite{Blumenthal:1953,Crippen:1988}, e.g.
\begin{equation} \label{eq:cmd3}
4 \MAG{ABC}^2 \,\longleftrightarrow~ \CMD_D[\msf A, \msf B, \msf C] \,\coloneq -\tfrac14\, \mrm{det}\big( \bar{\mbf D}[\msf A, \msf B, \msf C] \big) \,\coloneq -\raisebox{-1pt}{$\sfrac14$}\, \mrm{det}\! \left[ \begin{smallmatrix} \hit[2ex] 0~&1&1&1 \\[3pt] 1~&0&D_\msf{AB}&D_\msf{AC} \\[2pt] 1~&D_\msf{AB}&0&D_\msf{BC} \\[2pt] 1~&D_\msf{AC}&D_\msf{BC}&0 \dit[1ex] \end{smallmatrix} \right] ,
\end{equation}
where $\bar{\mbf D}$ denotes the bordered matrix of squared distances $\mbf D \coloneq [D_\msf{ab}]$ on the right-hand side of this equation.
Via Lagrange's identity and the usual law of cosines, the squared areas of the interior faces may likewise be written as polynomials in the squared edge lengths \cite{Yetter:1998} or, using Cayley-Menger determinants, as e.g.
\begin{equation} \label{eq:talata}
\fit[-0.5em] 16 \MAG{AB \tmv CD}^2 \,\longleftrightarrow~ \CMD_D[\msf A, \msf B]\, \CMD_D[\msf C, \msf D] \,-\, {\CMD_D[\msf A, \msf B; \msf C, \msf D]}^2 ~\eqcolon~ \CMD_D[ \msf A, \msf B \,\tmv\, \msf C, \msf D ] ~,
\end{equation}
wherein the non-symmetric two-point Cayley-Menger determinant is defined as
\begin{equation}
\! \CMD_D[\msf A, \msf B; \msf C, \msf D] \,\coloneq\, \sfrac12\fit \mrm{det}\! \left[ \begin{smallmatrix} \hit[1.5ex] 0~&1&1 \\[3pt] 1~&D_\msf{AC}&D_\msf{AD} \\[2pt] 1~&D_\msf{BC}&D_\msf{BD} \dit[1ex] \end{smallmatrix} \right] =~ \tfrac12 ( D_\msf{AD} + D_\msf{BC} - D_\msf{AC} - D_\msf{BD} ) \,\leftrightarrow~ \VEC{AB} \,\Dprod\, \VEC{CD} \fit, \label{eq:talata2}
\end{equation}
while $\CMD_D[\msf A, \msf B] \coloneq \CMD_D[\msf A, \msf B; \msf A, \msf B] = D_\msf{AB}$ and similarly $\CMD_D[\msf C, \msf D] = D_\msf{CD}$.
Because they are a special case of a determinantal formula which expresses the squared \textquote{hyper-areas} of the medial sections of $n$-simplices in terms of their edge lengths that was discovered by Istv\'an Talata \cite{Talata:2003}, we shall call the poly\-nomials $\CMD_D[\msf A, \msf B \,\tmv\, \msf C, \msf D]$ \TDEF{Talata determinants}.

These relations allow us to convert polynomials in the squared facial areas into polynomials in the squared edge lengths by simple substitution.
\begin{proposition}
Given a Euclidean tetrahedron $\GEO{ABCD}$, the Gramians $\Gamma_{\!F}[\msf a]$ with $\msf a \in \{ \msf A, \msf B, \msf C, \msf D \}$, when evaluated at $F_\msf{ABC} = 4 \smash{\MAG{ABC}}^2,\, \ldots, F_\msf{AD|BC} = 16 \smash{\MAG{AD \tmv BC}}^2$, are all equal to the fourth power of $3!$ times its volume, henceforth $t^4 \coloneq \smash{\raisebox{\upbgdl}{$\big(\fit$} 6 \MAG{ABCD} \raisebox{\upbgdl}{$\!\big)$}}{\hit[1.8ex]}^{\smash4}$.
The $2\times2$ principal minors of the Gram matrices, $\Gamma_{\!F}[\msf a; \msf b]$ with $\msf a, \msf b \in \{ \msf A, \msf B, \msf C, \msf D \}$ and $\msf a \ne \msf b$, likewise evaluated at these multiples of the squared areas in $\GEO{ABCD}$, are equal to $\smash{\MAG{ab}}^2\, t^2$.
\end{proposition}
\begin{proof}
On substituting for the six squared areas in $\Gamma_{\!F}[\msf A]$ using Eqs.~such as (\myref{eq:cmd3}), (\myref{eq:talata}) \& (\myref{eq:talata2}), one obtains (preferably with the aid of computer algebra) the square of the four-point Cayley-Menger determinant $\CMD_D[\msf A, \msf B, \msf C, \msf D] \coloneq \mrm{det}\big( \bar{\mbf D}[\msf A, \msf B, \msf C, \msf D] \big) / 8$, or
\begin{equation} \label{eq:gamma2d}
\Gamma_{\!F}[\msf A] \,\big|_{F=\Delta_\text{\tiny$\scriptstyle D$}} \,=~  \CMD_D[\msf A, \msf B, \msf C, \msf D]^2 ~\longleftrightarrow ~ t^4 ~,
\end{equation}
and likewise for the Gramians at the other three vertices (a considerably more complicated proof of the $n$-dimensional version may be found in Ref.~\cite{Veljan:1995}).
Similarly, on substituting for the squared areas in the $2\times2$ principal minor $\Gamma_{\!F}[\msf A; \msf B]$, one obtains
\begin{equation} \label{eq:minor2d}
\Gamma_{\!F}[\msf A; \msf B] \,\big|_{F=\Delta_\text{\tiny$\scriptstyle D$}} \,=~  \CMD_D[\msf A, \msf B]\, \CMD_D[\msf A, \msf B, \msf C, \msf D] ~\longleftrightarrow\, \MAG{AB}^2 t^2 ~,
\end{equation}
with analogous results for the other $2\times2$ principal minors of  $\mbf G_F[\msf A]$ as well as those of the Gram matrices at the remaining three vertices.
\end{proof}
\noindent Equations (\myref{eq:2x2minors}) \& (\myref{eq:minor2d}) of course constitute a polynomial version of the areal law of sines (\myref{eq:sines}).
It is also easily shown that the Cayley-Menger \& Talata determinants them\-selves satisfy Yetter's identity, meaning that upon expansion 
\begin{equation*}
\breve\Xi\big(\, \CMD_D[\msf A, \msf B, \msf C]\,,\,\ldots\,,\, \CMD_D[\msf A, \msf D \,|\, \msf B, \msf C] \,\big) =\, 0 \qquad\text{(the $0$ polynomial).}
\end{equation*}
\begin{remark}
The example $D_\msf{AB} = D_\msf{AC} = D_\msf{BC} = 12$, $D_\msf{AD} = D_\msf{BD} = 4$ \& $D_\msf{CD} = 3$ shows that all the three-point Cayley-Menger \& Talata determinants can be positive while the four-point Cayley-Menger determinant $\CMD_D[\msf A, \msf B, \msf C, \msf D] = -39$ is negative.
Nevertheless, on applying Eqs.~(\myref{eq:2x2minors}) \& (\myref{eq:minor2d}) to each of the six distinct minors $ \Gamma_{\!F}[\msf a; \msf b]$ ($\msf a, \msf b \in \{ \msf A, \msf B, \msf C, \msf D \}$, $\msf a \ne \msf b$), one finds easily that the areas in any such metric space necessarily violate a tetrahedron inequality in every one of the six triples thereof.
This shows that any metric space wherein the Talata determinants in every quadruple are all non-negative and the corresponding areas satisfy all the tetrahedron inequalities also fulfils the Euclidean four-point property \cite[Def.~50.1]{Blumenthal:1953} (because it saturates if \& only if $\GEO{AB} \parallel \GEO{CD}$, the relation $\CMD_D[\msf A, \msf B \,|\, \msf C, \msf D] \ge 0$ among the distances will be called the \TDEF{trapezoid inequality}).
It is well-known that even when all the five-point Cayley-Menger determinants vanish, such a metric space need not be realizable in three-dimensional Euclidean space (cf.~Refs.~\cite[Sec.~44]{Blumenthal:1953} or \cite{Bowers:2017}).
\end{remark}

\smallskip 
By first computing the volume $t$ via Eq.~(\myref{eq:gamma2d}) and then the edge lengths $\MAG{ab}$ from Eq.~(\myref{eq:minor2d}) and its permutations ($\msf a, \msf b \in \{ \msf A, \msf B, \msf C, \msf D \}$, $\msf a \ne \msf b$), one obtains a simple proof that the areas of the seven faces of a non-degenerate Euclidean tetrahedron determine it up to isometry.
This proof was first given, to this author's knowledge, in an unpublished paper posted on a remarkable online Blog by an amateur but dedicated geometer named Billy Don Sterling McConnell, apparently around 2012 (at the time of writing, this Blog was accessible at \url{http://hedronometry.com}).
McConnell also noted that this calculation would succeed only if the Gramians were strictly positive and the $2\times2$ principal minors non-negative (and hence likewise strictly positive), i.e.~the Gram matrices were all positive definite.
Crane \& Yetter subsequently also derived the edge lengths from the areas using spherical trigonometry \cite{Crane:2020}, but did not carefully identify the conditions that the putative areas must satisfy in order for their calculation to succeed.

Once the edge lengths, however obtained, are available coordinates for the vertices can be computed by standard \textquote{multi-dimensional scaling} techniques based on the \textquote{lineal} Gram matrices of dot products among the vectors along the edges originating at any vertex (see e.g.~Refs.~\cite{Borg:2005,Crippen:1978,Gower:2004}).
The proof given here instead computes vertex coordinates which reproduce the given areas directly from the areas themselves, without explicitly determining the edge lengths first.
In essence, it exploits the fact that the areal Gram matrix at $\GEO{A}$ is just the adjugate of the lineal Gram matrix at $\GEO{A}$.
\begin{theorem}[B.~D.~S.~McConnell] \label{thm:bdsmc}
Any seven real numbers $f_\msf{ABC}\fit$, $f_\msf{ABD}\fit$, $f_\msf{ACD}\fit$, $f_\msf{BCD}\fit$, $f_\msf{AB|CD}\fit$, $f_\msf{AC|BD}\fit$, $f_\msf{AD|BC}$ $\ge 0$ are equal to the areas of the exterior (times $2$) and interior (times $4$) faces of a non-degenerate Euclidean tetrahedron $\GEO{ABCD}$ if \& only if they satisfy Yetter's identity $\Xi_f = 0$, the $18$ tetrahedron inequalities $\Tau_{\!f} \ge 0$, and yield a Gramian at $\GEO{A}$ (or any other vertex) $\Gamma_{\!f^2}[\msf A] > 0$.
This tetrahedron is unique up to isometry.
\end{theorem}
\begin{proof}
The necessity of the stated conditions was established above.
To prove sufficiency, note these conditions together with Lemma \myref{thm:2x2minors} show that the Gram matrix $\mbf G_\msf A \coloneq \mbf G_F[\msf A]$ computed from the areas via Eq.~(\myref{eq:gramatAa}) with $F \coloneq f^2$ is positive definite by Sylvester's criterion.
Hence coordinates for the cross-products it represents are obtained by diagonalizing it as $\mbf G_\msf A = \mbf U\fit \mathbi\Lambda\fit \mbf U^\top$, letting $\mbf V \coloneq \mathbi\Lambda^{1/2\,} \mbf U^\top$, and setting
\begin{equation}
\mbf p ~\coloneq~ \mbf v_1 \leftrightarrow \VCP{AB}{AC} ~,\quad 
\mbf q ~\coloneq~ \mbf v_2 \leftrightarrow \VCP{AD}{AB} ~,\quad 
\mbf r ~\coloneq~ \mbf v_3 \leftrightarrow \VCP{AC}{AD} ~, \vspace{-0.5ex}
\end{equation}
where $\mbf v_1, \mbf v_2, \mbf v_3$ are the columns of $\mbf V$.
The dot products among these coordinate vectors will then reproduce the matrix $\mbf G_\msf A$ exactly.
To convert the cross products' coordinates into those of their component vectors, observe first that the cross products among any three vectors $\mbf b, \mbf c, \mbf d \in \mbb R^3$ are the columns $\mbf r, \mbf q\fit[0.1em], \mbf p$ of the adjugate matrix $\mathbf{Adj}[\mbf b, \mbf c, \mbf d]$.
Thus the well-known fact that the adjugate of the adjugate of a non-singular square matrix is the original matrix times its determinant, together with the fact that the determinant of the adjugate of a $3\times3$ matrix is the square of the determinant of the original matrix, establishes that the coordinates of the vertices of the tetrahedron $\mbf a, \mbf b, \mbf c, \mbf d$ are given by
\vspace{-0.5ex}
\begin{equation}
\mbf a ~=~ \mbf 0 ~,\quad \mbf b ~=~ \mbf p \Xprod \mbf q \,/\, t ~,\quad \mbf c ~=~ \mbf r \Xprod \mbf p \,/\, t ~,\quad \mbf d ~=~ \mbf q \Xprod \mbf r \,/\, t
\end{equation}
where $t \coloneq \sqrt{|\mrm{det}(\mbf V)|} = \sqrt[4\;]{\Gamma_{\!F}[\msf{A}]} > 0$.
Uniqueness up to isometry follows because the vectors $\mbf b - \mbf a \leftrightarrow \VEC{AB}$, $\mbf c - \mbf a \leftrightarrow \VEC{AC}$ \& $\mbf d - \mbf a \leftrightarrow \VEC{AD}$ are translation independent, while\linebreak[4] rot\-at\-ing them rotates their cross products identically without changing the matrix of dot products $\mbf G_\msf A\fit$.
Equation (\myref{eq:gramatA}), however, shows that $\mbf G_\msf A$ uniquely determines the six squared areas in $\mbf G_F[\msf A]$, which also determine the seventh $F_\msf{BCD}$ via Yetter's identity.%
\footnote{
Uniqueness can also be proven by showing that the squared distances computed from the coordinates, when inserted into Eqs.~(\myref{eq:cmd3}), (\myref{eq:talata}) and their permutations, reproduce the given areas squared, and that the Jacobian of this map\-ping from squared distances to squared areas satisfies $\mrm{det}\big( \mbf J^\top \mbf J \big) = 28\, \CMD_D[\msf A, \msf B, \msf C, \msf D]^4$ $>$ $0$ (as will be shown in Appendix \myref{sec:pmap} of Part IV).
} 
\end{proof}
\noindent Note that $\Gamma_{\!f^2}[\msf A] = 0$ if any tetrahedron inequality saturates, so the given conditions imply $\Tau_{\!f} > 0$.
The example $f_\msf{ABC} = 9$, $f_\msf{ABD} = 10$, $f_\msf{ACD} = 17$, $f_\msf{BCD} = 14$ \& $f_\msf{AB|CD} = \smash{\sqrt{261}}$, $f_\msf{AC|BD} = \smash{\sqrt{76}}$, $f_\msf{AD|BC} = \smash{\sqrt{329}}$ shows the Gramians can be negative even when Yetter's identity and all $18$ tetrahedron inequalities are strictly satisfied. \smallskip

In the course of his pursuit of \textquote{hedronometry,} McConnell has also over a span of better than three decades used Yetter's identity (which he independently rediscovered) to rewrite the polynomial $\Gamma_{\!F}[\msf A]$ in a variety of ingenious ways so as to make it symmetric under vertex permutations and look in some sense more like Heron's formula.
The analogies between his formulae and Heron's are however not compelling, in that they do not build upon the intimate connection between Heron's formula and the in-circle of the triangle seen in Fig.~\myref{fig:heron}.
This in turn is the basis for the extension of Heron's formula to be derived in Part II of this series of papers.

{\small

} 

\enlargethispage{0.75in} 

\begin{center}
\medskip\rule{70pt}{2pt}\medskip
\end{center}

\begin{displayquote}\small
As to the need of improvement there can be no question whilst the reign of Euclid continues.
My own idea of a useful course is to begin with arithmetic, and then not Euclid but algebra.
Next, not Euclid, but practical geometry, solid as well as plane; not demonstration, but to make acquaintance.
Then not Euclid, but elementary vectors, conjoined with  algebra, and applied to geometry.
Addition first; then the scalar product. Elementary calculus should go on simultaneously, and come into the vector algebraic geometry after a bit.
Euclid might be an extra course for learned men, like Homer.
But Euclid for children is barbarous.

\smallskip\raggedleft\textit{Oliver Heaviside, 1893}
\end{displayquote}

\medskip
\begin{displayquote}\small
Geometry without algebra is dumb! Algebra without geometry is blind!

\raggedleft\textit{David Hestenes \& Garret Sobczyk, 1984}
\end{displayquote}

\mypart{\mytitle} 

\begin{abstract}
This is the second part of a series of four papers in The $\Pi\mrm{ME}$ Journal.
It begins by showing that the \textquote{in-touch points} at which a tetrahedron's in-sphere touches its exterior faces divide those faces into six pairs of congruent triangles.
The \TDEF{natural parameters} of a tetra\-hedron are then defined as the common areas of these six congruent pairs.
This geometric definition may be expressed algebraically by six simple rational functions of the areas of the tetrahedron's seven faces (as defined in Part I).
The denominators of these rational functions are the exterior surface area of the tetrahedron, while each numerator is a product of two factors, one of which is the deviation from saturation of one of the tetrahedron inequalities of Part I and the other of which is the associated non-degeneracy condition.
This leads to an algebraic definition of the corresponding \TDEF{inverse natural parameters} as the rational functions obtained by replacing those numerators by the product of the deviations of the complementary pair of tetrahedron inequalities from saturation.
They turn out to be related to the areas of the triangles into which the exterior faces are divided by the \textquote{ex-touch points} of the tetrahedron's ex-spheres.
The product of such a complementary pair of natural and inverse natural parameters equals the  square of the product of the corresponding inter-vertex distance with the in-radius, and since the inverse parameters can also be expressed as rational functions of the natural parameters, this yields the desired extension of Heron's formula to tetrahedra.
That formula gives the fourth power of the volume as the product of the squared exterior surface area with a homogeneous quartic polynomial in the natural parameters, which may in turn be expressed as the negative of a simple  $4\times4$ symmetric determinant.
A series of remarks follows which give various perspectives on the formula, and present some other ways of expressing the fourth power of the volume as a polynomial or rational function of the natural and inverse parameters together.
The paper closes with a conjecture as to how the formula extends to $n$-dimensional simplices for all $n > 3$.
\end{abstract}

\section{The natural parameters of a tetrahedron}
The first step towards a formula for the volume of a tetrahedron that can justly be called a \underline{natural} extension of Heron's formula is to find parameters which determine its facial areas in much the same way that the Heron parameters were shown to determine the edge lengths of a triangle in Section \myref[1]{sec:intro}.\footnote{A Roman numeral followed by a colon specifies a reference to a numbered entity in another paper of this series, in this case Section \myref{sec:intro} of Part I.}
In analogy with the in-touch points of a triangle shown in Fig.~\myref[1]{fig:heron}, we shall denote the in-touch points of the in-sphere of a tetrahedron $\GEO{ABCD}$ by $\GEO{J}$, $\GEO{K}$, $\GEO{L}$ \& $\GEO{N}$ (\textquote{\textsf{M}} was reserved for the Monge point, although it plays no role here).
The centrality of the in-sphere and in-touch points to the geometry of tetrahedra may be demonstrated by a construction which parameterizes the
set of all non-degenerate tetrahedra as follows:
\begin{enumerate}
\item Choose a sphere of radius $r > 0$, centered on e.g.~the origin, as the in-sphere;
\item Choose four non-coplanar points $\GEO{J}$, $\GEO{K}$, $\GEO{L}$ \& $\GEO{N}$ on this sphere such that the plane through any two of them and the center of the sphere separates the remaining two; these four points will become the in-touch points of the tetrahedron;
\item Take the planes tangent to the sphere at these four points, and intersect them three-at-a-time to get the vertices $\GEO{A}$, $\GEO{B}$, $\GEO{C}$ \& $\GEO{D}$ of the tetrahedron.
\end{enumerate}
The results of this construction, carried out in the \texttt{GeoGebra} online dynamic geometry system \cite{GeoGebra:2021}, are shown in Fig.~\myref{fig:insphere} above.

\begin{figure}
\input{Fig_II-1.tex}
\label{fig:insphere}
\end{figure}
\smallskip
The first item of business is to establish the following:
\begin{lemma} \label{thm:intouch}
The twelve triangles into which the in-touch points divide the exterior faces of a tetrahedron occur in six congruent pairs, where each pair shares a common edge of the tetrahedron.
Moreover, the line segment between each pair of in-touch points is perpendicular to the common edge of the two faces those in-touch points lie in.
\end{lemma}
\begin{proof}
Using  the point labels in Fig.~\myref{fig:insphere}, the first part of the lemma may be proven by noting that the vector from $\GEO{A}$ (say) to the in-center $\GEO{\fit[0.1em] I\fit[0.1em]}$ can be written in two ways, i.e.
\begin{equation}
\VEC{AL} \,\Vplus\, \VEC{L\fit[0.1em] I\fit[0.1em]} ~=~ \VEC{AI\fit} ~=~ \VEC{AN} \,\Vplus\, \VEC{N\fit[0.1em] I\fit[0.1em]} ~,
\end{equation}
from which it follows that
\begin{equation}
\MAG{AL}^2 \,+\, \MAG{L\fit[0.1em] I\fit[0.1em]}^2 +\,\, 2\,  \VEC{AL} \Dprod \VEC{L\fit[0.1em] I\fit[0.1em]} ~=~ \MAG{AN}^2 \,+\, \MAG{N\fit[0.1em] I\fit[0.1em]}^2 \,+\, 2\, \VEC{AN} \Dprod \VEC{N\fit[0.1em] I\fit[0.1em]} ~.
\end{equation}
But $\VEC{AL} \Dprod \VEC{L\fit[0.1em] I\fit[0.1em]}$ $= 0 =$ $\VEC{AN} \Dprod \VEC{N\fit[0.1em] I\fit[0.1em]}$ since $\VEC{AL}$ \& $\VEC{AN}$ lie in the planes of the faces $\GEO{ABD}$ \& $\GEO{ABC}$ resp., while \smash{$\VEC{L\fit[0.1em] I\fit[0.1em]}$ \& $\VEC{N\fit[0.1em] I\fit[0.1em]}$} are perpendicular to those faces with a common length equal to the in-radius $r$ by definition.
This shows that $\MAG{AL}$ $=$ $\MAG{AN}$, and similarly $\MAG{BL}$ $=$ $\MAG{BN}$, so that $\GEO{ABL}$ is congruent to $\GEO{ABN}$ as claimed.
In an analogous fashion, one finds all the distances from each vertex to its three adjacent in-touch points are equal, i.e.
\begin{equation} \begin{aligned}
\MAG{AK} ~=~ \MAG{AL} ~=~ \MAG{AN} & ~,\quad \MAG{BJ} ~=~ \MAG{BL} ~=~ \MAG{BN} ~, \\[3pt]
\MAG{CJ} ~=~ \MAG{CK} ~=~ \MAG{CN} & ~,\quad \MAG{DJ} ~=~ \MAG{DK} ~=~ \MAG{DL} ~,
\end{aligned} \end{equation}
which implies the congruence of all the remaining pairs of triangles, where each pair meets in an edge of the tetrahedron and the triangles in each pair are spanned by that edge together with its adjacent in-touch points.

The second part of the lemma is likewise easily proven using the orthogonality of the vectors \smash{$\VEC{\fit[0.1em] I\fit[0.1em] L}$ \& $\VEC{\fit[0.1em] I\fit[0.1em] N}$} to the faces $\GEO{ABD}$ \& $\GEO{ABC}$, resp., and hence to their common edge $\GEO{AB}$:
\begin{equation}
\VEC{AB} \Dprod \VEC{NL} ~=~ \VEC{AB} \Dprod \big( \VEC{N\fit[0.1em] I\fit[0.1em]} \,\Vplus\, \VEC{\fit[0.1em] I\fit[0.1em]L}  \big) ~=~ \VEC{AB} \Dprod \VEC{N\fit[0.1em] I\fit[0.1em]} \,+\, \VEC{AB} \Dprod \VEC{\fit[0.1em] I\fit[0.1em] L} ~=~ 0 + 0 ~=~ 0 ~.
\end{equation}
The proofs for the pairs of ``contact triangles'' (as we shall call them) at the other five edges are similar. 
\end{proof}
\noindent This leads to the following tetrahedral analogues of the Heron parameters of a triangle.
\begin{definition} \label{def:natpar}
The \TDEF{natural parameters} of a tetrahedron are twice the common areas of each of these six pairs of congruent triangles, specifically:
\begin{align}
u ~\coloneq~ 2 \MAG{ABL} ~=~ 2 \MAG{ABN} ~,\quad v ~\coloneq~ 2 \MAG{ACK} ~=~ 2 \MAG{ACN} ~, \nonumber \\[3pt]
w ~\coloneq~ 2 \MAG{ADK} ~=~ 2 \MAG{ADL} ~,\quad x ~\coloneq~ 2 \MAG{BCJ\fit} ~=~ 2 \MAG{BCN} ~, \\[3pt] \nonumber
y ~\coloneq~ 2 \MAG{BDJ\fit} ~=~ 2 \MAG{BDL} ~,\quad z ~\coloneq~ 2 \MAG{CDJ\fit} ~=~ 2 \MAG{CDK} ~.
\end{align}
\end{definition}
\noindent Then because each exterior face of the tetrahedron is partitioned into three contact triangles by its in-touch point, the natural parameters satisfy the following system of linear equations:
\begin{equation} \begin{aligned}
u + v + x ~=~ 2 \MAG{ABC} & ~,\qquad  u + w + y ~=~ 2 \MAG{ABD} ~, \\[3pt]
v + w + z ~=~  2 \MAG{ACD} & ~,~\qquad x + y + z ~=~ 2 \MAG{BCD} ~.
\end{aligned} \label{eq:natparlinsys} \end{equation}
The problem is that, unlike the triangle where the three Heron parameters are connected to the edge lengths by a non-singular system of three linear equations, here there are only four equations in the six unknowns $u, v, w, x, y$ \& $z$.

\smallskip
To obtain their values, let the areal vectors of the contact triangles (times $2$) be:
\begin{equation}
\fit[-1em] \begin{array}{llll}
\mbf u_\msf{AB;C} \coloneq \VCP{NB}{NA} \,, &\fit[-0.5em] \mbf u_\msf{AB;D} \coloneq \VCP{LA}{LB} \,,
&\fit[-0.5em] \mbf v_\msf{AC;B} \coloneq \VCP{NA}{NC} \,, &\fit[-0.5em] \mbf v_\msf{AC;D} \coloneq \VCP{KC}{KA} \,, \\[9pt]
\mbf w_\msf{AD;C} \coloneq \VCP{KA}{KD} \,, &\fit[-0.5em] \mbf w_\msf{AD;B} \coloneq \VCP{LD}{LA} \,,
&\fit[-0.5em] \mbf x_\msf{BC;A} \coloneq \VCP{NC}{NB} \,, &\fit[-0.5em] \mbf x_\msf{BC;D} \coloneq \VCP{JB}{JC} \,, \\[9pt]
\mbf y_\msf{BD;A} \coloneq \VCP{LB}{LD} \,, &\fit[-0.5em] \mbf y_\msf{BD;C} \coloneq \VCP{JD}{JB} \,,
&\fit[-0.5em] \mbf z_\msf{CD;A} \coloneq \VCP{KD}{KC} \,, &\fit[-0.5em] \mbf z_\msf{CD;B} \coloneq \VCP{JC}{JD} \,.
\end{array} \fit[-0.5em] \label{eq:natvec}
\end{equation}
Note that the order of the factors in each cross-product has been chosen so as to ensure that these are all outwards-pointing vectors if the oriented volume of $\GEO{ABCD}$ is positive, or inwards-pointing if it is negative.
Then the sum of the areal vectors of the two contact triangles sharing a common edge is e.g.
\begin{align}
\fit[-1em] \mbf u_\msf{AB;D} \Vplus \mbf u_\msf{AB;C} ~=~ &
\VEC{LA} \Xprod \VEC{LB} \,\Vdiff\, \VEC{NA} \Xprod \VEC{NB} ~=~ \big( \VEC{LN} \Vplus \VEC{NA} \big) \Xprod \big( \VEC{LN} \Vplus \VEC{NB} \big) \Vdiff\, \VEC{NA} \Xprod \VEC{NB} \nonumber \\[6pt] =~  &
\VEC{NA} \Xprod \VEC{LN} \,\Vplus\, \VEC{LN} \Xprod \VEC{NB} ~=~ \VEC{LN} \Xprod \big( \VEC{NB} \,\Vdiff\, \VEC{NA} \big) \,=~ \VEC{LN} \Xprod \VEC{AB} ~.
\end{align} 
Since $\lVert \mbf u_\msf{AB;D} \rVert = \lVert \mbf u_\msf{AB;C} \rVert = u$ and \smash{$\VEC{LN}$} $\perp$ \smash{$\VEC{AB}$} by Lemma \myref{thm:intouch}, it follows that
\pagebreak[2]
\begin{align}
& \fit[-2.5em] \NCP{LN}{AB}^2 \,=~ \NMV{LN}^2 \NMV{AB}^2 \,=~ \lVert \mbf u_\msf{AB;D} \Vplus \mbf u_\msf{AB;C} \rVert^2 ~=
\label{eq:lnxab2} \\[2pt] \nonumber &
\lVert \mbf u_\msf{AB;C} \rVert^2 +\, \lVert \mbf u_\msf{AB;D} \rVert^2 -\, 2\, \lVert \mbf u_\msf{AB;C} \rVert \lVert \mbf u_\msf{AB;D} \rVert \cos( \varphi_\msf{AB} ) ~=~ 2\, u^2 \big( 1 -  \cos( \varphi_\msf{AB} ) \big) ~,
\end{align}
where $\varphi_\msf{AB}$ is the dihedral angle between $\GEO{ABC}$ \& $\GEO{ABD}$, and the \textquote{$-$} in front of the cosine is because $\varphi_\msf{AB}$ is the angle between $\mbf u_\msf{AB;C}$ \& $\Vdiff \mbf u_\msf{AB;D}$ (or vice versa).
By the areal law of cosines (\myref[1]{eq:cosines}), however, this cosine is equal to
\begin{subequations} \begin{align}
&  \cos( \varphi_\msf{AB} ) ~=~ 
\frac{\NCP{AB}{AC}^2 +\, \NCP{AB}{AD}^2 -\, \NCP{AB}{CD}^2} 
{2 \NCP{AB}{AC} \NCP{AB}{AD}\hit[2.5ex]}
\\[3pt] \Longleftrightarrow\quad &
1 \,-\, \cos( \varphi_\msf{AB} ) ~=~ 
\frac{\NCP{AB}{CD}^2 -\, \raisebox{0.5\upbgdl}{$\big($} \! \NCP{AB}{AC} \,-\, \NCP{AB}{AD} \! \raisebox{0.5\upbgdl}{$\big)$}^{\!2}} 
{2 \NCP{AB}{AC} \NCP{AB}{AD}\hit[2.5ex]} ~,
\end{align} \end{subequations}
and plugging that into Eq.~(\myref{eq:lnxab2}) then solving for $u^2$ gives
\begin{equation}
u^2 ~=~ \frac{\NMV{\fit LN}^2 \NMV{AB}^2 \NCP{AB}{AC} \NCP{AB}{AD}}{\NCP{AB}{CD}^2 -\, \raisebox{0.5\upbgdl}{$\big($} \! \NCP{AB}{AC} \,-\, \NCP{AB}{AD} \! \raisebox{0.5\upbgdl}{$\big)$}^{\!2}} ~.
\end{equation}
To finish the job a formula for $\NMV{\fit LN}^2$ is needed, and it is
\begin{equation} \begin{aligned}
\NMV{\fit LN\fit}^2 =~ &\> \LVert \VEC{\fit[0.1em] I\fit[0.1em] N} \Vdiff \VEC{\fit[0.1em] I\fit[0.1em] L\fit}\, \RVert^2 =~ \NMV{\fit[0.1em] I\fit[0.1em] N\fit}^2 +\fit \NMV{\fit[0.1em] I\fit[0.1em] L}^2 -\fit 2\, \VEC{\fit[0.1em] I\fit[0.1em] L} \Dprod \VEC{\fit[0.1em] I\fit[0.1em] N} ~=~ 2\, r^2\, \big( 1 + \cos( \varphi_\msf{AB} ) \big) \\[2pt] 
=~ &\> 2\, r^2\, \frac{\raisebox{0.5\upbgdl}{$\big($} \! \NCP{AB}{AC} \,+\, \NCP{AB}{AD} \! \raisebox{0.5\upbgdl}{$\big)$}^{\!2} -\, \NCP{AB}{CD}^2} 
{2 \NCP{AB}{AC} \NCP{AB}{AD}\hit[2.5ex]} ~, 
\end{aligned} \end{equation}
where $r = \NMV{\fit I\fit L\fit} = \NMV{\fit I\fit N\fit}$ is the in-radius and the change in the sign of the cosine has the same explanation as above.
This leads to the following relatively simple formulae:
\begin{align}
u ~=~ {}&{} r\, \NMV{AB}\, \cot( \varphi_\msf{AB} / 2 ) ~=~ r\, \NMV{AB}\, \sqrt{\frac{1 \,+\, \cos( \varphi_\msf{AB} )}{1 \,-\, \cos( \varphi_\msf{AB} )}} \label{eq:ucotan} \nonumber \\[1pt] 
=~ {}&{} r\, \NMV{AB}\, \sqrt{ \frac 
{\raisebox{0.5\upbgdl}{$\big($} \! \NCP{AB}{AC} \,+\, \NCP{AB}{AD} \! \raisebox{0.5\upbgdl}{$\big)$}^{\!2} -\, \NCP{AB}{CD}^2}
{\NCP{AB}{CD}^2 -\, \raisebox{0.5\upbgdl}{$\big($} \! \NCP{AB}{AC} \,-\, \NCP{AB}{AD} \!  \raisebox{0.5\upbgdl}{$\big)$}^{\!2}} } \\[3pt] 
=~ {}&{} r\, \NMV{AB}\, \sqrt{ \frac{\Tau_0[\msf A, \msf B]\, \Tau_1[\msf A, \msf B]} 
{\Tau_2[\msf A, \msf B]\, \Tau_3[\msf A, \msf B]} } \nonumber 
\end{align}
Here \smash{$\Tau_0[\msf A, \msf B] \ge 0$} is the non-degeneracy factor and \smash{$\Tau_k[\msf A, \msf B] \ge 0$} ($k = 1,2,3$) are the deviations of the tetrahedron inequalities from saturation defined in Eq.~(\myref[1]{eq:taus}), so the quantity in the square root is non-negative.
Similar expressions can of course be derived for the other parameters $v, w, x, y$ \& $z$ via the definitions given in Eq.~(\myref{eq:natvec}).

This expression may be further simplified via the trigonometric identity $\cot(\vartheta/2) = \csc(\vartheta) + \cot(\vartheta)$, where the sine in $\csc( \varphi_\msf{AB} )$ is obtained from the areal law of sines (\myref[1]{eq:sines}).
It then follows from Eq.~(\myref{eq:ucotan}) that
\pagebreak[2]
\begin{align}
u ~=~ \big(\, r  {}&{} \NMV{AB} \!\bigm/\! \sin ( \varphi_\msf{AB} ) \,\big) \big( 1 + \cos( \varphi_\msf{AB} ) \big) ~= \\[4pt]
r\, \NMV{AB} {}&{} 
\frac{\NCP{AB}{AC} \! \NCP{AB}{AD}} 
{\hit[2.5ex] \NMV{AB} \raisebox{0.5\upbgdl}{$\big|$} \VEC{AB} \Dprod \raisebox{0.5\upbgdl}{$\big($} \VCP{AC}{AD} \raisebox{0.5\upbgdl}{$\big) \big|$}}
\,  \Bigg( 1 \,+\, 
\frac{\NCP{AB}{AC}^2 +\, \NCP{AB}{AD}^2 -\, \NCP{AB}{CD}^2} 
{\hit[2.5ex] 2 \NCP{AB}{AC} \NCP{AB}{AD}} \Bigg) \nonumber \\[4pt] 
=~ {}&{} r\: \frac{2 \NCP{AB}{AC} \! \NCP{AB}{AD} \,+\, \NCP{AB}{AC}^2 +\, \NCP{AB}{AD}^2 -\, \NCP{AB}{CD}^2} 
{\hit[2.5ex] 2\; \raisebox{0.5\upbgdl}{$\big|$} \VEC{AB} \Dprod \raisebox{0.5\upbgdl}{$\big($} \VCP{AC}{AD} \raisebox{0.5\upbgdl}{$\big) \big|$}} \nonumber \\[3pt]  
=~ {}&{} \frac{\raisebox{0.5\upbgdl}{$\big($} \! \NCP{AB}{AC} \fit+\fit \NCP{AB}{AD} \! \raisebox{0.5\upbgdl}{$\big)$}^{\!2} -\, \NCP{AB}{CD}^2}{2\, s} ~=~ \frac{\Tau_0[\msf A, \msf B]\, \Tau_1[\msf A, \msf B]}{2\,s} ~, \nonumber
\end{align}
where the well-known \cite{Apostol:2012} relation $r = t \fit/\fit s$ $\coloneq$ \smash{$\raisebox{0.5\upbgdl}{\big|}\fit \VEC{AB} \Dprod \raisebox{0.5\upbgdl}{$\big($} \VCP{AC}{AD} \raisebox{0.5\upbgdl}{$\big)$} \raisebox{0.5\upbgdl}{\big|} \raisebox{0.25\upbgdl}{$\,\bigm/\,$} s$} was used to obtain the last line with $s \coloneq 2 \MAG{ABC} + 2 \MAG{ABD} + 2 \MAG{ACD} + 2 \MAG{BCD}\hit[2ex]$ equal to twice the exterior surface area.

\smallskip 
A little more generally, we obtain the following expressions for all six of the natural parameters in terms of the $\Tau\,$'s as defined in Eq.~(\myref[1]{eq:taus}):
\begin{proposition} \label{thm:natpar}
The natural parameters of a tetrahedron $\GEO{ABCD}$ with $s > 0$ are
\begin{equation} \begin{aligned}
\fit[-0.5em] u \:=\: \frac{\Tau_0[\msf A, \msf B]\, \Tau_1[\msf A, \msf B]}{2\,s} \,,\fit[0.5em] v \:=\: \frac{\Tau_0[\msf A, \msf C]\, \Tau_1[\msf A, \msf C]}{2\,s} \,,\fit[0.5em] w \:=\: \frac{\Tau_0[\msf A, \msf D]\, \Tau_1[\msf A, \msf D]}{2\,s} \,, \\[3pt]
\fit[-0.5em] z \:=\: \frac{\Tau_0[\msf C, \msf D]\, \Tau_1[\msf C, \msf D]}{2\,s} \,,\fit[0.5em] y \:=\: \frac{\Tau_0[\msf B, \msf D]\, \Tau_1[\msf B, \msf D]}{2\,s} \,,\fit[0.5em] x \:=\: \frac{\Tau_0[\msf B, \msf C]\, \Tau_1[\msf B, \msf C]}{2\,s} \,,
\end{aligned} \label{eq:natparratfun} \end{equation}
and $0$ if $s = 0$.
\end{proposition}
\begin{proof}
This follows simply from the preceding formula for $u$ together with analogous calculations starting from the definitions given in Eq.~(\myref{eq:natvec}).
\end{proof}

\section{The inverse natural parameters of a tetrahedron}
In light of these formulae for the natural parameters, we will also define:
\begin{definition} \label{def:invnatpar}
The \TDEF{inverse natural parameters\/} complementary to $u, v, w, x, y, z$ are:\pagebreak[2]
\begin{equation} \label{eq:invpar} \begin{aligned}
\fit[-0.5em] & \tilde u \:\coloneq\: \frac{\Tau_2[\msf A, \msf B]\, \Tau_3[\msf A, \msf B]}{2\,s} \,,\fit[0.5em] \tilde v \:\coloneq\: \frac{\Tau_2[\msf A, \msf C]\, \Tau_3[\msf A, \msf C]}{2\,s} \,,\fit[0.5em] \tilde w \:\coloneq\: \frac{\Tau_2[\msf A, \msf D]\, \Tau_3[\msf A, \msf D]}{2\,s} \,, \\[3pt]
\fit[-0.5em] & \tilde z \:\coloneq\: \frac{\Tau_2[\msf C, \msf D]\, \Tau_3[\msf C, \msf D]}{2\,s} \,,\fit[0.5em] \tilde y \:\coloneq\: \frac{\Tau_2[\msf B, \msf D]\, \Tau_3[\msf B, \msf D]}{2\,s} \,,\fit[0.5em] \tilde x \:\coloneq\: \frac{\Tau_2[\msf B, \msf C]\, \Tau_3[\msf B, \msf C]}{2\,s} \;
\end{aligned} \end{equation}
for $s > 0$, and $0$ otherwise.
\end{definition}
\noindent The following corollary to Proposition \myref{thm:natpar} justifies the \textquote{inverse} qualifier: 
\begin{corollary} \label{thm:r2d2}
The squared in-radius $r^2$ times the squared edge lengths are equal to the products of complementary pairs of natural \& inverse natural parameters, specifically:
\begin{equation} \label{eq:r2d2} \begin{aligned}
r^2 \MAG{AB}^2 ~=~ u\tilde u ~,\qquad r^2 \MAG{AC}^2 ~=~ & v\tilde v ~,\qquad r^2 \MAG{AD}^2 ~=~ w\tilde w ~, \\
r^2 \MAG{CD}^2 ~=~ z\tilde z ~,\qquad r^2 \MAG{BD}^2 ~=~ & y\tilde y ~,\qquad r^2 \MAG{BC}^2 ~=~ x\tilde x ~.
\end{aligned} \end{equation}
\end{corollary}
\begin{proof}
These relations follow easily from Proposition \myref{thm:natpar} and Definition \myref{def:invnatpar} together with Eq.~(\myref{eq:ucotan}) and the analogous equations for $v, w, x, y$ \& $z$.
\end{proof}
\noindent Note that $t^2\, \smash{\MAG{AB}}^2 = s^2 u \tilde u = \tfrac14\, {\Tau_0[\msf A, \msf B]\, \Tau_1[\msf A, \msf B]}\, {\Tau_2[\msf A, \msf B]\, \Tau_3[\msf A, \msf B]}$ etc.~are the $2\times2$ principal minors of the Gram matrices by Lemma \myref[1]{thm:2x2minors}, and that these also equal $\big( 4 \MAG{ABC} \! \MAG{ABD}\fit \sin(\varphi_\msf{AB}) \smash{\big)}^{2}$ etc.~by the areal law of sines (\myref[1]{eq:sines}).

\smallskip
The following further corollary summarizes some of the algebraic identities which connect the natural and inverse natural parameters with the seven areas.
\begin{corollary} \label{thm:ids}
With everything defined as above, the following identities hold:
\begin{subequations}
\begin{gather}
2\, (u + v + w + x + y + z) ~=~ s \,+\, 2\, \Xi\, / \,s ~; \label{eq:ida} \\[0.5ex]
\begin{alignedat}{4} \label{eq:idb}
u + v + x & ~=~ & 2\, \MAG{ABC} \,+\, \Xi \!\bigm/\! (2s) \,,\fit[0.5em] 
u + w + y & ~=~ & 2\, \MAG{ABD} \,+\, \Xi \!\bigm/\! (2s) \,, \\[-1pt]
v + w + z & ~=~ & 2\, \MAG{ACD} \,+\, \Xi \!\bigm/\! (2s) \,,\fit[0.75em]
x + y + z & ~=~ & 2\, \MAG{BCD} \,+\, \Xi \!\bigm/\! (2s) \,;
\end{alignedat} \\[1ex] 
\begin{aligned}
u + v + w ~=~ & \MAG{ABC} \,+\, \MAG{ABD} \,+\, \MAG{ACD} \,-\, \MAG{BCD} \,+\, \Xi \!\bigm/\! s ~, \\[-2pt]
u + x + y ~=~ & \MAG{BCD} \,+\, \MAG{ABC} \,+\, \MAG{ABD} \,-\, \MAG{ACD} \,+\, \Xi \!\bigm/\! s ~, \\[-2pt]
v + x + z ~=~ & \MAG{ACD} \,+\, \MAG{BCD} \,+\, \MAG{ABC} \,-\, \MAG{ABD} \,+\, \Xi \!\bigm/\! s ~, \\[-2pt]
w + y + z ~=~ & \MAG{ABD} \,+\, \MAG{ACD} \,+\, \MAG{BCD} \,-\, \MAG{ABC} \,+\, \Xi \!\bigm/\! s ~;
\end{aligned} \label{eq:idc} \\[0.5ex] 
u \,-\, z ~=~ \MAG{ABC} \,+\, \MAG{ABD} \,-\, \MAG{ACD} \,-\, \MAG{BCD} ~, \nonumber \\[-2pt] 
v \,-\, y ~=~ \MAG{ABC} \,+\, \MAG{ACD} \,-\, \MAG{ABD} \,-\, \MAG{BCD} ~, \label{eq:idd} \\[-2pt]
w \,-\, x ~=~ \MAG{ABD} \,+\, \MAG{ACD} \,-\, \MAG{ABC} \,-\, \MAG{BCD} ~; \nonumber \\
\big( v + w + x + y - \Xi / s \big)^{2} -\, 4\, uz ~=~ 16\, \MAG{AB \tmv CD}^2 ~, \nonumber \\[-2pt]
\big( u + w + x + z - \Xi / s \big)^{2} -\, 4\, vy ~=~ 16\, \MAG{AC \tmv BD}^2 ~,  \label{eq:ide} \\[-2pt]
\big( u \fit+\fit v \fit+\fit y \fit+\fit z -\, \Xi / s \big)^{2} - 4\, wx ~=~ 16\, \MAG{AD \tmv BC}^2 ~; \nonumber \\[0.5ex]
\begin{alignedat}{4} \label{eq:idf}
(u + \tilde u)\, s \, &=& \>\> 8\, \MAG{ABC} \MAG{ABD} \,,\fit[0.25em]
(z \,+\, \tilde z)\, s \, &=& \>\> 8\, \MAG{ACD} \MAG{BCD} \,, \\[0pt]
(v + \tilde v)\, s \, &=& \>\> 8\, \MAG{ABC} \MAG{ACD} \,,\fit[0.5em]
(y + \tilde y)\, s \, &=& \>\> 8\, \MAG{ABD} \MAG{BCD} \,, \\[0pt]
(w + \tilde w)\, s \, &=& \>\> 8\, \MAG{ABD} \MAG{ACD} \,,\fit[0.5em]
(x + \tilde x)\, s \, &=& \>\> 8\, \MAG{ABC} \MAG{BCD} \,;
\end{alignedat} \\[1ex] 
\fit[-0.5em] \begin{alignedat}{4}
(u - \tilde u)\, s/2 \, &=& \fit[0.5em] \big( \VCP{AB}{AC} \big) \Dprod \big( \VCP{AB}{AD} \big) ,\fit[0.5em] (z - \tilde z)\, s/2 \, &=& \fit[0.5em] \big( \VCP{AC}{AD} \big) \Dprod \big( \VCP{BC}{BD} \big) , \fit[0.1em] \\[1pt]
(\tilde v - v)\, s/2 \, &=& \fit[0.5em] \big( \VCP{AB}{AC} \big) \Dprod \big( \VCP{AC}{AD} \big) ,\fit[0.4em] (\tilde y - y)\, s/2 \, &=& \fit[0.5em] \big( \VCP{AB}{AD} \big) \Dprod \big( \VCP{BC}{BD} \big) , \\[1pt]
(w - \tilde w)\, s/2 \, &=& \fit[0.5em] \big( \VCP{AB}{AD} \big) \Dprod \big( \VCP{AC}{AD} \big) ,\fit[0.4em] (x - \tilde x)\, s/2 \, &=& \fit[0.5em] \big( \VCP{AB}{AC} \big) \Dprod \big( \VCP{BC}{BD} \big) ;
\end{alignedat} \label{eq:idg} \\[0.0ex] 
\fit[-0.5em] s^2 -\, 2\, (\tilde u + \tilde v + \tilde w + \tilde x + \tilde y + \tilde z)\, s ~=~
32\, \Big( \MAG{AB \tmv CD}^2 \!+ \MAG{AC \tmv BD}^2 \!+ \MAG{AD \tmv BC}^2 \Big) + \, 4\, \Xi \,. \label{eq:idh}
\end{gather}
\end{subequations} \label{thm:algids}
\end{corollary}
\vspace{-5ex} 
\begin{proof}
The identities in Eq.~(\myref{eq:ida}), (\myref{eq:idc}) \& (\myref{eq:idd}) follow upon substituting for the areas therein using the identities in Eq.~(\myref{eq:idb}).
The latter identities just reiterate Eq.~(\myref{eq:natparlinsys}), but with the addition of the multiple of $\Xi$ needed to make them hold even if Yetter's identity is not satisfied and the natural parameters therein are replaced by the rational functions in Eq.~(\myref{eq:natparratfun}).
The first identity in Eq.~(\myref{eq:ide}) may likewise be proven by substituting for the natural parameters therein using Eq.~(\myref{eq:natparratfun}), which yields
\begin{align*}
&{} \fit[-1.2em] (v+w+x+y - \Xi/s)^2 -\, 4\, uz \\
=~ {}&{} \begin{aligned}[t] &{} 16\, \Big( \Big( \MAG{ABC} + \MAG{ABD} \Big) \Big( \MAG{ACD} + \MAG{BCD} \Big) +\, 4 \MAG{AB|CD}^2 \Big)^{\!2} \!\bigm/\! s^2 \\
{}&{} -\, 16\, \Big( \Big( \MAG{ABC} + \MAG{ABD} \Big)^{\!2} -\, 4 \MAG{AB|CD}^2 \Big) \Big( \Big( \MAG{ACD} + \MAG{BCD} \Big)^{\!2} -\, 4 \MAG{AB|CD}^2 \Big) \!\bigm/\! s^2 \end{aligned} \\
=~ {}&{} 64\, \Big( \MAG{ABC} + \MAG{ABD} + \MAG{ACD} + \MAG{BCD} \Big)^{\!2}\, \MAG{AB|CD}^2 \!\bigm/\! s^2 ~=~ 16\, \MAG{AB|CD}^2 .
\end{align*}
The other identities in Eq.~(\myref{eq:ide}) may be proven similarly.
The identities in Eq.~(\myref{eq:idg}) just restate the last line of Eq.~(\myref[1]{eq:cosines}), while Eq.~(\myref{eq:idf}$\,$) restates the complementary relation given following Eq.~(\myref[1]{eq:sines0}).
The last identity (\myref{eq:idh}) will be left as an exercise.
\end{proof}
\begin{remark} \label{rem:extouch}
The analogues of the inverse natural parameters for a triangle are
\begin{align} 
\tilde u ~=~ \frac{( a - b + c)\, ( a + b - c )}{4\, s} ~,\quad
\tilde v ~=~ & \frac{( -a + b + c )\, ( a + b - c )}{4\, s} ~, \\[3pt] \text{and}\quad
\tilde w ~=~ & \frac{( -a + b + c )\, ( a - b + c )}{4\, s} ~, \nonumber
\end{align}
where $a, b, c \ge 0$ are its edge lengths and $s = (a+b+c)/2$ is its semi-perimeter (as in Section \myref[1]{sec:intro}).
These may be shown to satisfy
\begin{equation}
\tilde u ~=~ v r / r_\msf C ~=~ w r / r_\msf B ~,\fit[1em]
\tilde v ~=~ u r / r_\msf C ~=~ w r / r_\msf A ~,\fit[1em]
\tilde w ~=~ u r / r_\msf B ~=~ v r / r_\msf A ~,
\end{equation}
where $r$, $r_\msf A = rs/u$, $r_\msf B = rs/v$, $r_\msf C = rs/w$ are the radii of the in-circle and ex-circles tangent to the edge of length $a, b, c$ respectively, and $u, v, w$ are the Heron parameters (cf.~Fig.~\myref[1]{fig:heron}).
\pagebreak[2]

For a tetrahedron, the relations between the in-radius and ex-radii given in e.g.~Refs.~\cite{Hajji:2018, Richardson:1902, Toda:2014}, together with Eq.~(\myref{eq:idc}), show that these quantities satisfy
\begin{equation} \begin{aligned}
r_\msf A ~=~ \tfrac12\, rs / (u+v+w) ~,\quad &~ r_\msf B ~=~ \tfrac12\, rs / (u+x+y) ~, \\[1pt]
r_\msf C ~=~ \tfrac12\, rs / (\fit v+\fit x+\fit z) ~,\quad &~ r_\msf D ~=~ \tfrac12\, rs / (w+y+z) ~.
\end{aligned} \end{equation}
where $r_\msf A$ is the radius of the ex-sphere tangent to the exterior face opposite $\GEO{A}$, etc.
Furthermore, Lemma \myref{thm:natinvnat} (below) can be used to show that
\begin{equation} \begin{alignedat}{2}
(\tilde x \!+\! \tilde y \!+\! \tilde z)\, s ~=~ & 2\, (u \!+\! v \!+\! w) (x \!+\! y \!+\! z) \,,\fit[1em]
(\tilde v \!+\! \tilde w \!+\! \tilde z)\,s  ~=~ & 2\, (u \!+\! x \!+\! y) (v \!+\! w \!+\! z) , \\[3pt] 
(\tilde u \!+\! \tilde w \!+\! \tilde y)\, s ~=~ & 2\, (v \!+\! x \!+\! z) (u \!+\! w \!+\! y) \,,\fit[1.1em]
 (\tilde u \!+\! \tilde v \!+\! \tilde x)\, s ~=~ & 2\, (w \!+\! y \!+\! z) (u \!+\! v \!+\! x) ,
\end{alignedat} \end{equation}
and hence $2 \MAG{ABC} = u + v + x =  (\tilde u + \tilde v + \tilde x)\, r_\msf D / r$, etc.
This suggests that the inverse natural parameters of a tetrahedron scaled by the ratios of the ex-radii to the in-radius may be interpreted as twice the areas of the triangles into which the \TDEF{ex-touch points\/} $\GEO{J_A\!}\,$, $\GEO{K_B\!}\,$, $\GEO{L_C\!}\,$ \& $\GEO{N_D\!}\,$ divide the tetrahedron's exterior faces, e.g.
\begin{equation}
\frac{r_{\msf D\hit[1.5ex]}}{r}\, \tilde u ~=~ 2 \MAG{ABN_\msf D} ~,\quad \frac{r_{\msf D\hit[1.5ex]}}{r}\, \tilde v ~=~ 2 \MAG{ACN_\msf D} ~,\quad \frac{r_{\msf D\hit[1.5ex]}}{r}\, \tilde x ~=~ 2 \MAG{BCN_\msf D} ~.
\end{equation}
This hypothesis has been shown to hold numerically in randomly generated tetrahedra, thereby obtaining a "generic" proof of its correctness.
\end{remark}

It is readily verified that the non-negativity of all the natural and inverse natural parameters is equivalent to all $18$ tetrahedron inequalities holding, as long as the areas from which they were obtained are likewise non-negative.
This is entirely analogous to the way in which the non-negativity of the Heron parameters of a triangle assure that the triangle inequalities are satisfied.
In the case of the tetrahedron one also has Yetter's identity (Proposition \myref[1]{thm:yid}) to deal with, but it turns out that this likewise occasions no difficulties.
\begin{proposition} \label{thm:yidok}
The seven areas calculated from any values for the natural parameters via Corollary \myref{thm:algids} always satisfy Yetter's identity.
\end{proposition}
\begin{proof}
Simply use Eqs.~(\myref{eq:idb}) \& (\myref{eq:ide}) to substitute for the squared areas of the exterior \& interior faces in the polynomial $\Xi$ of Yetter's identity (\myref[1]{eq:yetter}) and simplify the result to get the $0$ polynomial.
\end{proof}
Next, a similar process will be used to express the inverse natural parameters as rational functions of the natural parameters.
\begin{lemma} \label{thm:natinvnat}
The inverse natural parameters of a tetrahedron $\GEO{ABCD}$ are given in terms of the natural parameters themselves as
\begin{eqnarray}
\fit[-2em] \tilde u ~= & \!\! \frac{\textstyle 2 \big((v+x)(w+y) - uz \big)}{\textstyle s\hit[2ex]} \,, \qquad \tilde z ~= & \!\! \frac{2 \big( (v+w)(x+y) - uz \big)}{s} \,, \nonumber \\ \label{eq:natinvnat}
\fit[-2em] \tilde v ~= & \!\! \frac{\textstyle2 \big( (u+x)(w+z) - vy \big)}{\textstyle s\hit[2ex]} \,, \qquad \tilde y ~= & \!\! \frac{2 \big( (u+w)(x+z) - vy \big)}{s} \,, \\
\fit[-2em] \tilde w ~= & \!\! \frac{\textstyle2 \big( (u+y)(v+z) - wx \big)}{\textstyle s\hit[2ex]} \,, \qquad \tilde x ~= & \!\! \frac{2 \big( (u+v)(y+z) - wx \big)}{s} \,. \nonumber
\end{eqnarray}
\end{lemma}
\begin{proof}
We will derive only the first of these formulae, since the others may be obtained in much the same fashion.
From Eqs.~(\myref{eq:ida}), (\myref{eq:idb}) \& (\myref{eq:idf}$\,$), we obtain
\begin{equation*}
\tilde u\, s ~=~ 8 \MAG{ABC} \MAG{ABD} \,-\, u\, s ~=~
2\, (u\!+\!v\!+\!x) (u\!+\!w\!+\!y) \,-\, 2\, u\, (u\!+\!v\!+\!w\!+\!x\!+\!y\!+\!z) ~,
\end{equation*}
which is readily verified to be $s$ times the first formula in Eq.~(\myref{eq:natinvnat}).
\end{proof}
\begin{remark}
Note that even when the natural parameters are all non-negative, these formulae can give negative values for one or more of the inverse natural parameters.
Therefore, unlike the Heron parameters of a triangle, the natural parameters of a tetrahedron cannot be chosen arbitrarily subject to being merely non-negative.
\end{remark}

\section{A natural extension of Heron's formula to tetrahedra}
The centerpiece of this paper, and indeed the entire series, may now be stated as follows.
\begin{theorem}[Heron's Formula for Tetrahedra] \label{thm:myform}
With everything defined as above, the volume $\hit[2.5ex] \MAG{ABCD} = t / 3!$ of a tetrahedron $\GEO{ABCD}$ is determined by its natural parameters $\hit[2.5ex] u, v, w, x, y, z$ according to the formula
\pagebreak[2]
\begin{align}
\fit[-1em] t^4 ~=~ & \> s^2\,  \big( 2\, v w x y + 2\, u w x z + 2\, u v y z - u^2 z^2 - v^2 y^2 - w^2 x^2 \big)
\label{eq:myform} \\[3pt] \nonumber
\eqcolon~ & \, s^2\, \Omega(u,v,w,x,y,z) ~=\, -4\: (u\!+\!v\!+\!w\!+\!x\!+\!y\!+\!z)^2\; \mathrm{det\!} \begin{bmatrix}~ 0 & u & v & w~ \\ ~u & 0 & x & y~ \\ ~v & x & 0 & z~ \\ ~w & y & z & 0~ \end{bmatrix} ,
\end{align}
where $s = 2\, (u \!+\! v \!+\! w \!+\! x \!+\! y \!+\! z)$ is twice the exterior surface area, and its in-radius is given by $r^4 =\, t^4 / s^4 =\, \Omega(u,v,w,x,y,z) / s^2\hit[2.5ex]$.
\end{theorem}
\begin{proof}
On dividing the formulae in Eq.~(\myref{eq:r2d2}) through by $r^2$ and substituting for the inverse natural parameters therein using the formulae from Eq.~(\myref{eq:natinvnat}), one obtains
\begin{equation} \fit[-1em] \begin{aligned}
\MAG{AB}^2 ~=~ \, \frac{2\, u\fit \big(  (v\!+\!x) (w\!+\!y) - uz \big)}{s\, r^2} ~,\quad & \MAG{CD}^2 ~=~ \frac{2\, z\fit \big(  (v\!+\!w) (x\!+\!y) - uz \big)}{s\, r^2} ~, \\
\MAG{AC}^2 ~=~ \, \frac{2\, v\fit \big(  (u\!+\!x) (w\!+\!z) - vy \big)}{s\, r^2} ~,\quad & \MAG{BD}^2 ~=~ \frac{2\, y\fit \big(  (u\!+\!w) (x\!+\!z) - vy \big)}{s\, r^2} ~, \\
\MAG{AD}^2 ~=~ \frac{2\, w \big(  (u\!+\!y) (v\!+\!z) - wx \big)}{s\, r^2} ~,\quad & \MAG{BC}^2 ~=~ \frac{2\, x\fit \big(  (u\!+\!v) (y\!+\!z) - wx \big)}{s\, r^2} ~.
\end{aligned} \label{eq:dsq-from-uu2} \end{equation}
Substituting these expressions for the squared distances $D_\msf{ab} \leftrightarrow \MAG{ab}^2$ in the usual $4$-point Cayley-Menger determinant $\CMD_D[\msf A, \msf B, \msf C, \msf D] \leftrightarrow t^2$ and factorizing the result (with the aid of a computer algebra system) then yields
\begin{equation*}
t^2 ~=~ \frac{2\, (u \!+\! v \!+\! w \!+\! x \!+\! y \!+\! z)\, \Omega(u,v,w,x,y,z)^2}{\big( s\, r^2 \big)^{3}} ~.
\end{equation*}
Multiplying this equation through by $s^3r^6$ and using the relation $r = t/s$ thus implies
\begin{equation*}
t^8 / s^3 ~=~ 2\, (u \!+\! v \!+\! w \!+\! x \!+\! y \!+\! z)\,  \Omega(u,v,w,x,y,z)^2  ~.
\end{equation*}
Multiplying through by $s^3$, using the relation $s = 2\, (u\!+\!v\!+\!w\!+\!x\!+\!y\!+\!z) $ from Eq.~(\myref{eq:ida}), and taking the square roots of both sides thus gives Eq.~(\myref{eq:myform}) as desired.
\end{proof}
\begin{remark} \label{rem:mfcontext}
Equation (\myref{eq:myform}) can also be derived by using Eqs.~(\myref{eq:idb}) \& (\myref{eq:ide}) to substitute for the squared areas in the Gramian $\Gamma_{\!F}[\msf A]$, which is a little messier but has the advantage of also being valid in the degenerate case.
On the other hand, if one converts $\Omega(u,v,w,x,y,z)$ into a rational function in the seven areas using Proposition \myref{thm:natpar}, the numerator turns out to be a polynomial of total degree $8$ in the areas containing $420$ terms, which does not factorize.
Since some of these terms contain odd powers of the exterior areas, those areas cannot be eliminated using Yetter's identity, but the interior areas occur in only even powers and hence can be.
If for example one eliminates $\smash{\MAG{AD \tmv BC}}\hit[1.6ex]^2$, the resulting polynomial factorizes into the product of $s^2$ and a polynomial  of total degree $6$ in the remaining six areas containing mere $22$ terms.
(Curiously, this degree $6$ polynomial is equal to $4$ times a four-point Cayley-Menger determinant in the remaining six \underline{areas}, wherein the exterior areas occupy the positions of the edge lengths in a quadrilateral and the two interior areas occupy the positions of its diagonals.)
It turns out that this $22$-term polynomial is the same as that which is obtained on eliminating $\smash{\MAG{AD \tmv BC}}\hit[1.6ex]^2$ from the Gramian of the areal vectors of the interior faces, namely
\begin{equation} \label{eq:intgramat}
\mrm{det}\big( \mbf G_\msf{int} \big) ~\coloneq~ \mrm{det\!} \left[ \begin{smallmatrix} \rule{0pt}{1.5ex}
\big\| \vv{\msf{AB}} \,\Xprod\, \vv{\msf{CD}} \big\|^2 & \big( \vv{\msf{AB}} \,\Xprod\, \vv{\msf{CD}} \big) \,\Dprod\, \big( \vv{\msf{AC}} \,\Xprod\, \vv{\msf{BD}} \big) & \big( \vv{\msf{AB}} \,\Xprod\, \vv{\msf{CD}} \big) \,\Dprod\, \big( \vv{\msf{AD}} \,\Xprod\, \vv{\msf{BC}} \big) \\[2pt] 
\big( \vv{\msf{AB}} \,\Xprod\, \vv{\msf{CD}} \big) \,\Dprod\, \big( \vv{\msf{AC}} \,\Xprod\, \vv{\msf{BD}} \big) & \big\| \vv{\msf{AC}} \,\Xprod\, \vv{\msf{BD}} \big\|^2 & \big( \vv{\msf{AC}} \,\Xprod\, \vv{\msf{BD}} \big) \,\Dprod\, \big( \vv{\msf{AD}} \,\Xprod\, \vv{\msf{BC}} \big) \\[2pt] 
\big( \vv{\msf{AB}} \,\Xprod\, \vv{\msf{CD}} \big) \,\Dprod\, \big( \vv{\msf{AD}} \,\Xprod\, \vv{\msf{BC}} \big) & \big( \vv{\msf{AC}} \,\Xprod\, \vv{\msf{BD}} \big) \,\Dprod\, \big( \vv{\msf{AD}} \,\Xprod\, \vv{\msf{BC}} \big) & \big\| \vv{\msf{AD}} \,\Xprod\, \vv{\msf{BC}} \big\|^2
\rule[-1ex]{0pt}{1ex} \end{smallmatrix} \right] ,
\end{equation}
which may be constructed from the areas using the areal law of cosines for the interior faces given in Eq.~(\myref[1]{eq:cosine3}).
This determinant in turn equals the sum of the exterior Gramians at the four vertices plus \smash{$\big( \smash{\MAG{ABC}}\hit[1.6ex]^2 +  \smash{\MAG{ABD}}\hit[1.6ex]^2 +  \smash{\MAG{ACD}}\hit[1.6ex]^2 +  \smash{\MAG{BCD}}\hit[1.6ex]^2 \big)\!$} $\Xi^2$.
In this way one can convert the polynomial in Eq.~(\myref{eq:myform}) into one that is equivalent modulo $\Xi$ to the Gramians, but at the expense of losing the symmetry under vertex permutations or having to impose Yetter's identity as a constraint (or both).
\end{remark}
This remark and Proposition \myref{thm:yidok} ensure that the areas obtained from the natural parameters via Corollary \myref{thm:ids} fulfill the conditions of Theorem \myref[1]{thm:bdsmc}, which shows that:
\begin{corollary} \label{thm:iff}
Any given $u, v, w, x, y, z \in \mbb R$ are the natural parameters of a non-degenerate Euclidean tetrahedron if \& only if they are positive, the corresponding inverse parameters obtained via Lemma \myref{thm:natinvnat} are positive, and the fourth power of the volume as calculated from Eq.~(\myref{eq:myform}), or equivalently $\Omega(u,v,w,x,y,z)$ by itself, is likewise positive. 
\end{corollary}
\begin{remark} \label{rem:ptolemy}
Under the correspondence $u \leftrightarrow D_\msf{AB\,}$, $\ldots\,$, $z \leftrightarrow D_\msf{CD\,}$, the non-negativity of the negative determinant $\Omega$ in Eq.~(\myref{eq:myform}) is algebraically analogous to the well-known inequality $-\mrm{det}\big( \mbf D \big) \ge 0$ with $\mbf D \coloneq [ D_\msf{ab} ]$.
Factorization of that determinant shows this in turn implies Ptolemy's three inequalities among the distances \smash{$\sqrt{D_\msf{ab}}$} between four points in Euclidean space (see e.g.~Ex.~2 on pg.~80 of Ref.~\cite{Blumenthal:1953}).
In a similar fashion, $\Omega$ factorizes into a product wherein each factor is linear in the products of the square-roots of \textquote{opposite} pairs of natural parameters, i.e.
\begin{equation} \label{eq:ptolemy} \begin{aligned}
\Omega(\hat u^2, \hat v^2, \hat w^2, \hat x^2, \hat y^2, \hat z^2) ~= {}&{} \\
\big( \hat u \hat z + \hat v \hat y + \hat w \hat x \big) {}&{} \big( \hat v \hat y + \hat w \hat x - \hat u \hat z \big) \big( \hat w \hat x + \hat u \hat z - \hat v \hat y \big) \big( \hat u \hat z + \hat v \hat y - \hat w \hat x \big) ~,
\end{aligned} \end{equation}
where $\hat u \coloneq \sqrt{u}\,, \,\ldots,\fit \hat z \coloneq \sqrt{z}\,$.
Nevertheless, even when they determine a non-degenerate Euclidean tetrahedron, the natural parameters are \underline{not} necessarily equal to the squared distances among four Euclidean points, because their square-roots can violate the triangle inequality or give a negative four-point Cayley-Menger determinant (as happens, for example, when $[u,v,w,x,y,z] = [2,4,1,10,5,6]$).
Because Ptolemy's inequalities saturate if \& only if the four points in question lie on a circle in a plane or are collinear, this algebraic analogy provides a way to visualize some of the zeros of $\Omega$, but it does not extend to a geometrically meaningful relationship.
The factorization (\myref{eq:ptolemy}) will nonetheless play a central role in Part III of this series.
\end{remark}
\enlargethispage{0.125in} 
\begin{remark} \label{rem:hyperbolic}
The expression of the quartic polynomial $\Omega(u,v,w,x,y,z)$ as a determinant nonetheless suggests that the formula (\myref{eq:myform}) does have a \EMPH{non-Euclidean} geometric interpretation.
Specifically, it is well known that the matrix of squared distances among a set of points in Euclidean space can be interpreted as the Gram matrix of a set of vectors on the null cone of an indefinite space with signature $[-1, -1, \ldots\,, -1, +1 ]$, and normalized so that their inner product with a fixed null vector, which serves as the point-at-infinity of inversive geometry, is unity; this corresponds to the border of $1$'s in Cayley-Menger determinants \cite{Dress:1993, Seidel:1952, Seidel:1995}.
Although this normalization is not applicable in the present situation, the rest of that geometric interpretation holds, in that the signature of the matrix in Eq.~(\myref{eq:myform}) is $[-1,-1,-1,+1]$.
Because the interior of the null cone, projectively viewed, constitutes a model of hyperbolic space \cite{Coxeter:1998}, it is likely that hyperbolic geometry, and the inversive geometry of its boundary at infinity, will give deeper insights into the meaning of the formula (\myref{eq:myform}).
\end{remark}
\begin{remark} \label{rem:mixforms}
Euler's theorem on homogeneous functions shows that $3\, \CMD_D[\msf A, \msf B, \msf C, \msf D] = \mbf d \,\Dprod\, \nabla_{\mbf d\,} \CMD_D[\msf A, \msf B, \msf C, \msf D]$, where $\mbf d \coloneq [D_\msf{AB\,},\, \ldots,\, D_\msf{CD}]$ is a vector of squared distances.
The derivatives of the determinant in $\nabla_{\mbf d\,} \CMD_D\,$, in turn, are the cofactors of the corresponding matrix, which Lemma \myref[1]{thm:aloc} shows are the dot products of the areal vectors of pairs of exterior faces.
It thus follows from Corollary \myref{thm:r2d2}, $r = t/s$ and Eq.~(\myref{eq:idg}) that
\begin{align}
\fit[-1em] t^2 ~=~ \tfrac13\, \mbf d \,\Dprod\, \nabla_{\mbf d\,} \CMD_D[\msf A, \msf B, \msf C, \msf D] ~=~ & \frac{s^3}{6\, t^2}\, \big( u\tilde u\, (z - \tilde z) + v\tilde v\, (y - \tilde y) + w\tilde w\, (x - \tilde x) \\
& \fit[2em] +\fit x\tilde x\, (w - \tilde w) + y\tilde y\, (v - \tilde v) + z\tilde z\, (u - \tilde u) \big) ~. \nonumber
\end{align}
From this, one obtains another formula for $t^4$ that is (outside of $s$) antisymmetric w.r.t.~interchange of the natural and inverse natural parameters.
On substituting for the inverse parameters using Lemma \myref{thm:natinvnat}, one again arrives at Eq.~(\myref{eq:myform}).

Alternatively, one can write the off-diagonal entries of $\mbf G_\msf A$ in terms of $u,\tilde u, v, \tilde v, w, \tilde w$ \& $s$ using Eq.~(\myref{eq:idg}), and its diagonal entries as $8\, \smash{\MAG{ABC}}^2 = s\, (u + \tilde u) (v + \tilde v) \fit/\fit (w + \tilde w)$, $8\, \smash{\MAG{ABD}}^2 = s\, (u + \tilde u) (w + \tilde w) \fit/\fit (v + \tilde v)$, $8\, \smash{\MAG{ACD}}^2 = s\, (v + \tilde v) (w + \tilde w) \fit/\fit (u + \tilde u)$ by Eq.~(\myref{eq:idf}$\,$).
The Gramian then becomes a rational function with a numerator which factorizes into a product of four cubic factors in the square-roots of those natural and of the corresponding inverse natural parameters $\check u \coloneq \sqrt{\tilde u}$, $\check v \coloneq \sqrt{\tilde v}$, $\check w \coloneq \sqrt{\tilde w}$, thereby showing that:
\begin{align}
& t^4 ~=~ \frac{s^3}{2\, (u + \tilde u)(v + \tilde v)(w + \tilde w)}\, \big( -\hat u \hat v \hat w + \check u \check v \hat w + \check u \hat v \check w + \hat u \check v \check w \big) \,\cdots \\
& \cdots\, \big( \hat u \hat v \hat w - \check u \check v \hat w + \check u \hat v \check w + \hat u \check v \check w \big) \big( \hat u \hat v \hat w + \check u \check v \hat w - \check u \hat v \check w + \hat u \check v \check w \big) \big( \hat u \hat v \hat w + \check u \check v \hat w + \check u \hat v \check w - \hat u \check v \check w \big) \nonumber
\end{align}
Analogous expressions can of course be derived from the Gramians at the other three vertices.
\end{remark}
\begin{remark} \label{rem:xfactor}
A rather different expression which also relates the natural parameters to the volume of the tetrahedron is
\begin{equation}
s\, (u - z)\fit (v - y)\fit (w - x)\fit (u \!+\! v \!+\! w)\fit (u \!+\! x \!+\! y)\fit (v \!+\! x \!+\! z)\fit (w \!+\! y \!+\! z)
\end{equation}
(cf.~Eqs.~(\myref{eq:ida}), (\myref{eq:idc}) \& (\myref{eq:idd})).
Unlike Eq.~(\myref{eq:myform}) this can be expanded into a polynomial in only the squared exterior areas, and hence may be written as a polynomial in the squared distances by replacing them by the corresponding three-point Cayley-Menger determinants (cf.~Eq.~(\myref[1]{eq:cmd3})).
This later polynomial factorizes into a product of the four-point Cayley-Menger determinant and another factor, dubbed the \textquote{$X$-factor,} of total degree $5$ in the squared distances.
By construction the $X$-factor vanishes whenever any one exterior area equals the sum of the other three, or the sum of any two exterior areas equals the sum of the other two (and hence, in particular, for equi-facial tetrahedra), but it is not necessarily non-negative even in the Euclidean case, and its full geometric interpretation remains an open problem.
\end{remark}

The determinantal form of Eq.~(\myref{eq:myform}) immediately suggests a further extension to Euclidean spaces of dimension $n > 3$, as well.
Clearly for a Euclidean $n$-simplex the analogues of the contact triangles of a tetrahedron are the \textquote{contact $(n-1)$-simplices} into which its facets are divided by their respective in-touch points.
There are $n\fit(n+1)$ of these and, via the Pythagorean theorem, it is easily shown that they again come in congruent pairs.
Taking $(n-1)!$ times the \textquote{hyper-areas} of these pairs of $(n-1)$-simplices as the natural parameters of the $n$-simplex then leads to the following:
\begin{conjecture} \label{thm:alldim}
The hyper-volume of an $n$-simplex $\OL{\msf{ABC} \cdots}$ is given in terms of its\linebreak $(n+1)\,n/2$ natural parameters  $u, v, \ldots\,, w\fit, x\fit, \ldots\, y\fit, \ldots\, z\fit, \ldots$ by
\begin{equation} \begin{aligned}
\big(\fit n!\, \big|\, \GEO{ABC \,\cdots\hit[1.7ex]} \,\big| {}&{} \fit\big)^{\!2(n-1)} ~= \\
&{} {(-1)\hit[1.8ex]}^n\, \big( 2\, (u + v + \cdots + z + \cdots ) \big)^{\!n-1}\; \mathrm{det\!} \begin{bmatrix} ~0 & u & v & \cdots & w~ \\ ~u & 0 & x &\! \cdots \!& y~ \\ ~v & x & 0 &\! \cdots \!& z~ \\[-3pt] \vdots\!\! & \vdots & \vdots &\! \ddots \!& \!\!\vdots \\ ~w & y & z &\! \cdots \!& 0~ \end{bmatrix} .
\end{aligned} \label{eq:conjecture} \end{equation}
\end{conjecture}
\noindent Note the linear factor on the right contributes $(n-1)^2$ to the dimensionality, while the determinantal factor adds another $(n-1) (n+1)$, matching the total of $2\, n\, (n-1)$ on the left. 
\begin{remark}
By computing the relevant quantities in multiple random $4$-simplices, it has been numerically confirmed that the volumes of the $20$ contact tetrahedra into which the in-touch points divide their facets are equal when they share a $2$-face, and that putting these $10$ numbers into Eq.~(\myref{eq:conjecture}) does indeed give the sixth power of $4!$ times their hyper-volumes.
\textit{Note added in proof:} An essentially geometric proof of this conjecture has just been published by the author \cite{Havel:2024}.
\end{remark}

\paragraph*{\indent\textit{Note added in proof but deleted by editor:}} The author has lately come across a question posted on the \texttt{StackExchange} forum's \texttt{Mathematics} site, dating from April 2021, where a user asked if the contact tri\-angle areas of a tetrahedron determined its volume (see \url{https://math.stackexchange.com/questions/4119104/volume-of-a-tetrahedron-given-}
\url{areas-of-6-triangles}).
In response, Billy Don Sterling McConnell (under his pseudonym \textquote{Blue}; cf.~Theorem \myref[1]{thm:bdsmc}) outlined an elementary though also fairly algebraically demanding derivation of the formula (\myref{eq:myform}), which he expressed using the factorization of $\Omega$ in Eq.~(\myref{eq:ptolemy}).
In this paper's notation, his derivation uses Eq.~(\myref{eq:idb}) and the relations $u = 4 \MAG{ABC}\! \MAG{ABD} (1 +\fit[0.1em] \cos(\varphi_\msf{AB})) \fit/\fit s$ etc.~from Eqs.~(\myref{eq:idf}) \& (\myref{eq:idg}) to rewrite the Gramian $\Gamma_{\!F}[\msf A]$ in terms of the contact triangle areas (aka natural parameters).
This makes him an essentially simultaneous co-discoverer of the formula (\myref{eq:myform}).
There is no indication, however, that he realized he had found a geometrically natural extension of Heron's formula to tetrahedra at that time.

\smallskip
\pagebreak[3]
{\small

} 


\section*{\indent Biography}
The author obtained his doctoral degree in Biophysics from the University of California at Berkeley in 1982, where he co-developed the \textquote{distance geometry} approach to the analysis and computation of molecular conformation with G.~M. Crippen \& I.~D. Kuntz (both then at U.~C.~San Francisco).
He then embarked upon a long and checkered career that took him back and forth across the Atlantic and through many nominally unrelated fields of mathematics, science and engineering.
He has written widely used software for computing protein structures from nuclear magnetic resonance (NMR) data (with Prof.~Kurt W\"uthrich \& Gerhard Wagner at the Swiss Federal Technical Institute in Z\"urich) as well as from multiple sequence alignments with homologues of known structure (at the Research Institute of Scripps Clinic in La Jolla, California and at the Univ.~of Michigan in Ann Arbor), work which ultimately contributed to the so-called solution of the Protein Folding Problem \cite{Marx:2022}.
He has also co-authored a book on the theory and applications of distance geometry to molecular conformation with Prof.~Crippen \cite{Crippen:1988}, and written papers on the combinatorial structure of metric and metric vector spaces with Prof.~Andreas W.~M.~Dress at the Univ.~of Bielefeld in Germany \cite{Dress:1986}.
He settled in Boston, Massachusetts in 1990, where he spent a number of years developing matrix methods for the analysis of biomolecular NMR data in the Dept.~of Biological Chemistry and Molecular Pharmacology the Harvard Medical School.
During this time he also began dabbling in geometric (aka Clifford) algebras and their applications to distance geometry \cite{Dress:1993} as well as mathematical physics more generally.
This latter interest (and his funding) led him to move to the Dept.~of Nuclear Science and Engineering at MIT around the turn of the millennium, where he worked with Prof.~David G.~Cory on applying geometric algebra techniques to the design and analysis of experiments demonstrating the principles of quantum computing by NMR spectroscopy.
After a stint attempting to commercialize thermochemical and nano-mechanical forms of energy storage, he unofficially retired circa 2016 and devoted much of his time to learning about neuroscience, artificial intelligence and applied category theory in the rich intellectual environment of the greater Boston area.
This series of papers is an outcome of his self-imposed isolation during the Coronavirus Pandemic of 2020-22.

\begin{displayquote}\medskip\centering
This series of papers is dedicated to the memory of my doctoral\\ co-advisor, Gordon M.~Crippen (1945--2022).
\end{displayquote}

\begin{center}
\medskip\rule{70pt}{2pt}\bigskip
\end{center}

%
%

\begin{displayquote}\small
If this [the Mysterium cosmographicum] is published, others will perhaps make discoveries I might have reserved for myself.
But we are all ephemeral creatures (and none more so than I).
I have, therefore, for the Glory of God, who wants to be recognized from the book of Nature, that these things may be published as quickly as possible.
The more others build on my work the happier I shall be. 

\raggedleft\textit{Johannes Kepler, 1595}
\end{displayquote}

\bigskip
\begin{displayquote}\small
Every formula which expresses a law of nature is a hymn of praise to God.

\smallskip\raggedleft\textit{Maria Mitchell, date unknown}
\end{displayquote}

\mypart{\mytitle} 

\begin{abstract}
This is the third part of a series of four papers in The $\Pi\mrm{ME}$ Journal.
It explores the structure of the zeros of the extension of Heron's formula to tetrahedra that was derived in Part II.
First, it is shown that almost all the zeros, specifically those for which the areal Gram matrices of Part I have rank $2$, are the limits of sequences of non-degenerate tetrahedra the vertices of which go off to infinity along a line while the areas of their interior and exterior faces remain finite.
Zeros with areal Gram matrices of rank $1$, in contrast, correspond to quadruples of points in the (finite) affine plane.
The full five-dimensional zero set is then shown to be canonically homeomorphic to a certain quotient of the well-known Klein quadric under a (generically) faithful action of the group $\mbb Z_2^4$.
The zero set also admits a natural stratification according to which $2\times 2$ minors of the areal Gram matrices vanish, and this gives it the structure of a finite graded lattice with the set of generic degenerate tetrahedra as its supremum and the set of all quadruples in the affine plane as its infimum. 
Finally, it is shown that the algebraic structure of the zeros in the affine plane naturally defines the associated four-element, rank $3$ chirotope, aka affine oriented matroid.
The paper closes with some remarks on the potential significance of these results in mathematics and physics, along with possible directions for future research.
\end{abstract}

\section{The projective nature of the zeros of $\boldsymbol\Omega$} \label{sec:zeros}
In what follows, a tetrahedron will be called \textquote{degenerate} whenever the polynomial $\Omega$ from the extension of Heron's formula to tetrahedra given in Theorem \myref[2]{thm:myform}\footnote{A Roman numeral followed by a colon specifies a reference to a numbered entity in another paper of this series, in this case Theorem \myref{thm:myform} of Part II.} vanishes.
Equivalently, a tetrahedron is degenerate when the determinant of the exterior areal Gram matrix from Part I at any (and hence all) of its vertices is zero.
It is geometrically clear that the areal Gram matrices of any quadruple of points in the Euclidean plane will have a rank of $1$, since the areal vectors of the seven faces in any such configuration are of course all collinear.
The infinitely more common class of zeros for which these Gram matrices have a rank of $2$, however, does \EMPH{not} correspond to planar configurations, nor to any other configuration of points heretofore considered in classical Euclidean geometry.

In order to gain some insight into what these are, we may express the squared distances in a three-point Cayley-Menger determinant in terms of the natural parameters just as was done with the four-point determinant in the proof of Theorem \myref[2]{thm:myform} (specifically, Eq.~(\myref[2]{eq:dsq-from-uu2})), obtaining e.g.~$\CMD_D[\msf A, \msf B, \msf C] \leftrightarrow$
\begin{align}
&\> -\!\sfrac14\, \mrm{det\!} \begin{bmatrix} ~0&1&1&1 \\ ~1 & 0 & u\tilde u / r^2 & v\tilde v / r^2 \\ ~1 & u\tilde u / r^2 & 0 & x\tilde x / r^2 \\ ~1 & v\tilde v / r^2 & x\tilde x / r^2 & 0 \end{bmatrix} ~=~ \tfrac14\, r^{-4}\, \Omega(u,v,x,\tilde x,\tilde v,\tilde u)
\\ =~ &\>
r^{-4}\, (u + v + x)^2\, \Omega(u,v,w,x,y,z)\, /\fit s^2 ~=~ (u + v + x)^2 ~=~ 4 \MAG{ABC}^2 \,, \nonumber
\end{align}
where $r = t/s$ is the in-radius as in Part II.
On multiplying through by $r^4$, we see that the three-point Cayley-Menger determinant $\Omega(u,v,x,\tilde x,\tilde v,\tilde u)$ in the \TDEF{complementary products} $u\tilde u,\, v\tilde v,\, x\tilde x$ vanishes in the limit of a degenerate tetrahedron for which $r^4 = \Omega / s^2 = 0$.
In fact all four such Cayley-Menger determinants $\CMD_{\widetilde D}[\msf a, \msf b, \msf c]$, wherein $\widetilde D_\msf{ab} \leftrightarrow r^2 D_\msf{ab}$ ($\msf a, \msf b, \msf c \in \{ \msf A, \msf B, \msf C, \msf D \}$ with $\msf a \ne \msf b \ne \msf c \ne \msf a$), vanish identically in that limit.
These three-point determinants correspond to
\begin{equation} \begin{alignedat}{2}
\Omega(u,v,x,\tilde x,\tilde v,\tilde u\,) / 4 ~=~ & 4\, r^4 \MAG{ABC}^2 ,\quad \Omega(u,w,y,\tilde y,\tilde w,\tilde u) / 4 ~=~ & 4\, r^4 \MAG{ABD}^2 , \\
\Omega(v,w,z,\tilde z,\tilde w,\tilde v) / 4 ~=~ & 4\, r^4 \MAG{ACD}^2 ,\quad \Omega(x,\,y,\,z,\,\tilde z,\,\tilde y,\,\tilde x) / 4 ~=~ & 4\, r^4 \MAG{BCD}^2 .
\end{alignedat} \end{equation}
Because the four-point Cayley-Menger determinant in the complementary products,\linebreak[2] $\CMD_{\widetilde D}[\msf A, \msf B, \msf C, \msf D] \leftrightarrow r^6\fit t^2$, also vanishes in this limit, it follows that in any degenerate tetrahedron these complementary products are (formally) the squared distances \smash{$\widetilde{D}$} among four points on a Euclidean line.
The tetrahedron's actual squared inter-vertex distances, however, diverge towards infinity as $r^{-2}$ times the corresponding complementary product, providing that product does not itself go to zero.
Thus these degenerate tetrahedra can be said to have collinear vertices separated by infinite distances, but with generically well-defined ratios $\MAG{AB}\raisebox{\upbgdl}{$/$}\MAG{AC} = u\tilde u/v\tilde v$ etc., and with interior \& exterior faces of finite, and generally non-zero, area.

One way to construct such zeros is to take a random non-degenerate tetrahedron in $\mbb R^3$, apply an affine transformation $\mathcal A_\sigma$ with diagonal matrix $\mbf{Diag}(\sigma^{-1},\, \sigma^{-1},\, \sigma)$ for some  $\sigma > 0$ to its vertices, and take the limit as $\sigma \rightarrow \infty$.
It is easily seen that these\linebreak[2] transformations act asymptotically on the (interior and exterior) faces as a two-dimen\-sional affine squeeze of the form $\mbf{Diag}(\sigma^{-1},\, \sigma)$.
Since such a squeeze preserves areas, it follows that as $\sigma \rightarrow \infty$ the seven areas will converge to well-defined finite values while the inter-vertex distances approach infinity and the volume goes to zero.
On putting the resulting limits of the areas into Eq.~(\myref[1]{eq:gramatA}), one obtains a Gram matrix at $\GEO{A}$ (or any other vertex) of rank $2$.
Areas corresponding to rank $1$ zeros, in contrast, are readily obtained by projecting the tetrahedron onto any plane, although these areas obviously determine the corresponding planar configuration only up to special (area preserving) affine transformations.
If the plane of projection is parallel to the third coordinate axis, $\mathcal A_\sigma$ acts on it as an area-preserving affine squeeze for all $\sigma > 0$ which sends the vertices to infinity as $\sigma \rightarrow \infty$, thereby showing how the rank $1$ and $2$ zeros of $\Omega$ are related.

These degenerate tetrahedra cannot be viewed as simply a quadruple of points on a line in\linebreak[2] the projective completion of Euclidean three-space, because it is possible for \underline{all} the vertices in such configurations to be at infinite distances from each other whereas a line in that completion has only \underline{one} point at infinity.
The proper interpretation of these un\-conventional Euclidean configurations within the framework of projective geometry, reaffirming Arthur Cayley's claim that \textquote{projective geometry is all geometry} \cite{ACampo:2015}, will be left as a challenge to the experts in that field (Refs.~\cite{Babson:2002, Babson:2006} might be a good place to\linebreak[2] start).
Instead, this paper will seek to motivate further study of such questions, by show\-ing that the zeros of $\Omega$ can be placed in a one-to-one correspondence with a certain quotient of the Klein quadric by an action of a discrete group of reflections on the\linebreak[2] Pl\"ucker coordinates.
It will further explore the combinatorial structure imposed on the\linebreak[2] set of all degenerate tetrahedra by the various possible combinations of vanishing com\-plementary products, and show (see also Appendix \myref[4]{sec:planar}) that the aforementioned rank $1$\linebreak[2] zeros in the affine plane are exactly those wherein all six complementary products vanish.

We end this section with a technical lemma which is needed to achieve these goals.
\begin{lemma} \label{thm:techlem}
Given any $u, v, w, x, y, z \in \mbb R$ with $s = 2\, (u\!+\!v\!+\!w\!+\!x\!+\!y\!+\!z) \ne 0$, let\linebreak[2] $\tilde u$, $\tilde v$, $\tilde w$, $\tilde x$, $\tilde y$, $\tilde z$ be the values of the rational functions in Eq.~(\myref[2]{eq:natinvnat}).
Then if any one of the complementary products $u\tilde u$, $v\tilde v$, $w\tilde w$, $x\tilde x$, $y\tilde y$ or $z\tilde z$ vanishes, the polynomial $\Omega$ in Eq.~(\myref[2]{eq:myform}) satisfies $\Omega(u,v,w,x,y,z) \le 0$. 
Conversely, given $u, v, w, x, y, z \ge 0$ with $\Omega(u,v,w,x,y,z)$ $\,=\,$ $0$, the rational functions' values $\tilde u, \tilde v, \tilde w, \tilde x, \tilde y, \tilde z$ from Eq.~(\myref[2]{eq:natinvnat}) are all non-negative, as are the corresponding complementary products.
\end{lemma}
\begin{proof}
To prove the first claim, suppose for example $u\tilde u = 0$ so that either $u = 0$ or $\tilde u = 0$.
In the former case we find that $\Omega|_{u = 0} = -\fit(w x - v y)^2 \le 0$, whereas if $u \ne 0$ we may solve $\tilde u = 2\, \big( (v+x)(w+y) - uz \big) / s = 0$ in Eq.~(\myref[2]{eq:natinvnat}) for $z^* = (v + x)(w + y) / u$, whence $\Omega|_{z=z^*} = -(v w - x y)^2 \le 0$ as well.
The proof if any of the other products vanish is of course similar.

To prove the second claim, we solve $\Omega(u,v,w,x,y,z) = 0$ for the product $u z$, obtaining (with $\hat v \coloneq \sqrt{v}$, $\hat w \coloneq \sqrt{w}$, etc.)
\begin{equation}
u z ~=~ w x + v y \,\pm\, 2\, \sqrt{v w x y} ~=~ \big( \hat w \hat x \,\pm\, \hat v \hat y \big)^{\!2} ~.
\end{equation}
If $s = 0$ the claim holds vacuously, and otherwise substituting this value of $u z$ into Eq.~(\myref[2]{eq:natinvnat}) for $\tilde u$ yields
\begin{equation}
\tilde u ~=~ 2\, \big( v w + x y \,\pm\, 2\, \sqrt{v w x y} \big)  / s ~=~ 2\, \big( \hat v \hat w \,\pm\, \hat x \hat y \big)^{\!2} / s ~\ge~ 0 ~,
\end{equation}
as desired.
The proof for the remaining five products is again similar.
\end{proof}
\begin{remark}
The second claim of this lemma suggests that $u, v\fit[-1pt], w\fit[-1pt], x, y\fit[-1pt], z \ge 0$ and $\Omega(u, v\fit[-1pt], w\fit[-1pt], x, y\fit[-1pt], z)$ $\ge$ $0$ implies $\tilde u, \tilde v\fit[-1pt], \tilde w\fit[-1pt], \tilde x, \tilde y\fit[-1pt], \tilde z \ge 0$.
Assuming these inequalities are strict, this strengthening of Corollary \myref[2]{thm:iff} can be proven by noting that $\Omega > 0$ implies the Gram matrix $\mbf G_\msf A$ at $\GEO{A}$ (say) has either zero or two negative eigenvalues (since $\mrm{det}\big( \mbf G_\msf A \big) = t^4 = s^2 \Omega > 0$).
Because the inverse natural parameters are all positive in the former case, we need only prove our claim in the latter.
To do so, first note that the diagonal entries of $\mbf G_\msf A$ are the squared exterior areas $4 \smash{\MAG{ABC}}^2 = (u + v + x)^2$, $4 \smash{\MAG{ABD}}^2 = (u + w + y)^2$, $4 \smash{\MAG{ACD}}^2 = (v + w + z)^2$ by Eq.~(\myref[2]{eq:idb}) and hence positive.
This in turn implies all three $2\times2$ principal minors $s^2 u\tilde u$, $s^2 v\tilde v$, $s^2 w\tilde w$ (see the note immediately preceding Corollary \myref[2]{thm:ids}) must be negative, since otherwise $\mbf G_\msf A$ would be positive definite by Sylvester's criterion.
It follows that $u, v, w > 0$ implies $\tilde u, \tilde v, \tilde w < 0$ and hence $\tilde u + \tilde v + \tilde w < 0$.
Upon substituting for $\tilde u, \tilde v$ \& $\tilde w$ in that sum using Lemma \myref[2]{thm:natinvnat}, however, we obtain the manifest contradiction $2 (uv\fit[-0.07em]+\fit[-0.07em]uw\fit[-0.07em]+\fit[-0.07em]vw\fit[-0.07em]+\fit[-0.07em]uz\fit[-0.07em]+\fit[-0.07em]vy\fit[-0.07em]+\fit[-0.07em]wx\fit[-0.07em]+\fit[-0.07em]xy\fit[-0.07em]+\fit[-0.07em]xz\fit[-0.07em]+\fit[-0.07em]yz) / s < 0$.
\end{remark}

\section{The connection to the Klein quadric} \label{sec:klein}
The analysis of the zeros of $\Omega$ that follows benefited greatly from an exposition of a connection between planar polygons and Grassmannians recently given by Cantarella \TDEF{et al.\/} \cite{Cantarella:2019} (see also Ref.~\cite{Eastwood:2000}).
Specifically, when $\mrm{rank}\big( \mbf G_\msf A \big) \le 2$ the areal vectors of the exterior faces of a tetrahedron are coplanar, and when signed correctly sum to zero by Minkowski's iden\-tity (\myref[1]{eq:minkow}) just like the edge vectors of a planar quadrilateral (this analogy was also briefly considered for non-degenerate simplices in Ref.~\cite{Khimshiashvili:2013}).
Thus it is possible to visualize these degenerate tetrahedra as planar quadrilaterals, though it should be noted that the vertices of such quadrilaterals are not those of the tetrahedra themselves, and that the cyclic order of their edges is arbitrary.

The approach used by Cantarella \TDEF{et al.\/} equates the components of four such vectors in $\mbb R^2$ to the real \& imag\-inary parts of the \EMPH{squares} of four complex numbers $m_\msf A \!+\!\fit i n_\msf A \fit$, $\ldots,\fit$ $m_\msf D \!+\!\fit i n_\msf D$, so that
\begin{equation} \label{eq:r2avec}
\fit[-0.5em] \begin{aligned}
\VCP{AB}{AC} \,=\, (m_\msf D^2 - n_\msf D^2)\, \mbf e_1 \!\Vplus 2\fit m_\msf D n_\msf D\, \mbf e_2 \,,{}&{}\, \Vdiff \VCP{AB}{AD} \,=\, (m_\msf C^2 - n_\msf C^2)\, \mbf e_1 \!\Vplus 2\fit m_\msf C n_\msf C\, \mbf e_2 \,, \\
\VCP{AC}{AD} \,=\, (m_\msf B^2 - n_\msf B^2)\, \mbf e_1 \!\Vplus 2\fit m_\msf B n_\msf B\, \mbf e_2 \,,{}&{}\, \Vdiff \VCP{BC}{BD} \,=\, (m_\msf A^2 - n_\msf A^2)\, \mbf e_1 \!\Vplus 2\fit m_\msf A n_\msf A\, \mbf e_2  \,,
\end{aligned}
\end{equation}
where $\mbf e_1$, $\mbf e_2$ is any orthonormal basis of their common plane (in practice, the basis obtained by diagonalizing $\mbf G_\msf A$ is a convenient one to use).
It follows that the norms of these cross products are given quite simply by
\begin{equation} \label{eq:degextareas}
m_\msf A^2 \!+\!\fit n_\msf A^2 = 2 \MAG{BCD} ,\fit[0.35em] m_\msf B^2 \!+\!\fit n_\msf B^2 = 2 \MAG{ACD} ,\fit[0.35em] m_\msf C^2 \!+\!\fit n_\msf C^2 = 2 \MAG{ABD} ,\fit[0.35em] m_\msf D^2 \!+\!\fit n_\msf D^2 = 2 \MAG{ABC} .
\end{equation}
The fact that the vectors in Eq.~(\myref{eq:r2avec}) sum to $\mbf 0$ further implies that the vectors $\mbf m$ $\coloneq$ \linebreak[2] $[m_\msf A$, $m_\msf B$, $m_\msf C$, $m_\msf D]^{\smash{\top}}\!$, $\mbf n$ $\coloneq$ $[n_\msf A$, $n_\msf B$, $n_\msf C$, $n_\msf D]^{\smash{\top}}\!$ in $\mbb R^4$ satisfy $\mbf m \Dprod \mbf n = 0$ \& $\| \mbf m \|^2 \!= \| \mbf n \|^2 \!= s/2$, where $s$ is twice the exterior surface area as usual.
These vectors are determined only up to improper rotations in the plane they span, and a rotation by an angle $\vartheta$ in their common plane corresponds to a rotation of the areal vectors by $2\fit\vartheta$.

We will now derive two equivalent formulae for each of the three interior facial areas in terms of the real and imaginary parts of these complex numbers.
\begin{lemma} \label{thm:degintareas}
With everything defined as above:
\begin{subequations} \begin{align}
16\, \MAG{AB|CD}^2 ~=~ &\>
\big( (m_\msf B + n_\msf A)^2 +\, (m_\msf A - n_\msf B)^2 \big)\, \big( (m_\msf B - n_\msf A)^2 +\, (m_\msf A + n_\msf B)^2 \big) \\[-1pt] =~ &\>
\big( (m_\msf D + n_\msf C)^2 +\, (m_\msf C - n_\msf D)^2 \big)\, \big( (m_\msf D - n_\msf C)^2 +\, (m_\msf C + n_\msf D)^2 \big) ; \nonumber \\[1pt]
16\, \MAG{AC|BD}^2 ~=~ &\>
\big( (m_\msf D + n_\msf B)^2 +\, (m_\msf B - n_\msf D)^2 \big)\, \big( (m_\msf D - n_\msf B)^2 +\, (m_\msf B + n_\msf D)^2 \big) \\[-1pt] =~ &\>
\big( (m_\msf C + n_\msf A)^2 +\, (m_\msf A - n_\msf C)^2 \big)\, \big( (m_\msf C - n_\msf A)^2 +\, (m_\msf A + n_\msf C)^2 \big) ; \nonumber \\[1pt]
16\, \MAG{AD|BC}^2 ~=~ &\>
\big( (m_\msf D + n_\msf A)^2 +\, (m_\msf A - n_\msf D)^2 \big)\, \big( (m_\msf D - n_\msf A)^2 +\, (m_\msf A + n_\msf D)^2 \big) \\[-1pt] =~ &\>
\big( (m_\msf C + n_\msf B)^2 +\, (m_\msf B - n_\msf C)^2 \big)\, \big( (m_\msf C - n_\msf B)^2 +\, (m_\msf B + n_\msf C)^2 \big) . \nonumber
\end{align} \label{eq:degintareas} \end{subequations}
\end{lemma}
\vspace{-5ex} 
\begin{proof}
To prove the first of the above formulae, we use Eq.~(\myref{eq:r2avec}) and the usual expression for the dot product of vectors in terms of their coordinates to obtain
\begin{equation*}
-\fit \big( \VCP{AC}{AD} \big) \Dprod\fit \big( \VCP{BC}{BD} \big) ~=~
(m_\msf A^2 - n_\msf A^2)\, (m_\msf B^2 - n_\msf B^2) \,+\, 4\, m_\msf A\fit n_\msf A\fit m_\msf B\fit n_\msf B ~.
\end{equation*}
By the areal law of cosines (\myref[1]{eq:cosines}) together with Eq.~(\myref{eq:degextareas}), however, $16 \MAG{AB|CD}^2 \,=$
\begin{align*}
\NCP{AB}{CD}^2 ~=~ &\>
\NCP{AC}{AD}^2 +\, \NCP{BC}{BD}^2 -\, 2\, \big( \VCP{AC}{AD} \big) \Dprod \big( \VCP{BC}{BD} \big) \\ \nonumber =~ &\>
\big( m_\msf A^2 + n_\msf A^2 \big)^{\!2} \!+ \big( m_\msf B^2 + n_\msf B^2 \big)^{\!2} \!+ 2\, (m_\msf A^2 - n_\msf A^2)\, (m_\msf B^2 - n_\msf B^2) + 8\, m_\msf A n_\msf A m_\msf B\fit n_\msf B \\ =~ &\> \nonumber
\big( (m_\msf B + n_\msf A)^2 +\, (m_\msf A - n_\msf B)^2 \big)\, \big( (m_\msf B - n_\msf A)^2 +\, (m_\msf A + n_\msf B)^2 \big) ~.
\end{align*}
A similar procedure, applied to the alternative expression $\NCP{AB}{CD}^2$ $\!=$ $\NCP{AB}{AC}^2$ $+\,$ $\NCP{AB}{AD}^2$ $-\,$ $2\, \big( \VCP{AB}{AC} \big) \Dprod \big( \VCP{AB}{AD} \big)$, yields the second formula for $16 \MAG{AB \tmv CD}^2$.
The re\-main\-ing formulae can be established in an analogous fashion.\hit
\end{proof}
These results allow us to express the natural and inverse natural parameters of a degenerate tetrahedron in terms of the $m$'s \& $n$'s quite simply as follows.
\begin{proposition} \label{thm:pc-to-np}
Given a tetrahedron with volume $\MAG{ABCD} = 0$ and exterior surface area (times $2$) $s > 0$, together with vectors $\mbf m, \mbf n \in \mbb R^4$ as above, the natural parameters are given by
\begin{equation} \begin{aligned}
u ~=~ 2\, (m_\msf C n_\msf D - m_\msf D n_\msf C)^2 / s ~,&\qquad
z ~=~ 2\, (m_\msf A n_\msf B - m_\msf B n_\msf A)^2 / s ~,\\
v ~=~ 2\, (m_\msf B n_\msf D - m_\msf D n_\msf B)^2 / s ~,&\qquad
y ~=~ 2\, (m_\msf A n_\msf C - m_\msf C n_\msf A)^2 / s ~,\\
w ~=~ 2\, (m_\msf B n_\msf C - m_\msf C n_\msf B)^2 / s ~,&\qquad
x ~=~ 2\, (m_\msf A n_\msf D - m_\msf D n_\msf A)^2 / s ~,
\end{aligned} \label{eq:mn2natpar} \end{equation}
while the inverse natural parameters are given by
\begin{equation} \begin{aligned}
\tilde u ~=~ 2\, (m_\msf C m_\msf D + n_\msf D n_\msf C)^2 / s ~,&\qquad
\tilde z ~=~ 2\, (m_\msf A m_\msf B + n_\msf B n_\msf A)^2 / s ~,\\
\tilde v ~=~ 2\, (m_\msf B m_\msf D + n_\msf D n_\msf B)^2 / s ~,&\qquad
\tilde y ~=~ 2\, (m_\msf A m_\msf C + n_\msf C n_\msf A)^2 / s ~,\\
\tilde w ~=~ 2\, (m_\msf B m_\msf C + n_\msf C n_\msf B)^2 / s ~,&\qquad
\tilde x ~=~ 2\, (m_\msf A m_\msf D + n_\msf D n_\msf A)^2 / s ~.
\end{aligned} \label{eq:mn2invpar} \end{equation}
If $s = 0$, of course, the natural and inverse natural parameters are all zero as well.
\end{proposition}
\begin{proof}
Upon substituting for the exterior \& interior areas in the expression for $u$ from Proposition \myref[2]{thm:natpar} using Eqs.~(\myref{eq:degextareas}) \& (\myref{eq:degintareas}), we obtain
\begin{align*}
& 2 su ~=~ \Big( 2 \MAG{ABC} + 2 \MAG{ABD} \Big)^{\!2} -\, 16 \MAG{AB|CD}^2 ~= \\[-2pt]
&\> \begin{aligned}[t] \big( m_\msf D^2 + n_\msf D^2 + m_\msf C^2 + n_\msf C^2 \big)^2 -\, \big( (m_\msf D + n_\msf C)^2 + (m_\msf C - n_\msf D)^2 \big) \big( (m_\msf D - n_\msf C)^2 + (m_\msf C + n_\msf D)^2 \big) & \\
=\> 4\, (m_\msf C n_\msf D - m_\msf D n_\msf C)^2 & ~, \end{aligned} \nonumber
\end{align*}
as desired.
Similarly, upon substituting for the areas in the expression for $\tilde u$ from Definition \myref[2]{def:invnatpar} (Eq.~(\myref[2]{eq:invpar})), we obtain
\begin{align*}
& 2 s\tilde u ~=~ 16 \MAG{AB|CD}^2 -\, \Big( 2 \MAG{ABC} - 2 \MAG{ABD} \Big)^{\!2} ~= \\[-2pt]
&\> \begin{aligned}[t]
\big( (m_\msf D + n_\msf C)^2 + (m_\msf C - n_\msf D)^2 \big) \big( (m_\msf D - n_\msf C)^2 + (m_\msf C + n_\msf D)^2 \big) - \big( m_\msf D^2 + n_\msf D^2 - m_\msf C^2 - n_\msf C^2 \big)^2 & \\
=\> 4\, (m_\msf C m_\msf D + n_\msf D n_\msf C)^2 & ~. \end{aligned} \nonumber
\end{align*}
The proofs of the expressions for the remaining natural and inverse natural parameters are analogous.
\end{proof}
The above expressions for the natural parameters involve the Pl\"ucker coordinates in the exterior product $\mbf m \wedge \mbf n$ of $\mbf m, \mbf n \in \mbb R^4$; these will henceforth be denoted by
\begin{equation} \label{eq:plucker} \begin{aligned}
& p_\msf{AB} \:\coloneq\: m_\msf A n_\msf B - m_\msf B n_\msf A \,,\fit[1.5em] p_\msf{AC} \:\coloneq\: m_\msf A n_\msf C - m_\msf C n_\msf A \,,\fit[1.5em] p_\msf{AD} \:\coloneq\: m_\msf A n_\msf D - m_\msf D n_\msf A \,, \\
& p_\msf{BC} \:\coloneq\: m_\msf B n_\msf C - m_\msf C n_\msf B \,,\fit[1.5em] p_\msf{BD} \:\coloneq\: m_\msf B n_\msf D - m_\msf D n_\msf B \,,\fit[1.5em] p_\msf{CD} \:\coloneq\: m_\msf C n_\msf D - m_\msf D n_\msf C \,,
\end{aligned} \end{equation}
so that $\mbf m \wedge \mbf n = [\fit p_\msf{AB}\, ,\ldots\,,\, p_\msf{CD}\fit]^\top \in \mbb R^6$.
As is well known (see e.g.~Refs.~\cite{Cantarella:2019, Michalek:2021}), these satisfy the Pl\"ucker identity $p_\msf{AB}\fit p_\msf{CD}\fit - p_\msf{AC}\fit p_\msf{BD}\fit + p_\msf{AD}\fit p_\msf{BC} = 0$, which in turn defines the Klein quadric $\mathcal K \coloneq \{ \mbf p \in \mbb R^6 \mid p_1 p_6 - p_2 p_5 + p_3 p_4 \,=\, 0 \}$.
Together with Proposition \myref{thm:pc-to-np}, this notation allows us to express the products of pairs of \textquote{opposite} natural parameters~as
\begin{equation} \label{eq:products}
u\, z ~=~ 4\, p_\msf{AB}^{\,2}\, p_\msf{CD}^{\,2} / s^2 ~,\quad
v\, y ~=~ 4\, p_\msf{AC}^{\,2}\, p_\msf{BD}^{\,2} / s^2 ~,\quad
w\, x ~=~ 4\, p_\msf{AD}^{\,2}\, p_\msf{BC}^{\,2} / s^2 ~,
\end{equation}
where the squared exterior surface area is $s^2/4 \,=\, \| \mbf m \|^2 \| \mbf n \|^2 =\, \| \mbf m \wedge \mbf n \|^2$.

The key to fully defining the correspondence between the natural parameters of degenerate tetrahedra and the Klein quadric, viewed non-projectively as a five-dimensional variety $\mathcal K \subset \mbb R^6$, is to observe the similarity between the factors in Eq.~(\myref[2]{eq:ptolemy}), which factorizes the polynomial $\Omega$ of Eq.~(\myref[2]{eq:myform}) into four quadratic factors in the square-roots $\hat u, \ldots, \hat z$ of the natural parameters, and the Pl\"ucker identity.
Which of these quadratic factors corresponds to the Pl\"ucker identity depends on the relative signs of the products of \textquote{opposite} pairs of Pl\"ucker coordinates.
Specifically, by Eq.~(\myref{eq:products}) together with the Pl\"ucker identity itself,
\begin{align} \label{eq:klein}
(s / 2)\, & \big( \mrm{sign}(p_\msf{AB}\fit p_\msf{CD})\, \hat u \hat z \,-\, \mrm{sign}(p_\msf{AC}\fit p_\msf{BD})\, \hat v \hat y \,+\, \mrm{sign}(p_\msf{AD}\fit p_\msf{BC})\, \hat w \hat x \big) ~= \\ \nonumber 
& \mrm{sign}(p_\msf{AB} p_\msf{CD})\, \lvert p_\msf{AB} p_\msf{CD} \rvert \,-\, \mrm{sign}(p_\msf{AC} p_\msf{BD})\, \lvert p_\msf{AC} p_\msf{BD} \rvert \,+\, \mrm{sign}(p_\msf{AD} p_\msf{BC})\, \lvert p_\msf{AD} p_\msf{BC} \rvert \\ \nonumber 
& \hskip0.4\textwidth =~  p_\msf{AB}\, p_\msf{CD} \,-\, p_\msf{AC}\, p_\msf{BD} \,+\, p_\msf{AD}\, p_\msf{BC} ~=~ 0 ~. \fit[-0.5em]
\end{align}
Clearly there are eight possible combinations of signs for the three terms in this equation, each of which corresponds to one of the four factors in Eq.~(\myref[2]{eq:ptolemy}) or (equivalently) its negative vanishing.
The two corresponding to the first factor in Eq.~(\myref[2]{eq:ptolemy}), namely those with $\mrm{sign}(p_\msf{AB}\fit p_\msf{CD}) = -\mrm{sign}(p_\msf{AC}\fit p_\msf{BD}) = \mrm{sign}(p_\msf{AD}\fit p_\msf{BC})$, only hold when all three products of pairs of opposite Pl\"ucker coordinates vanish simultaneously.
In the following, these factors will be denoted by
\begin{equation} \label{eq:omegafact} \begin{aligned}
\Omega_0(\hat u, \hat v, \hat w, \hat x, \hat y, \hat z) ~\coloneq~ \hat u \hat z + \hat v \hat y + \hat w \hat x \,,\quad & \Omega_1(\hat u, \hat v, \hat w, \hat x, \hat y, \hat z) ~\coloneq~ \hat v \hat y + \hat w \hat x - \hat u \hat z \,, \\
\Omega_2(\hat u, \hat v, \hat w, \hat x, \hat y, \hat z) ~\coloneq~ \hat w \hat x + \hat u \hat z - \hat v \hat y \,,\quad & \Omega_3(\hat u, \hat v, \hat w, \hat x, \hat y, \hat z) ~\coloneq~ \hat u \hat z + \hat v \hat y - \hat w \hat x \,.
\end{aligned} \end{equation}

The main complication to be dealt with in fully defining the correspondence between the zeros of $\Omega$ and $\mathcal K$ stems from the fact that the signs of the square-roots of the four complex numbers given by the components of the vectors in Eq.~(\myref{eq:r2avec}) are arbitrary.
To simplify the presentation we will restrict ourselves for now to zeros of $\Omega$ where the natural parameters are all strictly positive and none of the Pl\"ucker coordinates vanish, so they correspond to vectors in the \textquote{generic} Klein quadratic $\mathcal K^* \coloneq \{\, \mbf p \in \mathcal K \mid p_1,\,\ldots,\, p_6 \ne 0 \,\}$.
Then, starting from any given choice of signs for these four square-roots, all the others are obtained by reflecting the vectors $\mbf m, \mbf n \in \mbb R^4$ in the subspaces orthogonal to the four coordinate axes.
The group generated by these reflections is isomorphic to the direct product $\mbb Z_2^4$ of four cyclic groups $\mbb Z_2$ of order $2$, and the action of this group on the Pl\"ucker coordinates $\mbf m \wedge \mbf n$ either maintains the signs of the products of opposite pairs $p_\msf{AB} p_\msf{CD}$, $p_\msf{AC} p_\msf{BD}$, $p_\msf{AD} p_\msf{BC}$ thereof or else changes all those signs identically.

This group action, however, is not faithful because changing the signs of all four of the square-roots does not change the signs of the Pl\"ucker coordinates, although changing all the latter's signs identically also maintains the relative signs of the products of opposite pairs.
A faithful action of $\mbb Z_2^4$ on the generic Klein quadric $\mathcal K^*$ which does give all the combinations of signs consistent with a specific factor $\Omega_k$ ($k = 1,2,3$) of $\Omega$ vanishing is defined by $[\epsilon_0 \epsilon_1 p_\msf{AB}, \epsilon_0 \epsilon_2 p_\msf{AC}, \epsilon_0 \epsilon_3 p_\msf{AD}, \epsilon_3 p_\msf{BC}, \epsilon_2 p_\msf{BD}, \epsilon_1 p_\msf{CD}]$, where $\epsilon_0, \epsilon_1, \epsilon_2, \epsilon_3$ can assume any combination of values in $\{ -1, +1 \}$.
In what follows, the set of sixteen-vector orbits generated by this action of $\mbb Z_2^4$ on $\mathcal K^*$ will be denoted by $\mbb Z_2^4 \circ \mathcal K^*$, and the corresponding quotient by $\mathcal K^* /\, \mbb Z_2^4$.

\begin{theorem} \label{thm:kleinquad}
The natural parameters $u, v, w, x, y, z > 0$ of a tetrahedron $\GEO{ABCD}$ with $\Omega(u,v,w,x,y,z) = 0$ are in one-to-one correspondence with points in the quotient $\mathcal K^* /\, \mbb Z_2^4$ via Eq.~(\myref{eq:klein}), or equivalently, with sixteen-vector subsets of $\mathcal K^*$ where these subsets are the orbits of the group action $\smash{\mbb Z_2^4} \circ \mathcal K^*$.
\end{theorem}
\begin{proof}
First, note that the given conditions together with the second part of Lemma \myref{thm:techlem} imply that the inverse parameters $\tilde u, \tilde v, \tilde w, \tilde x, \tilde y, \tilde z$ are non-negative, so the parameters $u, v, w, x, y, z$ determine a proper but degenerate Euclidean tetrahedron.
We may now define a corresponding vector of Pl\"ucker coordinates (cf.~Eq.~(\myref{eq:omegafact})) as follows:
\begin{equation*}
\fit[-1em] [p_\msf{AB}, p_\msf{AC}, p_\msf{AD}, p_\msf{BC}, p_\msf{BD}, p_\msf{CD}] ~\coloneq~ \begin{cases}
\sqrt{s/2}\; [ \hat z, ~\hat y, ~\hat x, \fit[0.75em]\hat w,\, -\hat v,\, -\hat u] & \text{if~}~\Omega_1 ~=~ 0; \\
\sqrt{s/2}\; [ \hat z, ~\hat y, ~\hat x, \fit[0.85em]\hat w, \fit[0.75em]\hat v, \fit[0.75em]\hat u] & \text{if~}~\Omega_2 ~=~ 0; \\
\sqrt{s/2}\; [ \hat z, ~\hat y, ~\hat x,\fit[0.2em] -\hat w,\fit[0.2em] -\hat v, \fit[0.75em]\hat u] & \text{if~}~\Omega_3 ~=~ 0; \end{cases}
\end{equation*}
(note that $\Omega_0 > 0$ since $\hat u, \ldots, \hat z$ are strictly positive).
These coordinates will satisfy the Pl\"ucker identity, and letting $\mbb Z_2^4$ act on them as previously described will generate a set of sixteen distinct vectors on $\mathcal K^*$ corresponding to the same factor $\Omega_k = 0$ ($k = 1,2,3$).

Conversely, given any set of sixteen vectors in an orbit $\mbb Z_2^4 \circ \mathcal K^*$, the relative signs of the products of the three opposite pairs of Pl\"ucker coordinates therein will consistently determine which of the factors $\Omega_k$ ($k = 1,2,3$) of $\Omega$ vanishes, where the square-roots of the natural parameters in question are $\sqrt{\fit2/s}$ times the coordinates' absolute values as in Proposition \myref{thm:pc-to-np}.
It follows that these natural parameters satisfy $\Omega = 0$ as desired; the non-negativity of the inverse natural parameters then follows from the second part of Lemma \myref{thm:techlem} as before, so these natural parameters indeed correspond to a proper degenerate Euclidean tetrahedron.
\end{proof}
\begin{remark}
Clearly this correspondence between $\mathcal K^* / \mbb Z_2^4$ and the generic zeros of $\Omega$ is bicontinuous, so by continuity the zeros of $\Omega$ are canonically homeomorphic to $\mathcal K / \mbb Z_2^4$ in the usual (quotient) topology.
Note also that these sixteen-vector subsets of $\mathcal K^*$ will consist of eight \textquote{antipodal} pairs related by an overall change of sign.
Given any vector in $\mathcal K$ for which one or more of the Pl\"ucker coordinates are zero, these sign changes will of course have no effect on those coordinates and hence some of the sixteen vectors in the corresponding orbit will coincide.
While exactly two non-opposite natural parameters cannot vanish unless $\Omega(u,\,\ldots,\, z) < 0$, it is possible for three non-opposite parameters to vanish, in which case all four factors in Eq.~(\myref{eq:omegafact}) will vanish and the signs of the three non-zero Pl\"ucker coordinates can be chosen arbitrarily, giving rise to only four antipodal pairs of vectors in $\mathcal K$.
For example, if $[u,v,w,x,y,z] = [1,1,1,0,0,0]$, then $p_\msf{AB} = p_\msf{AC} = p_\msf{AD} = 0$ while $p_\msf{BC}, p_\msf{BD},p_\msf{CD} = \pm\smash{\sqrt3}$ and $\mrm{rank}\big( \mbf G_\msf A \big) = 2$. 
The example $[u,v,w,x,y,z] = [0,0,0,1,1,1]$, where $p_\msf{AB}, p_\msf{AC}, p_\msf{AD} = \pm\smash{\sqrt3}$ while $p_\msf{BC} = p_\msf{BD} = p_\msf{CD} = 0$, shows that this can also happen when $\mrm{rank}(\mbf G_\msf A) = 1$.
\end{remark}

\section{A stratification giving the zeros a lattice structure} \label{sec:lattice}
Recall that when $\Omega = 0$ the products of complementary natural and inverse natural parameters are the squared distances among four points on a line.
This may also be shown by extending the two-dimensional basis of Eq.~(\myref{eq:r2avec}) to an orthonormal three-dimensional basis $[\mbf e_1, \mbf e_2, \mbf e_3]$ and computing the cross products of the coplanar areal vectors therein.
These will all be parallel to $\mbf e_3$, and their components along that axis are simply
\begin{align} 
\begin{aligned} 
\mbf e_3 \Dprod \big( \big( \VCP{AB}{AC} \big) \Xprod \big(\! \Vdiff \VCP{AB}{AD} \big) \big) \,= & \\[-2pt] (m_\msf D^2 - n_\msf D^2)\, 2\, m_\msf C n_\msf C - {}&{} 2\, m_\msf D n_\msf D\, (m_\msf C^2 - n_\msf C^2) \:=\: \pm2\, p_\msf{CD}\, q_\msf{CD} \,, \end{aligned} & \nonumber \\[2pt]
\begin{aligned} 
\mbf e_3 \Dprod \big( \big( \VCP{AB}{AC} \big) \Xprod \big( \VCP{AC}{AD} \big) \big) \,= & \\[-2pt] (m_\msf D^2 - n_\msf D^2)\, 2\, m_\msf B n_\msf B - {}&{} 2\, m_\msf D n_\msf D\, (m_\msf B^2 - n_\msf B^2) \:=\: \pm2\, p_\msf{BD}\, q_\msf{BD} \,, \end{aligned} & \nonumber \\[2pt]
\begin{aligned} 
\mbf e_3 \Dprod \big( \big( \VCP{AB}{AC} \big) \Xprod \big(\! \Vdiff \VCP{BC}{BD} \big) \big) \,= & \\[-2pt] (m_\msf D^2 - n_\msf D^2)\, 2\, m_\msf A n_\msf A - {}&{} 2\, m_\msf D n_\msf D\, (m_\msf A^2 - n_\msf A^2) \:=\: \pm2\, p_\msf{AD}\, q_\msf{AD} \,, \end{aligned} & \nonumber \\[2pt]
\begin{aligned} 
\mbf e_3 \Dprod \big( \big(\! \Vdiff \VCP{AB}{AD} \big) \Xprod \big( \VCP{AC}{AD} \big) \big) \,= & \\[-2pt] (m_\msf C^2 - n_\msf C^2)\, 2\, m_\msf B n_\msf B - {}&{} 2\, m_\msf C\, n_\msf C\, (m_\msf B^2 - n_\msf B^2) \:=\: \pm2\, p_\msf{BC}\, q_\msf{BC} \,, \end{aligned} & \label{eq:anotherway} \\[2pt]
\begin{aligned} 
\mbf e_3 \Dprod \big( \big( \Vdiff \VCP{AB}{AD} \big) \Xprod \big( \Vdiff \VCP{BC}{BD} \big) \big) \,= & \\[-2pt] (m_\msf C^2 - n_\msf C^2)\, 2\, m_\msf A n_\msf A - {}&{} 2\, m_\msf C n_\msf C\, (m_\msf A^2 - n_\msf A^2) \:=\: \pm2\, p_\msf{AC}\, q_\msf{AC} \,, \end{aligned} & \nonumber \\[2pt]
\begin{aligned} 
\mbf e_3 \Dprod \big( \big( \VCP{AC}{AD} \big) \Xprod \big(\! \Vdiff \VCP{BC}{BD} \big) \big) \,= & \\[-2pt] (m_\msf B^2 - n_\msf B^2)\, 2\, m_\msf A n_\msf A - {}&{} 2\, m_\msf B n_\msf B\, (m_\msf A^2 - n_\msf A^2) \:=\: \pm2\, p_\msf{AB}\, q_\msf{AB} \,, \end{aligned} \nonumber &
\end{align}
where $p_\msf{ab}$ are the Pl\"ucker coordinates from Eq.~(\myref{eq:plucker}), $q_\msf{ab} \coloneq m_\msf a m_\msf b + n_\msf a n_\msf b$ are the corresponding inner products ($\msf a \ne \msf b \in \{ \msf A, \msf B, \msf C, \msf D \}$), and the signs depend on the orientation of $[\mbf e_1, \mbf e_2, \mbf e_3]$.
Observe that the sum of the cross products in the first three of these equations is minus the cross product of $\VCP{AB}{AC}$ with itself by Minkowski's identity (\myref[1]{eq:minkow}) and so equals $\mbf 0$; the sums of the other three triples of cross products sharing a common factor also vanish.

This shows that the quantities $\pm 2 p_\msf{ab} q_\msf{ab}$ in Eq.~(\myref{eq:anotherway}) are the signed distances among a quadruple of points along the $\mbf e_3$ axis.
Because their positions relative to the origin $\mbf 0$ are also determined by these relations but the orientation of the basis is arbitrary, they are better seen as a set of four collinear vectors modulo inversion in the origin.
Our next result gives the values of these vectors explicitly in terms of the natural parameters.
\begin{proposition} \label{thm:collinear}
Given a tetrahedron $\OL{\msf{ABCD}}$ for which $\Omega(u,v,w,x,y,z) = 0$, the complementary products $u\tilde u$, $\ldots\,$, $z\tilde z$ are the squared distances among at least one of the four quadruples of points given relative to the origin by the following four quadruples of vectors in $\mbb R^1 \approx \mbb R$ (wherein $\hat s \coloneq \sqrt{s/2}$):
\begin{equation} \label{eq:collinvecs} \begin{aligned}
\{\, \hat u \hat v \hat w / \hat s ,~\;\quad \hat u \hat x \hat y / \hat s ,\;\quad \hat v \hat x \hat z / \hat s ,\;\quad \hat w \hat y \hat z / \hat s \,\} \quad&\text{if }~ \Omega_0(\hat u, \hat v, \hat w, \hat x, \hat y, \hat z) ~=~ 0 \,; \\
\{\, \hat u \hat v \hat w / \hat s ,\fit[1.5em] \hat u \hat x \hat y / \hat s ,\;\, -\,\hat v \hat x \hat z / \hat s ,~ -\hat w \hat y \hat z / \hat s \,\} \quad&\text{if }~ \Omega_1(\hat u, \hat v, \hat w, \hat x, \hat y, \hat z) ~=~ 0 \,; \\
\{\, \hat u \hat v \hat w / \hat s ,\;~ -\hat u \hat x \hat y / \hat s ,\quad \hat v \hat x \hat z / \hat s ,~ -\hat w \hat y \hat z / \hat s \,\} \quad&\text{if }~ \Omega_2(\hat u, \hat v, \hat w, \hat x, \hat y, \hat z) ~=~ 0 \,; \\
\{\, \hat u \hat v \hat w / \hat s ,\,~ -\hat u \hat x \hat y / \hat s ,~ \!-\hat v \hat x \hat z / \hat s ,\fit[1.3em] \hat w \hat y \hat z / \hat s \,\} \quad&\text{if }~ \Omega_3(\hat u, \hat v, \hat w, \hat x, \hat y, \hat z) ~=~ 0 \,.
\end{aligned} \end{equation}
\end{proposition}
\begin{proof}
Let us define:
\begin{equation*} \begin{aligned}
\tilde u_\pm ~\coloneq~ \frac{2\, (\hat v \hat w \,\pm\, \hat x \hat y)^2}{s} ~;\quad & \tilde v_\pm ~\coloneq~ \frac{2\, (\hat u \hat w \,\pm\, \hat x \hat z)^2}{s} ~;\quad \tilde w_\pm ~\coloneq~ \frac{2\, (\hat u \hat v \,\pm\, \hat y \hat z)^2}{s} ~; \\
\tilde x_\pm ~\coloneq~ \frac{2\, (\hat u \hat y \,\pm\, \hat v \hat z)^2}{s} ~;\quad & \tilde y_\pm ~\coloneq~ \frac{2\, (\hat u \hat x \,\pm\, \hat w \hat z)^2}{s} ~;\quad\> \tilde z_\pm ~\coloneq~ \frac{2\, (\hat v \hat x \,\pm\, \hat w \hat y)^2}{s} ~.
\end{aligned} \end{equation*}
Then if we regard $\tilde u,\, \ldots,\, \tilde z$ as rational functions of the natural parameters $u,\, \ldots,\, z$ as in Eq.~(\myref[2]{eq:natinvnat}), it is easily shown that when e.g.~$\Omega_1(\hat u, \hat v, \hat w, \hat x, \hat y, \hat z) = 0$:
\begin{equation*} \begin{aligned}
\tilde u\, \rvert_{uz = (\hat v \hat y + \hat w \hat x)^2} ~=~ \tilde u_- ~;\qquad & \tilde z\, \rvert_{uz = (\hat v \hat y + \hat w \hat x)^2} ~=~ \tilde z_- ~; \\
\tilde v\, \rvert_{vy = (\hat u \hat z - \hat w \hat x)^2} ~=~ \tilde v_+ ~;\qquad & \tilde y\, \rvert_{vy = (\hat u \hat z - \hat w \hat x)^2} ~=~ \tilde y_+ ~; \\
\tilde w\, \rvert_{wx = (\hat u \hat z - \hat v \hat y)^2} ~=~ \tilde w_+ \,;\qquad & \tilde x\, \rvert_{wx = (\hat u \hat z - \hat v \hat y)^2} ~=~ \tilde x_+ ~.
\end{aligned} \end{equation*}
It is now a matter of inspection to verify that, on multiplying these expressions for the inverse parameters by the corresponding natural parameters, the resulting complementary products are indeed the squared distances among the four points given by the vectors on the second line of Eq.~(\myref{eq:collinvecs}).
The proofs for the other three cases are similar.
\end{proof}
\begin{remark} \label{rem:cubic}
Given any $\alpha, \beta, \gamma, \delta, \varsigma \in \mbb R$ with $\varsigma > 0$ and $\alpha\beta\gamma\delta > 0$ but not all with the same sign, the equations $\alpha^2 = 2\, uvw / \varsigma$, $\beta^2 = 2\, uxy / \varsigma$, $\gamma^2 = 2\, vxz / \varsigma$, $\delta^2 = 2\, wyz / \varsigma$ together with $\varsigma = s = 2\, (u\!+\!v\!+\!w\!+\!x\!+\!y\!+\!z)$ and $\Omega = 0$ can in principle be solved to obtain a relation between these five numbers and the natural parameters of degenerate tetrahedra.
The key is to note that when the first four of these equations hold, we have
\begin{equation} \label{eq:onprats}
(\beta \gamma) / (\alpha \delta) ~=~ x / w ~,\qquad (\beta \delta) / (\alpha \gamma) ~=~ y / v ~,\qquad (\gamma \delta) / (\alpha \beta) ~=~ z / u ~,
\end{equation}
and that those equations are valid regardless of the signs of $\alpha, \beta, \gamma, \delta$ as above.
Thus they may be used to eliminate three non-opposite natural parameters from the equations $\Omega = 0$ and $\varsigma = 2(u+v+w+x+y+z)$.
This yields two linear equations in three unknowns, say $u, v, w$, with which two of these three unknowns may be eliminated from $uvw = \varsigma\alpha^2 / 2$, say $v$ \& $w$, resulting in a cubic equation for the remaining parameter $u$.
Numerical examples suggest that when $\varsigma$ exceeds a certain positive lower bound determined by $\alpha, \beta, \gamma$ \& $\delta$ this cubic generically has three distinct real roots, two of which comprise positive values for $u$ and yield positive values for $v$ \& $w$ upon back-substitution; moreover the same parameters are obtained if $\alpha, \beta, \gamma$ \& $\delta$ are all negated, thereby establishing a $2\!:\!2$ relation between these two semi-algebraic sets.
Further discussion of this less redundant, but more complicated, parametrization of the set of generic degenerate tetrahedra may be found in Appendix \myref[4]{sec:2to2}.
\end{remark}
\noindent Note that if $\Omega_0 = 0$ some triple of non-opposite parameters must vanish (e.g.~$u, v, w = 0$; $u, v, x = 0$; etc.), in which case at most one of the vectors $\alpha, \beta, \gamma, \delta$ is non-zero, and that they all vanish whenever any pair of opposite parameters are both zero (e.g.~$u, z = 0$). \smallskip

We now turn to the non-generic case in which two or more of the four collinear vectors coincide, or equivalently, some of the complementary products $u\tilde u$, $v\tilde v$ etc.~vanish.
The fact that \textquote{distance equals zero} is an equivalence relation implies that only certain combinations of the complementary products can vanish simultaneously, which in turn implies the lattice of inclusions depicted in Fig.~\myref{fig:lattice} among the semi-algebraic sets defined by $\Omega = 0$ and such simultaneous subsets of the six equations $u\tilde u = 0\fit,\, \ldots\fit,\, z\tilde z = 0$.
Via Theorem \myref{thm:kleinquad}, this combinatorial structure leads to an apparently novel stratification of the Grassmannian corresponding to $\mathcal K$, although that will not be further explored here.

\begin{figure}
\input{Fig_III-1.tex}
\label{fig:lattice}
\end{figure}

While the semi-algebraic set in $\mbb R_{\ge0}^{\,6}$ defined by $\Omega = 0$ is five-dimensional, also requiring any one complementary product to vanish reduces the dimension of that solution set to four.
The four exterior areas can then be used as local coordinates on each of these six four-dimensional sets since, as will now be illustrated, the three interior areas can be computed from them.
Clearly one of the interior areas is given by the sum or a difference of two of the exterior areas, depending on which tetrahedron inequality saturates (i.e.~$\Tau_k[\msf a, \msf b] = 0$) in order to make the corresponding complementary product vanish.
To see how to obtain the other two, suppose for example that $u = 0 ~\Rightarrow~ 2 \MAG{AB \tmv CD} = \MAG{ABC} + \MAG{ABD}$.
Upon using Yetter's identity $\Xi = 0$ to eliminate $\MAG{AD \tmv BC}$ from the Gramian $\mrm{det}\big( \mbf G_\msf A \big)$ and then eliminating $\MAG{AB \tmv CD}$ from the result using this linear relation, we obtain
\begin{equation}
0 ~\le~ \mrm{det}\big( \mbf G_\msf A \big) ~=~ \begin{aligned}[t] -16\, \Big( & \MAG{ABC} \big( \MAG{ABD}^2 - \MAG{BCD}^2 \big) + \MAG{ABD} \big( \MAG{ABC}^2 - \MAG{ACD}^2 \big) \\[-3pt]
& +\, 4\, \MAG{AC|BD}^2 \big( \MAG{ABC} + \MAG{ABD} \big) \Big)^{\!2} ~. \\[-3pt]
\end{aligned} \end{equation}
Thus the polynomial inside the main parentheses vanishes and can readily be solved to obtain $\smash{\MAG{AC \tmv BD}}^2$, whereupon $\smash{\MAG{AD \tmv BC}}^2$ may be obtained from Yetter's identity.
In a similar fashion one can show that the exterior areas uniquely determine the interior whenever any single one of the tetrahedron inequalities saturates.

Requiring more than one complementary product to vanish immediately reduces the generic dimensionality of the subset of degenerate tetrahedra defined by that requirement to three.
In contradistinction with intersection theory over algebraically closed fields, the dimension remains three even as additional complementary products are required to vanish; this, of course, is because these semi-algebraic sets intersect non-transversely in $\mbb R^6$.
In all of these cases Yetter's identity, a vanishing Gramian and the linear equations corresponding to those tetrahedron inequalities that saturate can be solved to express the exterior areas in terms of the \underline{interior}.
Although a full proof would involve a fairly intricate case-by-case analysis, it appears that whenever the equations are consistent their solutions (the exterior areas) are simply signed sums of the three interior areas, i.e.~linear combinations with coefficients equal to $\pm1$.
These coefficients are determined by which combinations of the eighteen tetrahedron inequalities are saturated.

At the lowest level of the hierarchy of three-dimensional solutions all six complementary products vanish, meaning the $1\mrm D$ vectors of Proposition \myref{thm:collinear} all coincide and at least one $\Tau_k[\msf a, \msf b] = 0$ in each of the six triples of tetrahedron inequalities.
These six linear equations in the four unknown exterior areas generally do not have a simultaneous solution, but when they do it is a signed sum of the interior areas as above which also satisfies Yetter's identity as well as $uvw = uxy = vxz = wyz = 0 = \Omega$.
We close this section by proving, as indicated in Fig.~\myref{fig:lattice}, that this level of the lattice corresponds to those configurations for which $\mrm{rank}\big( \mbf G_\msf A \big) = 1$.
\begin{proposition} \label{thm:planar}
Given a (not-necessarily-degenerate) Euclidean tetrahedron $\GEO{ABCD}$ with natural parameters $u, v, w, x, y, z$ and inverse natural parameters $\tilde u, \tilde v, \tilde w, \tilde x, \tilde y, \tilde z$, the rank of the Gram matrix $\mbf G_\msf A$ at vertex $\GEO{A}$ (or any other vertex) is at most $1$ if \& only if
\begin{equation} \label{eq:complementarity}
u\tilde u ~=~ v\tilde v ~=~ w\tilde w ~=~ x\tilde x ~=~ y\tilde y ~=~ z\tilde z ~=~ 0 ~.
\end{equation}
Thus the planar situation is characterized by these complementarity relations between the natural and inverse natural parameters.
\end{proposition}
\begin{proof}
In proving this proposition, it is convenient to consider the full $4\times4$ Gram matrix of twice the outwards-pointing areal vectors of the exterior faces $\mbf G_\msf{ext\,}$.
Letting $F_\msf{abc;abd} \coloneq (F_\msf{ab|cd} - F_\msf{abc} - F_\msf{abd}) / 2$ for all $\{ \msf a, \msf b, \msf c, \msf d \} = \{ \msf A, \msf B, \msf C, \msf D \}$, where $F_\msf{ABC} \leftrightarrow 4\fit\smash{\MAG{ABC}}^2$, $\ldots\,$,  $F_\msf{AD|BC} \leftrightarrow 16\fit\smash{\MAG{AD \tmv BC}}^2$ are indeterminates representing the squared areas as usual, the matrix of polynomials corresponding to $\mbf G_\msf{ext}$ is
\begin{equation*}
\mbf G_F[\msf A, \msf B, \msf C, \msf D] ~\coloneq~
\begin{bmatrix}
F_\msf{ABC} & F_\msf{ABC;ABD} & F_\msf{ABC;ACD} & F_\msf{ABC;BCD} \\
F_\msf{ABD;ABC} & F_\msf{ABD} & F_\msf{ABD;ACD} & F_\msf{ABD;BCD} \\
F_\msf{ACD;ABC} & F_\msf{ACD;ABD} & F_\msf{ACD} & F_\msf{ACD;BCD} \\
F_\msf{BCD;ABC} & F_\msf{BCD;ABD} & F_\msf{BCD;ACD} & F_\msf{BCD}
\end{bmatrix} ~.
\end{equation*}
Here the upper-left $3\times3$ block is the Gram matrix $\mbf G_F[\msf A]$ at $\msf A$ from Eq.~(\myref[1]{eq:gramatA}), while the last row/column are $F_\msf{ABC;BCD} = F_\msf{BCD;ABC} = (F_\msf{AD|BC} - F_\msf{ABC} - F_\msf{BCD}) / 2$, $F_\msf{ABD;BCD} = F_\msf{BCD;ABD} = (F_\msf{AC|BD} - F_\msf{ABD} - F_\msf{BCD}) / 2$, $F_\msf{ACD;BCD} = F_\msf{BCD;ACD} = (F_\msf{AB|CD} - F_\msf{ACD} - F_\msf{BCD}) / 2$.

As a Gram matrix amongst vectors in a Euclidean space $\mbf G_\msf{ext}$ is assured of being positive semi-definite, as are its $3\times3$ principal submatrices $\mbf G_\msf A$, $\mbf G_\msf B$, $\mbf G_\msf C$ \& $\mbf G_\msf D$.
It is also easily seen that $\mbf 1 = [1,1,1,1]^\top$ is an eigenvector of $\mbf G_F[\msf A, \msf B, \msf C, \msf D]$ with eigenvalue \smash{$-\fit\breve\Xi_F / 2 = 0$}, so $\mrm{det}\big( \mbf G_\msf{ext} \big) = 0$.
If $\mrm{rank}\big( \mbf G_\msf A \big) = 1$, then the areal vectors of the three exterior faces meeting at \smash{$\GEO{A}$} are collinear, and Minkowski's identity (\myref[1]{eq:minkow}) requires that the areal vector of the fourth exterior face $\GEO{BCD}$ also be collinear with those vectors.
It follows that $\mrm{rank}\big( \mbf G_\msf{ext} \big) = 1 = \mrm{rank}\big( \mbf G_\msf B \big) = \mrm{rank}\big( \mbf G_\msf C \big) = \mrm{rank}\big( \mbf G_\msf D \big) $, as well.
Because the $2\times2$ principal minors of $\mbf G_\msf{ext}$ are $\Tau_{\!0}[\msf A, \msf B]\fit\Tau_{\!1}[\msf A, \msf B]\fit\Tau_{\!2}[\msf A, \msf B]\fit\Tau_{\!3}[\msf A, \msf B] / 4$ $=$ $s^2 u\tilde u$ etc., this establishes that Eq.~(\myref{eq:complementarity}) holds if $s > 0$.
It also holds, of course, if $s = 2(u+\cdots+z) = 0$ and hence in general, as claimed.

\begin{figure}[!b]
\input{Fig_III-2.tex}
\label{fig:regions}
\end{figure}

\smallskip Conversely, if Eq.~(\myref{eq:complementarity}) holds then all the $2\times2$ principal minors of $\mbf G_\msf{ext}$ vanish, and since its determinant also vanishes this matrix will have a rank of $1$ if all of its $3\times3$ principal minors vanish, i.e.~$\mrm{det}\big( \mbf G_\msf A \big) = \mrm{det}\big( \mbf G_\msf B \big) = \mrm{det}\big( \mbf G_\msf C \big) = \mrm{det}\big( \mbf G_\msf D \big) = 0$.
These determinants are equal by Lemma \myref[1]{thm:gramequiv} since $\Xi = 0$ by Proposition \myref[2]{thm:yidok}, while $\mrm{det}\big( \mbf G_\msf A \big) = t^4$ by Eq.~(\myref[1]{eq:gramatA}) and $t^4 = s^2\, \Omega(u,v,w,x,y,z)$ by Theorem \myref[2]{thm:myform}.
But Lemma \myref{thm:techlem} shows that if any one of the complementary products $u\tilde u\,, \ldots\,, z\tilde z$ vanish then $\Omega(u,v,w,x,y,z) \le 0$, whereas $\Omega(u,v,w,x,y,z) \ge 0$ in the Euclidean case assumed here.
It follows that both $\Omega$ and all the $3\times3$ principal minors vanish as desired.
\end{proof}

\begin{remark} \label{rem:rank1}
As illustrated in Fig.~\myref{fig:regions}, only seven of the $2^6 = 64$ possible combinations of complementary natural or inverse natural (but not both) parameters vanishing so as to satisfy Eq.~(\myref{eq:complementarity}) can occur in any planar configuration; these are exactly those for which the vectors $\alpha, \beta, \gamma, \delta \in \mbb R^1$ in Remark \myref{rem:cubic} are all zero.
In analogy to the way in which the $3! / 2 = 3$ linear orders (up to inversion) for three distinct points on a line correspond to the zeros of Heron's formula, these combinations correspond to the seven uniform rank $3$ \TDEF{chirotopes\/} (aka affine oriented matroids) of four-point configurations in the affine plane \cite{Bjorner:1999,Richter-Gebert:2017}.
The three combinations wherein every vertex is an extreme point of their convex hull may be further divided into four subclasses each, depending on exactly which combinations of tetrahedron inequalities saturate, and these in turn may be distinguished by their allowable sequences as defined by Goodman \& Pollack \cite{Goodman:1983, Goodman:1990}.
Although distances realizing the given areas (or equivalently natural parameters) do not exist when $\mrm{rank}\big( \mbf G_\msf{ext} \big) = 2$, in the planar ($\mrm{rank}\big( \mbf G_\msf{ext} \big) = 1$) situation they exist and determine the associated chirotope \cite[\S3.6]{Crippen:1988}.
They are, however, not unique because the distances in any special (area preserving) affine transform thereof will also realize those same areas via Eqs.~(\myref[1]{eq:cmd3}) \& (\myref[1]{eq:talata}).
\end{remark}
The details of this analysis of the combinatorics of the rank $1$ zeros, the existence of canonical distances realizing the areas in the rank $1$ case, an extension of Theorem \myref[1]{thm:bdsmc} that characterizes the algebraic relations among the polynomials involved even if the indeterminates therein are not equal to Euclidean invariants in a non-degenerate tetrahedron, and the equations involved in the $2\!:\!2$ parametrization of the generic rank $2$ zeros from Remark \myref{rem:cubic}, may be found in Appendices \myref[4]{sec:planar} -- \myref[4]{sec:2to2} of the last part of this series, respectively; Appendix \myref[4]{sec:invol} gives the definitions and a preliminary analysis of two non-commuting involutions on the set of all Euclidean tetrahedra, non-degenerate and otherwise.\pagebreak[2]

\section{Closing remarks} \label{sec:closing}
Collectively, Parts I--III of this series contained two rather substantial surprises.
The first is that anything therein was new, especially considering the classical nature of the subject matter and the elementary techniques used to derive the results.
This can be explained, at least in part, by the advent of computer algebra systems (such as the \texttt{SageMath} software package used for the calculations presented herein) which now enable an average student to accomplish many feats beyond even the grand masters of that distant era when low-dimensional Euclidean geometry was still at the cutting edge of mathematics and physics research \cite{Davis:1995}.
Vector algebra techniques, which the geometers of the late $19^\text{th}$ and early $20^\text{th}$ centuries seem to have largely relegated to the domain of physics, also proved enabling here.
The most salient reason, however, is probably the fact that the Hamiltonians of classical physics depend on the distances between pairs of particles but not on the areas spanned by triples thereof (let alone the areas of medial parallelograms).
As a result, many of the algebraic relations dealt with herein are not inherent in most people's physical intuition.
The possibility that these relations have a role to play in quantum mechanics \cite{Bengtsson:2017,Ferro:2021,Snider:2018} or quantum gravity \cite{Loll:2020,Rovelli:2015} remains to be explored.

The second, not entirely unrelated, surprise is of course the \textquote{projective} nature of the zeros of our extension of Heron's formula (\myref[2]{eq:myform}).
Given that this extension has a pretty strong claim to being intrinsic to classical Euclidean geometry, it is difficult to argue that these collinear tetrahedra with vertices separated by infinite distances are not part of three-dimensional Euclidean space.
Are they then also part of the Euclidean plane, and if not, where do they belong?
(One might almost be tempted to ask \textquote{is Flatland a bigger place than A.~Square ever knew?} \cite{Stewart:2002}!)
As noted in Section \myref{sec:zeros} above, they cannot reasonably be embedded even in the projective completion of Euclidean space, and furthermore, having certain well-defined metrical properties, namely the areas of the seven faces, or equiv\-alently, the six natural parameters, it seems unlikely they can be fully rationalized in purely projective terms.
If, as suggested in Remark \myref[2]{rem:hyperbolic}, hyperbolic geometry has a role to play in the geometric interpretation of the zeros of the determinant $-\Omega$, then the well-known fact that the inertial world-lines through an event in space-time can be viewed as points in the hyperboloid model of hyperbolic space \cite{Dray:2021,Ungar:2022} also hints at future applications in relativistic physics \cite{Penrose:1984,Regge:2000}.

Such interpretational issues aside, the results presented herein suggest a number of new lines of inquiry, the first of which was to work out a proof of Conjecture \myref[2]{thm:alldim} regarding the extension of the tetrahedron's formula (\myref[2]{eq:myform}) to higher dimensions.
Vector algebra only works in three dimensions, but higher dimensional generalizations are available, some of which actually predate it \cite{Crowe:1967}.
Today these are usually called Clifford algebras by mathematicians, although most of their users in the engineering and physics communities prefer the appellation \TDEF{geometric algebra\/}, as did W.~K.~Clifford himself (even though these algebras were considered only briefly in E.~Artin's much more recent book by that name \cite{Artin:1988}).
They acquire particular power when applied to a vector space model of inversive geometry, wherein the group of Euclidean similarities corresponds to the stabilizer of the point-at-infinity \cite{Doran:2003, Dorst:2009, Hestenes:1984, Li:2008, Li:2001, Sobczyk:2013}.
In that form, they proved ideally suited to the task of expressing the higher-dimensional analogues of the natural parameters as rational functions of the \textquote{hyper-areas} of the $n$-simplices' facets and $(2,n-1)$-medial sections (cf.~Proposition \myref[2]{thm:natpar}), as well as proving the conjecture itself \cite{Havel:2024}.

Beyond that, it might also be interesting to study the level sets defined by the equations $T = s^2\fit\Omega (u,\, \ldots,\, z)$ for $T > 0$, which can be shown to be likewise unbounded \cite{Havel:2024}, and in particular how they are stratified by the ratios of pairs of inter-vertex distances.
A purely algebraic challenge would be to invert the system of equations (\myref[2]{eq:natinvnat}) so as to obtain the natural parameters as (roots of?) rational functions of the inverse parameters; the corresponding problem in the plane (cf.~Remark \myref[2]{rem:extouch}) has the solution $u^2 = s\, \tilde v \tilde w / \tilde u$ etc., but it does not seem that $s$ can be simply expressed in terms of $\tilde u$, $\tilde v$ \& $\tilde w$ alone even in triangles.
Some of the ideas herein might fruitfully be applied to the study of three-dimensional polyhedra based in their triangulations \cite{Lee:2017,Sabitov:1998} (in this regard, it should be noted that a triangulation of any medial parallelogram of a tetrahedron is not generated by its standard barycentric subdivision).
It may also be possible to derive \textquote{Heron-like} formulae for polyhedra and other convex solids which admit an in-sphere, similar to those for tangential (circumscriptable) polygons \cite[Chs.~4 \& 13]{Apostol:2012}. 
The connections between the rank $1$ zeros and the order-theoretic structure of planar four-point configurations noted above in Remark \myref{rem:rank1} (cf.~Fig.~\myref{fig:regions}) may also inspire new developments in discrete and combinatorial geometry \cite{Richter-Gebert:2017, Seidel:1995, Ziegler:2000}.
Certainly, it will not be long before the computational commutative algebra community finds new directions in which to extend and generalize the elementary results presented herein \cite{Michalek:2021, Miller:2005}.

{\small

} 


\bigskip
\begin{center}
\medskip\rule{70pt}{2pt}\bigskip
\end{center}

\medskip
\begin{displayquote}\small
Wir mussen wissen. Wir werden wissen.

\smallskip\raggedleft\textit{from the tomb of David Hilbert}
\end{displayquote}

\medskip
\begin{displayquote}\small
After reading the section, \textquote{Problem Solvers and Theorizers,} a mathematician friend (one of the most distinguished living mathematicians) wrote that he would not speak to the author ever again
\textellipsis The truth offends.

\smallskip\raggedleft\textit{Gian-Carlo Rota, \textquote{Indiscrete Thoughts,} 1997}
\end{displayquote}

\appendix 
\mypart{\mytitle} 

\begin{abstract}
This is the last part of a series of four papers in The $\Pi\mrm{ME}$ Journal.
It consists of five appendices \myref{sec:planar}--\myref{sec:invol}, each of which considers and partially resolves certain issues left open in the first three parts, closing in some cases with further conjectures.
Appendix \myref{sec:planar} shows that the set of all quadruples of points in general position in the affine plane decomposes naturally into $16$ classes according to which combinations of the $18$ tetrahedron inequalities from Part I saturate, and that these classes refine the seven cases distinguished by their \TDEF{chirotopes} (aka affine oriented matroids).
Appendix \myref{sec:canonical} further shows that every such quadruple is canonically associated with a unique four-point configuration in the \EMPH{Euclidean} plane which minimizes the radius of gyration about their centroid subject to preserving the seven areas, and conjectures that this correspondence can be extended to the full set of zeros of the polynomial $\Omega$ in our extension of Heron's formula (Theorem \myref{thm:myform} of Part II).
Appendix \myref{sec:pmap} is devoted to a single theorem proving that in non-degenerate tetrahedra the image of the polynomial map from squared distances to squared areas is defined by Yetter's identity (cf.~Proposition \myref{thm:yid} in Part I), and extends this characterization to indeterminates representing squared areas and distances which may not be realizable in Euclidean space (and may even be negative).
Appendix \myref{sec:2to2} studies the parametrization of the generic zeros of $\Omega$ that was touched upon in Remark \myref{rem:cubic} of Part III, presenting evidence that it induces a two-to-two correspondence between these five parameters and the zeros, and shows that in non-degenerate tetrahedra the absolute values of the four ``$1$-dimensional vector'' parameters are nothing but the vertex-to-in-touch-point distances times the in-radius.
Finally, Appendix \myref{sec:invol} delves into two involutions on the set of all tetrahedra, one of which has  previously been studied by M.~Fiedler and the other of which appears to be new, and posits that they may jointly generate a finite group.
\end{abstract}


\section{The combinatorial structure of the rank 1 zeros} \label{sec:planar}%
Recall that $\mbf G_\msf A$ is the Gram matrix at vertex $\GEO A$ of a tetrahedron $\GEO{ABCD}$, as defined in Eq.~(\myref[1]{eq:gramatA}),\footnote{
A Roman numeral followed by a colon specifies a reference to a numbered entity in another paper of this series, in this case Eq.~(\myref{eq:gramatA}) from Part I.} 
and that Proposition \myref[3]{thm:planar} showed the $\mrm{rank}\big( \mbf G_\msf A \big) \le 1$ situation is characterized by a complementarity relation between the natural and inverse natural parameters, as defined in Part II.
It was then noted in Remark \myref[3]{rem:rank1} that only seven of the $2^6 = 64$ complementary combinations can actually be realized in the affine plane.
This is because, a little more generally, only sixteen of the $3^6 = 729$ combinations of $\Tau_k[\msf a, \msf b] = 0$ (with $k = 1, 2$ or $3$ for all pairs $\{ \msf a, \msf b \} \subset \{ \msf A, \msf B, \msf C, \msf D \}$; see Eq.~(\myref[1]{eq:taus})) constitute a consistent system of linear equations connecting the seven areas.
As illustrated in Fig.~\myref[3]{fig:regions}, an example of each of these sixteen combinations may be obtained by fixing the triangle $\GEO{BCD}$, and placing $\GEO{A}$ in one of the sixteen regions into which the plane is divided by the lines through that triangle's edges together with the lines through each of its vertices and parallel to the opposite edge.
To facilitate the ensuing discussion, Fig.~\myref{fig:regions-labeled} of this part of the series again shows these $16$ regions, now labeled by the pool ball icons \textsf{0} through \textsf{15}.

The seven combinations of natural and inverse natural parameters vanishing so as to satisfy the complementarity relations will henceforth be referred to as \textquote{cases,} to distinguish them from the full set of sixteen \textquote{classes,} or combinations of equations $\Tau_k[\msf a, \msf b] = 0$ that are mutually consistent.
In the four (single class) cases realized by placing $\GEO{A}$ in the regions labeled by the pool balls \textsf{0} through \textsf{3}, three non-opposite natural and three non-opposite inverse natural parameters vanish, while in the remaining three cases a pair of opposite natural parameters vanishes along with the four non-complementary inverse natural parameters.
Geometrically, the first four cases are characterized by having one vertex in the convex span of the other three, while all four vertices are extreme points of their convex hull in the remaining three cases. 

In the first four cases it is easily seen that the vanishing of each triple of non-opposite natural parameters ensures that $\Omega(u,v,w,x,y,z) = 0$ (cf.~Eq.~(\myref[2]{eq:myform})), but in the latter three cases we have
\begin{equation} \begin{alignedat}{6}
u ~=~ z ~=~ 0  & ~\implies~ & \Omega(u,v,w,x,y,z) & ~=~ & -\fit(w x - v y)^2 & ~\le~ & 0 ~, \\
v ~=~ y ~=~ 0  &~\implies~ & \Omega(u,v,w,x,y,z) & ~=~ & -\fit(w x - u z)^2 & ~\le~ & 0 ~, \\
w ~=~ x ~=~ 0 & ~\implies~ & \Omega(u,v,w,x,y,z)&  ~=~ & -\fit(v y - u z)^2 & ~\le~ & 0 ~.
\end{alignedat} \end{equation}
This shows that in the latter three cases there are likewise only three natural parameters that can be freely varied without making $\Omega$ negative.
This is consistent with our earlier observation that for rank $1$ tetrahedra the natural parameters do \underline{not} determine a four-point configuration in the plane up to Euclidean isometries, but only up to area-preserving affine transformations.
Since the special affine group of the plane is a five-dimensional Lie group, the space of four-point configurations modulo its action has dimension $2 \cdot 4 - 5 = 3$ in accord with this expectation.

\begin{figure}
\input{Fig_IV-1.tex}
\label{fig:regions-labeled}
\end{figure}

We will now show how each of the sixteen consistent combinations of these equations corresponds to one of the sixteen classes of affine four-point configurations illustrated in Fig.~\myref{fig:regions-labeled}, and that each class can be parameterized by the areas of the three interior faces (better known, in planar quadrilaterals, as Varignon parallelograms).
We do this by expressing the areas of the four exterior faces as signed sums of the areas of the interior faces, where each such quadruple of linear relations is particular to the configuration class in question, and then showing how the barycentric coordinates of $\GEO{A}$ relative to $\GEO{BCD}$ can be obtained from these four (unsigned) areas.

When $u = v = w = 0 = \tilde z = \tilde y = \tilde x$, for example (specifically, the class with $\GEO{A} \in \GEO{BCD}$ as depicted in Fig.~\myref{fig:regions-labeled}), the seven areas satisfy the system of six linear equations
\begin{equation}
\Tau_1[\msf A, \msf B] ~=~ \Tau_1[\msf A, \msf C] ~=~ \Tau_1[\msf A, \msf D] ~=~ \Tau_3[\msf C, \msf D] ~=~ \Tau_3[\msf B, \msf D] ~=~ \Tau_3[\msf B, \msf C] ~=~ 0 ~.
\end{equation}
Written out in terms of matrices, this system of equations is
\begin{equation}
\left[ \begin{smallmatrix} \hit[1.5ex] ~1\>&\>~1\>&\>~0\>&\>~0~\\[5pt] ~1\>&\>~0\>&\>~1\>&\>~0~\\[5pt] ~0\>&\>~1\>&\>~1\>&\>~0~\\[5pt] ~0\>&\>~0\>&\>-1\>&\>~1~\\[5pt] ~0\>&\>-1\>&\>~0\>&\>~1~\\[5pt] -1\>&\>~0\>&\>~0\>&\>~1~\rule[-1ex]{0pt}{1ex} \end{smallmatrix} \right] \hspace{-0.5em} \begin{array}{c} \left[ \begin{smallmatrix} \hit[1ex] \tmv \OL{\msf{ABC}} \tmv \\[2pt] \tmv \OL{\msf{ABD}} \tmv \\[2pt] \tmv \OL{\msf{ACD}} \tmv \\[2pt] \tmv \OL{\msf{BCD}} \tmv \rule[-1ex]{0pt}{1ex} \end{smallmatrix} \right] \\[4ex] \hit[4ex] \end{array} ~=~\>\> 2\, \left[ \begin{smallmatrix} \hit[1.5ex] \tmv \OL{\msf{AB \tmv CD}} \tmv \\[1.2pt] \tmv \OL{\msf{AC \tmv BD}} \tmv \\[1.2pt] \tmv \OL{\msf{AD \tmv BC}} \tmv \\[1.2pt] \tmv \OL{\msf{AB \tmv CD}} \tmv \\[1.2pt] \tmv \OL{\msf{AC \tmv BD}} \tmv \\[1.2pt] \tmv \OL{\msf{AD \tmv BC}} \tmv \rule[-1ex]{0pt}{1ex} \end{smallmatrix} \right] ~.
\end{equation}
The left kernel of this matrix is spanned by $[1,0,-1,-1,0,1]$ \& $[0,1,-1,0,-1,1]$, and since these vectors are orthogonal to the right-hand side, the equations admit the exact solution:
\begin{equation}
\left[ \begin{smallmatrix} \hit[2ex] \tmv \OL{\msf{ABC}} \tmv \\[3pt] \tmv \OL{\msf{ABD}} \tmv \\[3pt] \tmv \OL{\msf{ACD}} \tmv \\[3pt] \tmv \OL{\msf{BCD}} \tmv \rule[-0.75ex]{0pt}{1ex} \end{smallmatrix} \right] ~=~ \left[ \begin{smallmatrix} \hit[1.5ex] \>\; \tmv \OL{\msf{AB \tmv CD}} \tmv \;+\; \tmv \OL{\msf{AC \tmv BD}} \tmv \;-\; \tmv \OL{\msf{AD \tmv BC}}  \tmv ~ \\[2pt] \>\; \tmv \OL{\msf{AB \tmv CD}} \tmv \;-\; \tmv \OL{\msf{AC \tmv BD}} \tmv \;+\; \tmv \OL{\msf{AD \tmv BC}}  \tmv ~ \\[2pt] -\; \tmv \OL{\msf{AB \tmv CD}} \tmv \;+\; \tmv \OL{\msf{AC \tmv BD}} \tmv \;+\; \tmv \OL{\msf{AD \tmv BC}}  \tmv ~ \\[2pt] \>\; \tmv \OL{\msf{AB \tmv CD}} \tmv \;+\; \tmv \OL{\msf{AC \tmv BD}} \tmv \;+\; \tmv \OL{\msf{AD \tmv BC}}  \tmv ~ \rule[-0.75ex]{0pt}{1ex} \end{smallmatrix} \right] ~.
\end{equation}
It is easily shown that these four linear relations among the seven areas imply both $\Xi = 0$ in Eq.~(\myref[1]{eq:yetter}) and $\Omega = 0$ in Eq.~(\myref[2]{eq:myform}).
By performing similar analyses for each of the sixteen classes illustrated in Fig.~\myref{fig:regions-labeled}, one arrives at the sixteen quadruples of linear relations among the areas specified by the first $4$ rows of Table \myref{tab:classes}.

\begin{table} \centering
\caption{
The first four rows give the signed sums $\Upsilon_k$ ($k = 1,2,3$) of the interior areas that equal the exterior areas specified in the leftmost column for each of the $16$ configuration classes illustrated in Fig.~\myref{fig:regions-labeled}, where $\Upsilon_k \coloneq \varepsilon_1\! \MAG{AB \tmv CD} + \varepsilon_2\! \MAG{AC \tmv BD} + \varepsilon_3\! \MAG{AD \tmv BC}$ with $\varepsilon_k = -1$, $\varepsilon_i = \varepsilon_j = +1$ ($\{ i, j, k \} = \{ 1, 2, 3 \}$) and $\Upsilon_0 \coloneq \Upsilon_1 + \Upsilon_2 + \Upsilon_3$.
The last three rows give the signs of the bary\-centric coordinates $\alpha_\msf B$, $\alpha_\msf C$, $\alpha_\msf D$ of $\GEO{A}$ versus the fixed triangle $\GEO{BCD}$ in each configuration class.
The vertical lines of the table separate classes with differing chirotopes (see text).\endgraf}
\smallskip 
\includegraphics[scale=0.8]{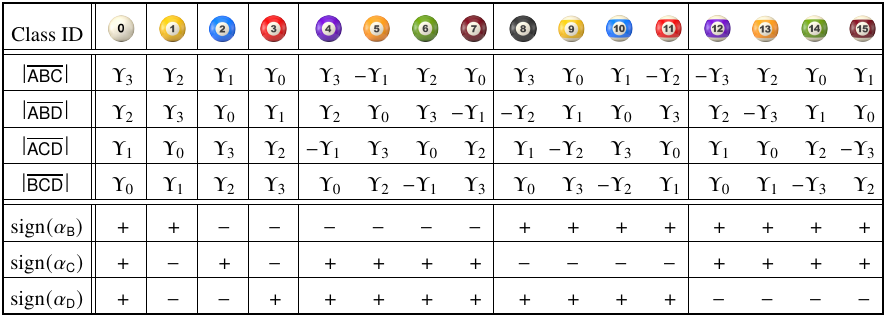} 
\label{tab:classes}
\end{table}

Thus for each configuration class, the exterior areas are given by the signed sum $\Upsilon_k$ ($k = 0,\, \ldots,\, 3$) of the interior, as defined in Table \myref{tab:classes}'s caption, in the corresponding column of that table.
Note that even though the signed sums $\Upsilon_1$, $\Upsilon_2$, $\Upsilon_3$ are not necessarily non-negative (and hence not the deviations from saturation of geometric inequalities), their signs are the same for all configurations in each quadruple of four consecutive columns in the table.
The signed sum that is negative is manifest in the negated $\Upsilon_k$ which appears the last twelve columns of the table.
The last three rows of the table give the signs of the barycentric coordinates $\alpha_\msf B$, $\alpha_\msf C$ \& $\alpha_\msf D$ of $\GEO{A}$ versus $\GEO{BCD}$, the absolute values of which are of course the ratios of the external areas, namely
\begin{equation} \label{eq:abs-bc}
\lvert \alpha_\msf B \rvert ~=~ \frac{\;\MAG{ACD}\;}{\MAG{BCD}} ~,\quad \lvert \alpha_\msf C \rvert ~=~ \frac{\;\MAG{ABD}\;}{\MAG{BCD}} ~,\quad \lvert \alpha_\msf D \rvert ~=~ \frac{\;\MAG{ABC}\;}{\MAG{BCD}} ~.
\end{equation}
Hence these are the signs one must give the absolute values of $\alpha_\msf B$, $\alpha_\msf C$ \& $\alpha_\msf D$, depending on which one of the $16$ configuration classes is being parameterized.

As is well known, the signs of the three barycentric coordinates also determine which of the seven realizable, uniform, four-point, rank $3$ \TDEF{chirotopes\/} (or affine oriented matroids) each class corresponds to \cite{Bjorner:1999, Richter-Gebert:2017}.
Thus these seven chirotopes correspond exactly to the seven cases separated by vertical lines in Table \myref{tab:classes}, which Fig.~\myref{fig:regions-labeled} shows are distinguished by the sets of natural and inverse natural parameters that vanish on them.
The full classification of four-point configurations into sixteen classes thus constitutes a refinement of the cases distinguished by their chirotopes, specifically a division of the three cases covered by the last $12$ columns of Table \myref{tab:classes} into four classes each.
Goodman \& Pollack \cite{Goodman:1983, Goodman:1990} also defined a classification of affine point configurations which is finer than the chirotope one, and which is based upon the concept of \textquote{allowable sequences} of permutations of point labels.
Each permutation is obtained by projecting the points onto an oriented line in the plane, while the sequence of permutations is obtained by rotating the line through $2\pi$ radians; such periodic sequences are uniquely defined by the configuration up to inversion of all the permutations therein and reversal of the overall sequence.

Based on an enumeration of the allowable sequences, again using the \texttt{GeoGebra} dynamic geometry system \cite{GeoGebra:2021}, it has been confirmed that the four classes into which each of the convex chirotope cases is divided (columns \textsf{4--7}, \textsf{8--11} \& \textsf{12--15} in Table \myref{tab:classes}) are distinguished by their allowable sequences.
This is illustrated by the allowable sequences for classes \textsf{4--7} which are shown in Table \myref{tab:permseqs}.
It also appears that those allow\-able sequences are well defined for each class, in that moving the vertex $\GEO{A}$ around within each of the regions separated by dotted lines in Fig.~\myref{fig:regions-labeled} merely changes the permutation to its predecessor or successor in the allowable sequence.
This analysis suggests that more generally there may be a connection between the allowable sequences of affine configurations and the lines parallel to those spanned by pairs of points in the configuration, but through points other than those in that pair.
If so, this may serve to make the \textquote{combinatorial types} distinguished by allowable sequences up to relabeling more amenable to analysis than they currently seem to be, at least in comparison to the better known \textquote{order types} distinguished by chirotopes \cite{Goodman:1993}.

\begin{table} \centering
\caption{
The periodic sequences of allowable permutations which are generated when $\GEO{A}$ is placed in the regions labeled by the pool balls $\msf 4$, $\msf 5$, $\msf 6$ \& $\msf 7$ in Fig.~\myref{fig:regions-labeled}, where each sequence starts from an orthogonal projection of the points onto a line parallel to $\GEO{CD}$ which is rotated up to the first inversion of the starting permutation (see text).
Permutations which differ between the sequences are underlined for emphasis.}
\includegraphics[scale=0.8]{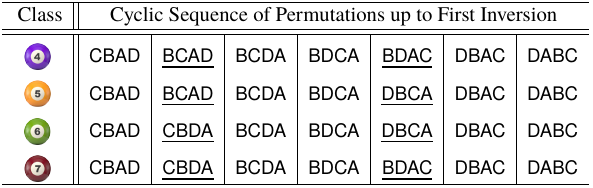}
 \label{tab:permseqs}
\end{table}

\section{Canonical distances for rank 1 configurations} \label{sec:canonical}
Given any values of the natural parameters for which the Gram matrix $\mbf G_\msf A$ has rank $1$, a specific \underline{Euclidean} configuration will be determined generically by any two Euclidean but not affine invariants associated with the planar tetrahedron, for example any two of the distances between its vertices.
It is also possible, however, to specify a \underline{canonical} planar Euclidean tetrahedron which realizes the same seven areas, which in turn are determined by the natural parameters via Corollary \myref[2]{thm:ids}.
This canonical planar tetrahedron is the one consistent with the given natural parameters and associated areas that minimizes the squared radius of gyration of the points about their centroid.
By a well-known theorem of Lagrange \cite{Crippen:1988,Gidea:2012,Havel:2024}, this in turn is equal to $1/16$ times the sum of the six squared distances among the four vertices.
It will now be shown how this problem can be solved by Lagrange's method of undetermined multipliers.
As an interesting by-product, this proves that even though the squared areas cannot be realized by Cayley-Menger and Talata (Eqs. (\myref[1]{eq:cmd3}) \& \linebreak[3] (\myref[1]{eq:talata})) determinants in the general rank $\le 2$ case, they always can be in the rank $1$ case.

To find such distances, consider the three-point instance of \textit{Sch\"onberg's quadratic form\/} \cite{Blumenthal:1953, Crippen:1988, Havel:2024, Schoenberg:1937}:
\begin{equation} \label{eq:schoenberg}
S_\msf{BCD}(\delta_\msf B, \delta_\msf C, \delta_\msf D) ~\coloneq~ -\sfrac12\fit[-0.25em] \begin{array}{c} \left[ \begin{smallmatrix} \hit[1.5ex] \delta_\msf B & \delta_\msf C & \delta_\msf D \dit[1ex] \end{smallmatrix} \right] \\[4ex] \end{array} \fit[-0.5em] \left[ \begin{smallmatrix} \hit[1.5ex] 0 & D_\msf{BC} & D_\msf{BD} \\[2pt] D_\msf{BC} & 0 & D_\msf{CD} \\[2pt] D_\msf{BD} & D_\msf{CD} & 0 \dit[1ex] \end{smallmatrix} \right] \left[ \begin{smallmatrix} \hit[1.5ex] \delta_\msf B \\[2pt] \delta_\msf C \\[2pt] \delta_\msf D \dit[1ex] \end{smallmatrix} \right] .
\end{equation}
It is well known [\textit{op.~cit.\/}] that $D_\msf{BC\,}, D_\msf{BD\,}, D_\msf{CD} \ge 0$ are the squared distances $\MAG{BC}^2\fit[-0.25em],\:$ $\MAG{BD}^2\fit[-0.25em],\:$ $\MAG{CD}^2$ among the vertices of a triangle $\GEO{BCD}$ in the Euclidean plane if \& only if $S_\msf{BCD}(\delta_\msf B, \delta_\msf C, \delta_\msf D)$ $\ge 0$ for all $\delta_\msf B + \delta_\msf C + \delta_\msf D = 0$.
It is however seldom mentioned that, in this case, if \smash{$\GEO{A}$ \& $\GEO{A'}$} are two points in the plane and $\alpha_{\hit[1.5ex]\msf B}, \alpha_{\hit[1.5ex]\msf C}, \alpha_{\hit[1.5ex]\msf D}$ \& $\alpha_\msf B', \alpha_\msf C', \alpha_\msf D'$ are their barycentric coordinates w.r.t.~\smash{$\GEO{BCD}$} ($\alpha_{\hit[1.5ex]\msf B} + \alpha_{\hit[1.5ex]\msf C} + \alpha_{\hit[1.5ex]\msf D} = 1 = \alpha_\msf B' + \alpha_\msf C' + \alpha_\msf D'$), then
\begin{equation}
\delta_\msf B \,=\, \alpha_{\hit[1.5ex]\msf B} -\, \alpha_\msf B' \,,\fit[0.5em] \delta_\msf C \,=\, \alpha_{\hit[1.5ex]\msf C} -\, \alpha_\msf C' \,,\fit[0.5em] \delta_\msf D \,=\, \alpha_{\hit[1.5ex]\msf D} -\, \alpha_\msf D' \>\implies\> \MAG{AA'}^2 =\, S_\msf{BCD}(\delta_\msf B, \delta_\msf C, \delta_\msf D) \,.
\end{equation}
This implies, in particular, that the squared distances from the vertex $\GEO{A}$ to the vertices of $\GEO{BCD}$ itself are given by
\begin{equation} \label{eq:a2bcd} \begin{aligned}
\MAG{AB}^2 ~=~ & S_\msf{BCD} (\alpha_{\msf B} - 1,\, \alpha_{\msf C\,},\, \alpha_{\msf D})  ~=~ S_\msf{BCD}(-\alpha_{\msf C} - \alpha_{\msf D\,},\, \alpha_{\msf C\,},\, \alpha_{\msf D}) ~, \\
\MAG{AC}^2 ~=~ & S_\msf{BCD} (\alpha_{\msf B\,},\, \alpha_{\msf C} - 1,\, \alpha_{\msf D})  ~=~ S_\msf{BCD}(\alpha_{\msf B\,},\, -\alpha_{\msf B} - \alpha_{\msf D\,},\, \alpha_{\msf D}) ~, \\
\MAG{AD}^2 ~=~ & S_\msf{BCD} (\alpha_{\msf B\,},\, \alpha_{\msf C\,},\, \alpha_{\msf D} - 1) ~=~ S_\msf{BCD}(\alpha_{\msf B\,},\, \alpha_{\msf C\,},\, -\alpha_{\msf B} - \alpha_{\msf C}) ~,
\end{aligned} \end{equation}
where the right-hand sides rewrite the left in homogeneous form by substituting for $1 = \alpha_{\hit[1.5ex]\msf B} + \alpha_\msf{\hit[1.5ex]\msf C} + \alpha_{\hit[1.5ex]\msf D}$.
This means that the sum of all six squared distances among the four vertices can be written as a function of the squared distances among the vertices of the triangle $\GEO{BCD}$ alone:
\begin{align}
16\, R_\msf G ~=~ &\> S_\msf{BCD}(-\alpha_{\msf C} - \alpha_{\msf D}, \alpha_{\msf C}, \alpha_{\msf D}) \,+\, S_\msf{BCD}(\alpha_{\msf B}, -\alpha_{\msf B} - \alpha_{\msf D}, \alpha_{\msf D}) \nonumber \\ &\>
+\; S_\msf{BCD}(\alpha_{\msf B}, \alpha_{\msf C}, -\alpha_{\msf B} - \alpha_{\msf C}) \,+ \big( \alpha_{\msf B} + \alpha_{\msf C} + \alpha_{\msf D} \big)^{\!2} \raisebox{1pt}{\large$\big($} \MAG{BC}^2 +\, \MAG{BD}^2 +\, \MAG{CD}^2 \raisebox{1pt}{\large$\big)$} \nonumber \\
=~ &\> 16\, \varrho_\msf{BC}\, \MAG{BC}^2 +\, 16\, \varrho_\msf{BD}\, \MAG{BD}^2 +\, 16\, \varrho_\msf{CD}\, \MAG{CD}^2 ~,
\end{align}
where $R_\msf G = \raisebox{3pt}{\ssmall$\Big($} \MAG{AB}^2 +\cdots+ \MAG{CD}^2 \raisebox{3pt}{\ssmall$\Big)$} \raisebox{2pt}{\small$\bigm/$} 16$ is the squared radius of gyration, and 
\begin{equation} \begin{aligned}
\varrho_\msf{BC} ~\coloneq~ & \tfrac1{16}\, \big( 2\, \alpha_\msf B^2 \,+\, 2\, \alpha_\msf C^2 \,+\, \alpha_\msf D^2 \,+\, \alpha_{\hit[1.5ex]\msf B} \alpha_{\hit[1.5ex]\msf C} \,+\, 3\, \alpha_{\hit[1.5ex]\msf B} \alpha_{\hit[1.5ex]\msf D} \,+\, 3\, \alpha_{\hit[1.5ex]\msf C} \alpha_{\hit[1.5ex]\msf D} \big) ~, \\
\varrho_\msf{BD} ~\coloneq~ & \tfrac1{16}\, \big( 2\, \alpha_\msf B^2 \,+\, \alpha_\msf C^2 \,+\, 2\, \alpha_\msf D^2 \,+\, 3\, \alpha_{\hit[1.5ex]\msf B} \alpha_{\hit[1.5ex]\msf C} \,+\, \alpha_{\hit[1.5ex]\msf B} \alpha_{\hit[1.5ex]\msf D} \,+\, 3\, \alpha_{\hit[1.5ex]\msf C} \alpha_{\hit[1.5ex]\msf D} \big) ~, \\
\varrho_\msf{CD} ~\coloneq~ & \tfrac1{16}\, \big( \alpha_\msf B^2 \,+\, 2\, \alpha_\msf C^2 \,+\, 2\, \alpha_\msf D^2 \,+\, 3\, \alpha_{\hit[1.5ex]\msf B} \alpha_{\hit[1.5ex]\msf C} \,+\, 3\, \alpha_{\hit[1.5ex]\msf B} \alpha_{\hit[1.5ex]\msf D} \,+\, \alpha_{\hit[1.5ex]\msf C} \alpha_{\hit[1.5ex]\msf D} \big) ~. \\
\end{aligned} \end{equation}

Regarding $R_\msf G$ now as a polynomial function of $D_\msf{BC}$, $D_\msf{BD}$ \& $D_\msf{CD}$, consider the Lagrangian for the minimization of $\frac12\, R_\msf G$ subject to the constraint on the three-point Cayley-Menger determinant $\CMD_D[\msf B, \msf C, \msf D] = f_\msf{BCD}^2$, namely
\begin{equation}
L_\lambda(D_\msf{BC}, D_\msf{BD}, D_\msf{CD}) \,\coloneq\, \tfrac12\, \big( \varrho_\msf{BC}\fit D_\msf{BC} + \varrho_\msf{BD}\fit D_\msf{BD} + \varrho_\msf{CD}\fit D_\msf{CD} \big) \,-\, \lambda \big( \CMD_D[\msf B, \msf C, \msf D] - f_\msf{BCD}^2 \big)
\end{equation}
where $\lambda$ is a Lagrange multiplier and $f_\msf{BCD}$ is the value of \smash{$\MAG{BCD}$} (times $2$).
Both this area and the barycentric coordinates are obtained from the natural parameters via Corollary \myref[2]{thm:ids} together with the relations between these coordinates and the areas given in Table \myref{tab:classes} \& Eq.~(\myref{eq:abs-bc}).
Setting the gradient of $L_\lambda$ w.r.t.~$D_\msf{BC}$, $D_\msf{BD}$, $D_\msf{CD}$ to $\mbf 0$ gives
\begin{equation} \label{eq:lgrad}
\mathbf 0 ~=~ \nabla L_\lambda ~=~ \frac12 \left[ \begin{smallmatrix} \hit[1.5ex] \varrho_\msf{BC} \\[4pt] \varrho_\msf{BD} \\[4pt] \varrho_\msf{CD} \dit[1ex] \end{smallmatrix} \right] \,\Vdiff\, \frac{\lambda}2 \left[ \begin{smallmatrix} \hit[1.5ex] D_\msf{BD} \,+\, D_\msf{CD} \,-\, D_\msf{BC} \\[2pt] D_\msf{CD} \,+\, D_\msf{BC} \,-\, D_\msf{BD} \\[2pt] D_\msf{BC} \,+\, D_\msf{BD} \,-\, D_\msf{CD} \dit[1ex] \end{smallmatrix} \right] ,
\end{equation}
and on adding these three equations together and solving for $\lambda$, one obtains
\begin{equation}
\lambda^* ~=~ \frac{\varrho_\msf{BC} + \varrho_\msf{BD} + \varrho_\msf{CD}}{D_\msf{BC} + D_\msf{BD} + D_\msf{CD}} ~.
\end{equation}
Substituting this value of $\lambda$ back into Eq.~(\myref{eq:lgrad}) and clearing denominators by multiplying through by $D_\msf{BC} + D_\msf{BD} + D_\msf{CD}$ then leads to the following linear system of equations:
\begin{equation}
\left[ \begin{smallmatrix} \hit[1.5ex]
2\, \varrho_\msf{BC} \,+\, \varrho_\msf{BD} \,+\, \varrho_\msf{CD} & -\,\varrho_\msf{BD} \,-\, \varrho_\msf{CD} & -\,\varrho_\msf{BD} \,-\, \varrho_\msf{CD} \\[2pt]
-\,\varrho_\msf{BC} \,-\, \varrho_\msf{CD} & \varrho_\msf{BC} \,+\, 2\, \varrho_\msf{BD} \,+\, \varrho_\msf{CD} & -\,\varrho_\msf{BC} \,-\, \varrho_\msf{CD} \\[2pt]
-\,\varrho_\msf{BC} \,-\, \varrho_\msf{BD} & -\,\varrho_\msf{BC} \,-\, \varrho_\msf{BD} & \varrho_\msf{BC} \,+\, \varrho_\msf{BD} \,+\, 2\, \varrho_\msf{CD} \dit[1ex] \end{smallmatrix} \right] \left[ \begin{smallmatrix} \hit[1.5ex] D_\msf{BC} \\[2pt] D_\msf{BD} \\[2pt] D_\msf{CD} \dit[1ex] \end{smallmatrix} \right] ~=~ \left[ \begin{smallmatrix} \hit[1.5ex] 0 \\[4pt] 0 \\[4pt] 0 \dit[1ex] \end{smallmatrix} \right] ~.
\end{equation}
The right-kernel of this matrix is spanned by $[ \varrho_\msf{BD} + \varrho_\msf{CD} \,,~ \varrho_\msf{BC} + \varrho_\msf{CD} \,,~ \varrho_\msf{BC} + \varrho_\msf{BD} ]^\top$, so the optimum squared distances are
\begin{equation}
D_\msf{BC}^* ~=~ \zeta\, (\varrho_\msf{BD} + \varrho_\msf{CD}) ~,\quad D_\msf{BD}^* ~=~ \zeta\, (\varrho_\msf{BC} + \varrho_\msf{CD}) ~,\quad D_\msf{CD}^* ~=~ \zeta\, (\varrho_\msf{BC} + \varrho_\msf{BD})
\end{equation}
for some $\zeta > 0$.
Its value may be found by substituting these values of the squared distances into the constraint, which gives
\begin{equation}
\CMD_D^*[ \msf B, \msf C, \msf D ] ~=~ \zeta^2\, (\varrho_\msf{BC} \varrho_\msf{BD} + \varrho_\msf{BC} \varrho_\msf{CD} + \varrho_\msf{BD} \varrho_\msf{CD}) ~=~ f_\msf{BCD}^2 ~,
\end{equation}
so that $\zeta = f_\msf{BCD}\, / \sqrt{\varrho_\msf{BC} \varrho_\msf{BD} + \varrho_\msf{BC} \varrho_\msf{CD} + \varrho_\msf{BD} \varrho_\msf{CD}}$.
Once the distances within $\GEO{BCD}$ are known, the remaining optimum distances can be computed from Eq.~(\myref{eq:a2bcd}).

Note that even though the process by which this optimum was obtained singled out $\GEO{BCD}$ as a barycentric basis, it remains canonical in that the same squared distances would've been found for any other choice of basis.
Just as importantly, the three-point Cayley-Menger and Talata determinants computed from those squared distances will all match the squares of the areas obtained directly from the natural parameters as above.
Finally, if so desired Cartesian coordinates which realize these squared distances can be computed via standard eigenvalue methods \cite{Crippen:1978,Crippen:1988}, and the areas computed from these coordinates will also match those obtained directly from the natural parameters.
It follows that in the $\mrm{rank}\big( \mbf G_\msf A \big) = 1$ case the natural parameters not only determine a four-point planar configuration up to special affine transformation, but also determine a unique Euclidean configuration which minimizes the squared radius of gyration $R_\msf G$ subject to reproducing the areas as calculated from those natural parameters.

Because the set of all sets of four points in the Euclidean plane modulo isometries is five-dimensional just like the set of all tetrahedra with $\Omega = 0$ (cf.~Fig.~\myref[3]{fig:lattice}), it is tempting to speculate that the above canonical map can be extended to a canonical bijection between these two five-dimensional sets.
This would be much more in accord with the intuition, which has been taken for granted by scientists and mathematicians throughout history, that the boundary of the six-dimensional set of non-degenerate tetrahedra is the set of all quadruples of points in the Euclidean plane.
This possibility seems interesting enough to be stated formally as a conjecture:
\begin{conjecture} \label{thm:canmap}
The bijection between the set of all quadruples in the special affine plane and those quadruples in the Euclidean plane for which the radius of gyration attains its unique minimum, subject to maintaining the areas of the triangles and Varignon parallelograms in the affine quadruple, can be extended to a bijection of the set of all degenerate tetrahedra, as defined by the zeros of $\Omega$, with the set of all quadruples in the Euclidean plane; moreover this mapping is or can be chosen to be canonical, so these two sets can be identified.
\end{conjecture}

One approach to proving this conjecture starts from the construction of the zeros of $\Omega$ as the limit of affine transformations $\lim_{\sigma \rightarrow \infty} \mathcal A_\sigma$ applied to a random tetrahedron, as described in Sec.~\myref[3]{sec:zeros}.
Observe that this limit is translation independent, and that $\mathcal P \coloneq \lim_{\sigma \rightarrow 0} \sigma \mathcal A_\sigma$ is the projection of that tetrahedron onto the plane orthogonal to the line which the zeros lie on.
Up to isometry, these limits are also independent of rotations about that line, so a tetrahedron together with its orientation w.r.t.~the line has eight degrees of freedom.
Recalling that the generic zeros have five degrees of freedom, the set of such oriented tetrahedra which go to the same zero of $\Omega$ as $\sigma \rightarrow \infty$ should have three generic degrees of freedom, as should its image under the projection $\mathcal P$.
The conjecture would be proved by constructing the unique planar configuration in that image which minimized the radius of gyration over all such projections for any given zero of $\Omega$, with those from the rank $1$ zeros satisfying Eq. (\ref{eq:lgrad}) above.

An alternative approach might be found by studying the critical points of the sum of squared deviations of the Cayley-Menger $\CMD_D[\msf a, \msf b, \msf c]$ \& Talata $\CMD_D[\msf a, \msf b \,|\, \msf c, \msf d]$ determinants from the squared areas in the zeros of $\Omega$, and showing that exactly one of them corresponds to a four-point configuration in the (finite) Euclidean plane.
In Appendix \myref{sec:pmap} below, we show that these critical points are uniquely determined by any given values of the seven areas in a non-degenerate (fully three-dimensional) tetrahedron.

\section{\hspace*{-7pt}The polynomial map from squared distances to squared areas} \label{sec:pmap}
\makebox[4.8em][r]{Theorem \myref[1]{thm:bdsmc}} showed that the squared areas of a non-degenerate tetrahedron \smash{$\GEO{ABCD}$} determine it uniquely up to isometry.
This appendix will prove a stronger result, namely that any vector $\mbf f \in \mbb R^7$ of \textquote{squared areas} is in the range of the polynomial map from \textquote{squared distances} $\mbf d \in \mbb R^6$ to squared areas (as defined by Eqs.~(\myref[1]{eq:cmd3}) \& (\myref[1]{eq:talata})) or of its negative if \& only if the squared areas satisfy Yetter's identity $\breve\Xi_{\,\mbf f} = 0$ and have a Gramian at $\msf A$ (or any other vertex) $\Gamma_\mbf f[\msf A] \ne 0$.
Moreover, it is in the range of this poly\-nomial map precisely when $\Gamma_\mbf f[\msf A] > 0$, or in the range of its negative when $\Gamma_\mbf f[\msf A] < 0$.

\begin{theorem} \label{thm:pmap}
Let $\mbf p_\pm$ be the two quadratic polynomial mappings from the semi-algebraic set $\mathcal D \coloneq \{\, \mbf d \in \mbb R^6 \mid \CMD_{\,\mbf d\,}[\msf A, \msf B, \msf C, \msf D] \ne 0 \,\}$ into $\mbb R^7$ that are given by
\vspace{-0.5ex}
\begin{multline}
\mbf p_\pm( \mbf d ) ~\coloneq~ \mbi\pm \big[ \CMD_{\,\mbf d\,}[\msf A, \msf B, \msf C]\,,\, \CMD_{\,\mbf d\,}[\msf A, \msf B, \msf D]\,,\,  \CMD_{\,\mbf d\,}[\msf A, \msf C, \msf D]\,,\, \CMD_{\,\mbf d\,}[\msf B, \msf C, \msf D]\,, \\
\CMD_{\,\mbf d\,}[\msf A, \msf B \, \tmv\, \msf C, \msf D]\,,\, \CMD_{\,\mbf d\,}[\msf A, \msf C \, \tmv\, \msf B, \msf D]\,,\, \CMD_{\,\mbf d\,}[\msf A, \msf D \, \tmv\, \msf B, \msf C]\, \big]^{\!\top} . \vspace{-0.5ex}
\end{multline}
Then if Yetter's identity $\breve\Xi_{\,\mbf f} = 0$ is satisfied by the given squared areas $\mbf f = [ F_\msf{ABC}\,,$ $\ldots\,,$ $F_\msf{AD|BC}]^\top$ $\in \mbb R^7$, there exist squared distances $\mbf d = [D_\msf{AB}\,, \ldots\,, D_\msf{CD}]^\top \in \mathcal D$ such that
\vspace{-0.5ex}
\begin{equation} 
\mbf f ~=~ \mbf p_+(\mbf d) \quad\text{if}\quad \Gamma_\mbf f[\msf A] ~>~ 0 ~,\quad\text{or}\quad \mbf f ~=~ \mbf p_-(\mbf d) \quad\text{if}\quad \Gamma_\mbf f[\msf A] ~<~ 0 ~, \vspace{-0.5ex}
\end{equation}
where $\Gamma_\mbf f[\msf A]$ is the Gramian at $\msf A\fit$, and otherwise $\mbf p_\pm(\mbf d) \ne \mbf f$ $\,\forall\,$ $\mbf d \in \mathcal D$.
In other words, the semi-algebraic set $\{\, \mbf f \in \mbb R^7 \mid \Gamma_\mbf f[\msf A] \ne 0\fit,\, \breve\Xi_{\,\mbf f} = 0 \,\} = \mbf p_+( \mathcal D ) \fit[0.25em]\dot{\cup}\fit[0.25em] \mbf p_-(\mathcal D)$ where \textquote{$\fit[0.2em] \dot{\cup}$} is the disjoint union of these two images of $\mathcal D$, and each component is explicitly parameterized by the corresponding mapping $\mbf p_{\pm\,}$.
\end{theorem}
\begin{proof} 
Let $\mathbf Q_\msf{\,abc}\!$ \& $\mathbf Q_\msf{\,ab|cd}$ $\big( \msf a, \msf b, \msf c, \msf d \in \{ \msf A, \msf B, \msf C, \msf D \}$ $\mid$ $|\{\msf a, \msf b, \msf c, \msf d\}|$ $= 4 \big)$ be the matrices of the quadratic forms defined by the Cayley-Menger \& Talata determinants relative to the order of the squared distances in $\mbf d$, e.g.
\begin{equation}
\mathbf Q_\msf{\,BCD} ~\coloneq~ \sfrac14 \left[ \begin{smallmatrix} \rule{0pt}{1.5ex} 0&0&~0&~0&~0&~0 \\[2pt] 0&~0&~0&~0&~0&~0 \\[2pt] 0&~0&~0&~0&~0&~0 \\[2pt] 0&~0&~0&-1&~1&~1 \\[2pt] 0&~0&~0&~1&-1&~1 \\[2pt] 0&~0&~0&~1&~1&-1 \rule[-1ex]{0pt}{1ex} \end{smallmatrix} \right] ~,\qquad
\mathbf Q_\msf{\,AB|CD} ~\coloneq~ \sfrac14 \left[ \begin{smallmatrix} \rule{0pt}{1.5ex} 0&~0&~0&~0&~0&~2 \\[2pt] 0&-1&~1&~1&-1&~0 \\[2pt] 0&~1&-1&-1&~1&~0 \\[2pt] 0&~1&-1&-1&~1&~0 \\[2pt] 0&-1&~1&~1&-1&~0 \\[2pt] 2&~0&~0&~0&~0&~0 \rule[-1ex]{0pt}{1ex} \end{smallmatrix} \right] ~.
\end{equation}
Then the vector-valued functions $\mbf p_\pm$ can be written as
\vspace{-0.5ex}
\begin{multline}
\mbf p_\pm (\mbf d) ~=~ \mbi\pm \big[ \mathbf d^\top \mathbf Q_\msf{ABC}\, \mathbf d \,,~ \mathbf d^\top \mathbf Q_\msf{ABD}\, \mathbf d \,,~ \mathbf d^\top \mathbf Q_\msf{ACD}\, \mathbf d \,,~ \mathbf d^\top \mathbf Q_\msf{BCD}\, \mathbf d \,,~ \\ 
\mathbf d^\top \mathbf Q_\msf{AB|CD}\, \mathbf d \,,~ \mathbf d^\top \mathbf Q_\msf{AC|BD}\, \mathbf d \,,~ \mathbf d^\top \mathbf Q_\msf{AD|BC}\, \mathbf d \big]^{\!\top} ~.
\end{multline}
Now consider the pair of least-squares problems:
\begin{equation}
\min_{\mbf d \in \mathcal D} \big( \tfrac12\, {\Sigma_\pm^2( \mathbf d )} \big) ~,\quad \text{where} \quad \Sigma_\pm^2( \mathbf d ) ~\coloneq~ \|\fit \mbf p_\pm( \mbf d ) \,\Vdiff\, \mbf f \fit\|^2 \qquad \big( \mbf d \in \mathcal D \big) .
\end{equation}
The Jacobians $\mbf J_\pm = \mbf J_\pm(\mbf d)$ of $\tfrac12\, \Sigma_\pm^2$ will have the vectors $\mbi\pm \mbf d^\top \mbf Q_\msf{ABC}$, $\mbi\pm \mbf d^\top \mbf Q_\msf{ABD}$, $\mbi\pm \mbf d^\top \mbf Q_\msf{ACD}$ \& $\mbi\pm \mbf d^\top \mbf Q_\msf{BCD}$ as their first four rows, and $\mbi\pm \mbf d^\top \mbf Q_\msf{AB|CD}$,  $\mbi\pm \mbf d^\top \mbf Q_\msf{AC|BD}$ \& $\mbi\pm \mbf d^\top \mbf Q_\msf{AD|BC}$ as their last three rows.
As is well-known \cite{Bjorck:1996}, the change in squared distances $\Del \mbf d$ obtained by applying the Gauss-Newton method to these least-squares problems satisfies the normal equations
\begin{equation}
\mbf J_\pm^\top \mbf J_{\hit[1.4ex]\pm}\, \Del \mbf d ~=~ \Vdiff \mbf J_\pm^\top\, \Del \mathbi \Sigma_{\hit[1.4ex]\pm} ~,
\end{equation}
where $\Del \mbi\Sigma_\pm = \Del \mbi\Sigma_\pm(\mbf d) \coloneq \mbf p_\pm(\mbf d) \Vdiff \mbf f$ are the vectors of residuals.
The right-hand side is the negative gradient of each sum-of-squares $\tfrac12 \nabla \Sigma_\pm^2\fit$, and it will vanish if \& only if $\Sigma_\pm^2(\mbf d) = 0$ or $\Del \mathbi \Sigma_\pm (\mbf d)$ lies in the generically one-dimensional right null space of $\mbf J_\pm^\top$.
Written out in full, these two matrices are given by:
\begin{multline}
\mbi\pm 2\, \mbf J_\pm^\top(\mbf d) ~=~ \left[ \begin{smallmatrix} \rule{0pt}{1.5ex}
-D_\msf{AB} \fit+\fit D_\msf{AC} \fit+\fit D_\msf{BC} & -D_\msf{AB} \fit+\fit D_\msf{AD} \fit+\fit D_\msf{BD} & 0 & 0 & \cdots \\[2pt]
\fit[0.5em]D_\msf{AB} \fit-\fit D_\msf{AC} \fit+\fit D_\msf{BC} & 0 & -D_\msf{AC} \fit+\fit D_\msf{AD} \fit+\fit D_\msf{CD} & 0 & \cdots \\[2pt]
0 & D_\msf{AB} \fit-\fit D_\msf{AD} \fit+\fit D_\msf{BD} & \fit[0.5em]D_\msf{AC} \fit-\fit D_\msf{AD} \fit+\fit D_\msf{CD} & 0 & \cdots \\[2pt]
\fit[0.5em]D_\msf{AB} \fit+\fit D_\msf{AC} \fit-\fit D_\msf{BC} & 0 & 0 & -D_\msf{BC} \fit+\fit D_\msf{BD} \fit+\fit D_\msf{CD} & \cdots \\[2pt]
0 & \fit[0.5em]D_\msf{AB} \fit+\fit D_\msf{AD} \fit-\fit D_\msf{BD} & 0 & \fit[0.5em]D_\msf{BC} \fit-\fit D_\msf{BD} \fit+\fit D_\msf{CD} & \cdots \\[2pt]
0 & 0 & \fit[0.5em]D_\msf{AC} \fit+\fit D_\msf{AD} \fit-\fit D_\msf{CD} & \fit[0.5em]D_\msf{BC} \fit+\fit D_\msf{BD} \fit-\fit D_\msf{CD} & \cdots
\rule[-1ex]{0pt}{1ex} \end{smallmatrix} \right. \\ 
\left. \begin{smallmatrix} \rule{0pt}{1.5ex}
\cdots & 2\, D_\msf{CD} & -D_\msf{AB} \fit+\fit D_\msf{AD} \fit+\fit D_\msf{BC} \fit-\fit D_\msf{CD} & -D_\msf{AB} \fit+\fit D_\msf{AC} \fit+\fit D_\msf{BD} \fit-\fit D_\msf{CD} \\[2pt]
\cdots & -D_\msf{AC} \fit+\fit D_\msf{AD} \fit+\fit D_\msf{BC} \fit-\fit D_\msf{BD} & 2\, D_\msf{BD} & \fit[0.5em]D_\msf{AB} \fit-\fit D_\msf{AC} \fit-\fit D_\msf{BD} \fit+\fit D_\msf{CD} \\[2pt]
\cdots & \fit[0.5em]D_\msf{AC} \fit-\fit D_\msf{AD} \fit-\fit D_\msf{BC} \fit+\fit D_\msf{BD} &\fit[0.5em]D_\msf{AB} \fit-\fit D_\msf{AD} \fit-\fit D_\msf{BC} \fit+\fit D_\msf{CD} & 2\, D_\msf{BC} \\[2pt]
\cdots & \fit[0.5em]D_\msf{AC} \fit-\fit D_\msf{AD} \fit-\fit D_\msf{BC} \fit+\fit D_\msf{BD} & \fit[0.5em]D_\msf{AB} \fit-\fit D_\msf{AD} \fit-\fit D_\msf{BC} \fit+\fit D_\msf{CD} & 2\, D_\msf{AD} \\[2pt]
\cdots & -D_\msf{AC} \fit+\fit D_\msf{AD} \fit+\fit D_\msf{BC} \fit-\fit D_\msf{BD} & 2\, D_\msf{AC} & \fit[0.5em]D_\msf{AB} \fit-\fit D_\msf{AC} \fit-\fit D_\msf{BD} \fit+\fit D_\msf{CD} \\[2pt]
\cdots & 2\, D_\msf{AB} & -D_\msf{AB} \fit+\fit D_\msf{AD} \fit+\fit D_\msf{BC} \fit-\fit D_\msf{CD} & -D_\msf{AB} \fit+\fit D_\msf{AC} \fit+\fit D_\msf{BD} \fit-\fit D_\msf{CD}
\rule[-1ex]{0pt}{1ex} \end{smallmatrix} \right] ~.
\end{multline}
\newcommand*{\dpmst}{\ensuremath\mbf d_{\raisebox{1.25pt}[0pt][-3pt]{$\scriptstyle\pm$}}^*}%
\newcommand*{\SigpmSq}{\ensuremath\Sigma_{\raisebox{1.25pt}[0pt][-3pt]{$\scriptstyle\pm$}}^2}%
\newcommand*{\JacpmT}{\ensuremath\mbf J_{\raisebox{1.25pt}[0pt][-3pt]{$\scriptstyle\pm$}}^\top}%
Using computer algebra, it is readily shown that $\mrm{det}\big( \JacpmT(\mbf d)\, \mbf J_\pm (\mbf d) \big) = 28\, \CMD_{\,\mbf d\,}[ \msf A, \msf B, \msf C, \msf D ]^4$ $\ne$ $0$, so the Jacobians $\mbf J_\pm$ are of full rank $6$ and the mappings $\mbf p_{\pm\,}: \mathcal D \rightarrow \mbb R^7$ are injective.
It is also easily seen that the vector $\mbf n \coloneq$ $[ 1, ~1, ~1, ~1, -1, -1, -1 ]^\top$ always lies in these Jacobians' left null space, meaning that $\JacpmT( \mbf d )\, \mbf n = \mbf 0$ for all $\mbf d \in \mathcal D$.
Thus the residual vector $\Del \mathbi \Sigma_\pm$ at any critical point $\dpmst$ of $\SigpmSq$ can be written as $\sigma_\pm\, \mbf n$ for some $\sigma_\pm \in \mbb R$, but it can also be written as $\Del \mathbi \Sigma_\pm( \dpmst ) \fit=\fit \mbf p_\pm( \dpmst )  \fit\Vdiff\fit \mbf f$.
By expanding the corresponding Cayley-Menger \& Talata determinants, it is easily shown that $\mbf p_\pm( \mbf d )$ always satisfies Yetter's identity, i.e.~$\mbf n^\top \mbf p_\pm = 0$, and hence
\begin{equation}
\Sigma_\pm^2 ~=~ \big\|\, \Del \mathbi\Sigma_{\hit[1.3ex]\pm}  \big\|^2 ~=~ \mbf n^\top \mbf n\, \sigma_\pm^2 ~=~ 7\, \sigma_\pm^2 ~=~ \sigma_{\hit[1.3ex]\pm}\, \mbf n^\top\, \big( \mbf p_{\hit[1.3ex]\pm} (\mbf d_\pm^*) \,\Vdiff\, \mbf f \big) ~=~  -\sigma_{\hit[1.3ex]\pm}\, \mbf n^\top\, \mbf f ~.
\end{equation}
It follows that the residual vanishes at every critical point of $\SigpmSq$ if \& only if $\mbf n^\top\, \mbf f = 0$, i.e.~the given squared areas $\mbf f$ satisfy Yetter's identity $\breve\Xi_{\,\mbf f} = 0$.

This shows that when $\breve\Xi_{\,\mbf f} = 0$ every critical point $\dpmst$ of $\SigpmSq$ is a global minimum with value $\SigpmSq(\dpmst) = 0$, but it does not establish the existence of a critical point, and indeed it is easy to find examples of $\mathbf f$'s for which numerical evidence indicates that either $\Sigma_+$ or $\Sigma_-$ has no critical points, in that minimization of $\Sigma_+^2$ or $\Sigma_-^2$ the drives entries in $\mathbf d$ towards $\pm\infty$.
It will now be shown, however, that for every $\mbf f \in \mathbb R^7$ with $\mbf n^\top \mbf f = 0$ and $\Gamma_{\mbf f\,}[\msf A] \ne 0$ exactly one of the functions $\Sigma_+^2$ or $\Sigma_-^2$ defined by $\mbf f$ has a critical point.
Towards that end, first note that by Euler's theorem on homogeneous functions,
\begin{equation}
\mbf p_\pm(\mbf d) ~=~ \tfrac12\, \mbf J_\pm(\mbf d)\, \mbf d ~=~ \mbi\pm \tfrac12\, \mbf J_+(\mbf d)\, \mbf d
\end{equation}
for all $\mbf d \in \mathcal D$, and since every critical point $\dpmst$ of $\SigpmSq$ achieves a residual $\Del\Sigma_\pm( \dpmst ) = \mbf 0$,
\begin{equation} \label{eq:euler}
\mbf f ~=~ \mbf p_\pm^* ~\coloneq~ \mbf p_{\hit[1.4ex]\pm}(\mbf d_\pm^*) ~=~ \tfrac12\, \mbf J_{\hit[1.3ex]\pm}(\mbf d_\pm^*)\, \mbf d_\pm^* ~\eqcolon~ \tfrac12\, \mbf J_\pm^*\, \mbf d_\pm^* ~.
\end{equation}
Next, we denote the \textquote{three-point} Cayley-Menger determinants in the squared \underline{areas} by
\begin{equation} \label{eq:acmd}
\CMD_{\fit\mbf f\fit}[\msf a, \msf b] ~\coloneq~ \tfrac12 \big( F_{\hit[1.3ex]\msf{abc}} F_{\hit[1.3ex]\msf{abd}} + F_{\hit[1.3ex]\msf{abc}} F_{\hit[1.3ex]\msf{ab|cd}}  + F_{\hit[1.3ex]\msf{abd}} F_{\hit[1.3ex]\msf{ab|cd}} \big) -\, \tfrac14 \big( F_\msf{abc}^2 + F_\msf{abd}^2 + F_\msf{ab|cd}^2 \big)
\end{equation}
for $\{\msf a, \msf b, \msf c, \msf d\} = \{\msf A, \msf B, \msf C, \msf D\}$.
Let $\mathbf Q_\msf{\,ab}$ be the $7\times7$ matrices of these six quadratic forms versus the order of the squared areas in $\mbf f$, e.g.
\begin{equation}
\mathbf Q_\msf{\,AB} ~\coloneq~ \sfrac14 \left[ \begin{smallmatrix} \rule{0pt}{1.5ex}
-1&~1&~0&~0&~1&~0&~0 \\[2pt] ~1&-1&~0&~0&~1&~0&~0 \\[2pt] ~0&~0&~0&~0&~0&~0&~0 \\[2pt]  ~0&~0&~0&~0&~0&~0&~0 \\[2pt] ~1&~1&~0&~0&-1&~0&~0 \\[2pt] ~0&~0&~0&~0&~0&~0&~0 \\[2pt]  ~0&~0&~0&~0&~0&~0&~0
\rule[-1ex]{0pt}{1ex} \end{smallmatrix} \right] ,~\ldots,~\> \mathbf Q_\msf{\,CD} ~\coloneq~ \sfrac14 \left[ \begin{smallmatrix} \rule{0pt}{1.5ex}
~0&~0&~0&~0&~0&~0&~0 \\[2pt]  ~0&~0&~0&~0&~0&~0&~0 \\[2pt] ~0&~0&-1&~1&~1&~0&~0 \\[2pt] ~0&~0&~1&-1&~1&~0&~0 \\[2pt] ~0&~0&~1&~1&-1&~0&~0 \\[2pt] ~0&~0&~0&~0&~0&~0&~0 \\[2pt]  ~0&~0&~0&~0&~0&~0&~0
\rule[-1ex]{0pt}{1ex} \end{smallmatrix} \right] ~.
\end{equation}
At a critical point $\dpmst$ these quadratic forms are equal to
\begin{equation} \label{eq:fdok}
\mbf f^\top \fit\mbf Q_\msf{\,ab\hit[1.4ex]}\, \mbf f ~=~ \mbf p_\pm^{*\top} \mbf Q_\msf{\,ab\hit[1.4ex]}\, \mbf p_\pm^* ~=~ \tfrac14\, \mbf d_\pm^{*\top} \mbf J_\pm^{*\top} \mbf Q_\msf{\,ab\hit[1.4ex]}\, \mbf J_\pm^*\, \mbf d_\pm^* ~=~ D_\msf{ab}^*\, \CMD_{\,\hit[1.5ex]\mbf d_{\smash{\pm}}^*}[\msf A, \msf B, \msf C, \msf D] ~,
\end{equation}
where $D_{\smash{\msf{ab}}}^*$ is the associated entry of $\dpmst$ and the last equality may readily be proven by applying any computer algebra program to the definitions (cf.~Eq.~(\myref[1]{eq:sines})).
Given any $\delta \ne 0$, we now set $D_\msf{ab} \coloneq \delta\, \mathbf f^\top \mathbf Q_\msf{\,ab}\, \mathbf f$ for $\msf a, \msf b \in \{ \msf A, \msf B, \msf C, \msf D \}$ ($\msf a \ne \msf b$), and denote the result of substituting these values of $D_\msf{ab}$ in $\CMD_{\,\mbf d\,}[\msf A, \msf B, \msf C, \msf D]$ by $\CMD_{\delta;\mathbf f}[\msf A, \msf B, \msf C, \msf D]$.
Then it is easily seen that equations (\myref{eq:fdok}) will be satisfied when
\begin{equation}
\delta ~=~ \delta^* ~\coloneq~ \frac{1}{\sqrt[4\,]{\CMD_{1;\, \mbf f\fit}[\msf A, \msf B, \msf C, \msf D]\,}} ~.
\end{equation}
which of course requires that $\CMD_{1;\, \mbf f\fit}[\msf A, \msf B, \msf C, \msf D] > 0$.
It turns out, however, that if one eliminates $F_\msf{BCD}$ from $\CMD_{1;\, \mbf f\fit}[\msf A, \msf B, \msf C, \msf D]$ using Yetter's identity $F_\msf{BCD} = F_\msf{AB|CD} + F_\msf{AC|BD} + F_\msf{AD|BC} - F_\msf{ABC} - F_\msf{ABD} - F_\msf{ACD}$, it factorizes to $\Gamma_{\mbf f\fit} [\msf A]^2$.
Therefore $\CMD_{1;\, \mbf f\fit}[\msf A, \msf B, \msf C, \msf D]$ is indeed strictly positive by hypothesis, and we have satisfied the necessary conditions for a critical point in Eq.~(\myref{eq:fdok}).

Finally, upon setting the entries of $\mathbf d^*$ to $D_{\smash{\msf{ab}}}^{\,*} = \delta^*\, \mathbf f^\top \mathbf Q_\msf{\,ab\,} \mathbf f$ in the given order, and again eliminating $F_\msf{BCD}$ using Yetter's identity as above, one finds that the Cayley-Menger \& Talata determinants factorize as
\begin{equation} \label{eq:critpt}
\mbi\pm \mbf p_\pm( \mbf d^* ) ~=~ \mbf p_+(\mathbf d^*) ~=~ {(\delta^*)}^2\, \Gamma_{\mbf f\fit}[\msf A]\, \mathbf f ~.
\end{equation}
Since $\big|\fit \Gamma_{\mbf f\,}[\msf A] \fit\big| = \CMD_{1; \mbf f\,}[\msf A, \msf B, \msf C, \msf D]^{1/2} = (\delta^*)^{-2}$, it follows that the defining equation (\myref{eq:euler}) for a critical point $\mbf p_\pm(\mbf d^*) = \mbf f$ is satisfied by $\mbf p_+$ when $\Gamma_{\mbf f\,} [\msf A] > 0$, or by $\mbf p_-$ when $\Gamma_{\mbf f\,}[\msf A] < 0$, as desired.
\end{proof} 

\begin{remark}
This proof also shows that the maps $\mbf p_\pm$ are diffeomorphisms onto their respective images.
Because $\Sigma_\pm^2(\mbf d) = \Sigma_\pm^2(\Vdiff\mbf d)$ for all $\mbf d \in \mathcal D$, the functions $\Sigma_\pm^2$ would have $\mbf 0$ as a saddle point but for $\mbf 0 \notin \mathcal D$.
Other critical points not in $\mathcal D$ are easy to construct simply by setting the entries of $\mbf d$ to the squared distances in a random planar tetrahedron, finding a left null vector $\mbf n_\perp$ of $\mbf J_\pm(\mbf d)$ orthogonal to $\mbf n$, and setting $\mbf f_\pm = \mbf p_\pm(\mbf d) \Vplus \nu\, \mbf n_\perp$ for any $\nu \in \mbb R$.
Such vectors $\mbf f_\pm$ will satisfy $\mbf n^\top \mbf f _\pm = 0$ \& $\Gamma_{\mbf f_\pm}[\msf A] \ne 0$ unless $\nu = 0$, and $\mbf d$ will generally also be a saddle point of $\Sigma_\pm^2$ extended to $\mbb R^6$.
Although most vectors $\mbf f \in \mbb R^7$ with $\mbf n^\top \mbf f = 0$ \& $\Gamma_{\,\mbf f\,}[\msf A] = 0$ will not be equal to the Cayley-Menger \& Talata determinants for any $\mbf d \in \mbb R^6$, it was shown in Appendix \myref{sec:canonical} that they are when the entries of $\mbf f$ are the squared areas of a tetrahedron in the (finite) Euclidean plane.
\end{remark}

\section{Towards a 2-to-2 parametrization of the generic zeros} \label{sec:2to2}
This appendix presents the details of the $2\!:\!2$ parametrization of the \textquote{generic} zeros of $\Omega$ that was briefly described in Remark \myref[3]{rem:cubic}.
It is based upon the parameters $\alpha, \beta, \gamma$ \& $\delta \in \mbb R$ introduced in that remark together with the exterior surface area $s$ which, when treated in this way as an independent parameter, will be denoted by $\varsigma$.
Here, \textquote{generic} means merely that $\alpha^2\!$, $\beta^2\!$, $\gamma^2\!$ \& $\delta^2$ are all non-zero and distinct, although they are subject to the constraint $\alpha\beta\gamma\delta > 0$ (see Proposition \myref[3]{thm:collinear}) along with $|\alpha| + |\beta| + |\gamma| + |\delta| > |\alpha + \beta + \gamma + \delta |$ (so they do not all have the same sign).
In addition, $\varsigma$ must exceed a positive lower bound determined by $\alpha, \beta, \gamma, \delta$ (to be derived below).

Before diving into the intricacies of this matter, however, we will show how the squared magnitudes of the \textquote{$1\mrm D$ vectors} represented by the first four parameters, namely
\begin{equation} \label{eq:abcd-defs}
\alpha^2 ~\coloneq~ 2\, uvw / s ~,\quad \beta^2 ~\coloneq~ 2\, uxy / s ~,\quad \gamma^2 ~\coloneq~ 2\, vxz / s ~,\quad \delta^2 ~\coloneq~ 2\, wyz / s ~,
\end{equation}
can be interpreted geometrically in non-degenerate tetrahedra.
\begin{proposition} \label{thm:itp2vtx}
The three equal squared distances from each vertex of a non-degen\-erate tetrahedron $\GEO{ABCD}$ to the three of the in-touch points $\GEO{J}, \GEO{K}, \GEO{L}, \GEO{N}$ which lie on the three exterior faces meeting in that vertex are
\begin{equation} \begin{aligned}
\MAG{AN}^2 \,=~ \MAG{AL}^2 \,=~ \MAG{AK}^2 \,=~ \alpha^2 / r^2 \,,&\quad
\MAG{BN}^2 \,=~ \MAG{BL}^2 \,=~ \MAG{BJ}^2 \,=~ \beta^2 / r^2 \,, \\
\MAG{CN}^2 \,=~ \MAG{CK}^2 \,=~ \MAG{CJ}^2 \,=~ \gamma^2 / r^2 \,,&\quad
\MAG{DL}^2 \,=~ \MAG{DK}^2 \,=~ \MAG{DJ}^2 \,=~ \delta^2 / r^2 \,,
\end{aligned} \end{equation}
where $s$ is twice the surface area and $r = t/s$ is the in-radius as usual.
\end{proposition}
\begin{proof}
As is well known and easily shown, the barycentric coordinates of the in-center $\GEO{\,I\,} \in \GEO{ABCD}$ are the ratios of the exterior areas to the total exterior surface area $s/2$.
Hence the squared distance from e.g.~vertex $\GEO{A}$ to $\GEO{\,I\,}$ may be obtained from Sch\"onberg's quadratic form (cf.~Eq.~(\myref{eq:schoenberg}) above) as:
\begin{equation}
\MAG{A\fit I\fit}^2 ~=~ \frac{-1}{2\, s^2} \left[ \begin{smallmatrix} \hit[2ex]
2 \MAG{BCD} \,-\, s \\[3pt] 2 \MAG{ACD} \\[3pt] 2 \MAG{ABD} \\[3pt] 2 \MAG{ABC}
\dit[1ex] \end{smallmatrix} \right]^\top \left[ \begin{smallmatrix} \hit[2ex]
0 & \MAG{AB}^2 & \MAG{AC}^2 & \MAG{AD}^2 \\[2pt]
\MAG{AB}^2 & 0 & \MAG{BC}^2 & \MAG{BD}^2 \\[2pt]
\MAG{AC}^2 & \MAG{BC}^2 & 0 & \MAG{CD}^2 \\[2pt]
\MAG{AD}^2 & \MAG{BD}^2 & \MAG{CD}^2 & 0
\dit[1ex] \end{smallmatrix} \right] \left[ \begin{smallmatrix} \hit[2ex]
2 \MAG{BCD} \,-\, s \\[3pt] 2 \MAG{ACD} \\[3pt] 2 \MAG{ABD} \\[3pt] 2 \MAG{ABC}
\dit[1ex] \end{smallmatrix} \right]
\end{equation}
On substituting for the squared distances on the right-hand side using Eq.~(\myref[2]{eq:dsq-from-uu2}) and for the areas using Eq.~(\myref[2]{eq:idb}), this squared distance becomes simply
\begin{equation}
\MAG{A\fit I\fit}^2 ~=~ \frac{2\, s\, uvw \,+\, \Omega(u,\, \ldots\,,\, z)}{r^2 s^2} ~=~ 2\, uvw / (r^2 s) \,+\, r^2 ~,
\end{equation}
where we have used the relation $r^4 = \Omega / s^2$ from Theorem \myref[2]{thm:myform} to get the right-hand side.
The result now follows from the orthogonality of the vectors $\VEC{\fit I\fit N}$, $\VEC{\fit I\fit L}$ \& $\VEC{\fit I\fit K}$, all of length $r$, to $\VEC{AN}$, $\VEC{AL}$ \& $\VEC{AK}$ respectively.
The proofs for the other three vertices are similar.
\end{proof}
\noindent Thus in the non-degenerate case the absolute values of the four \textquote{vectors} in Proposition \myref[3]{thm:collinear} are just the vertex-to-in-touch-point distances times $r$.
Moreover, in the degenerate case the ratios of those distances to one another as well as to the inter-vertex distances are also generically well defined.

This geometric interpretation of these quantities is of particular interest when combined with the following lemma, the proof of which will be left as an easy exercise in computer algebra (cf.~Lemma \myref[2]{thm:natinvnat}).
\begin{lemma}
The inverse natural parameters of a tetrahedron (non-degenerate or\linebreak[2] otherwise) can be written in terms of the square-roots $\hat u, \hat v, \,\ldots,\, \hat z$ of the natural par\-ameters as
\begin{subequations}
\begin{alignat}{4}
\tilde u\, s / 2 {}&{} ~=~ {}&{}  (\hat v \hat w - \hat x \hat y)^2 +\, \Omega_0 \Omega_1 
{}&{} ~=~ {}&{}  (\hat v \hat w + \hat x \hat y)^2 -\, \Omega_2 \Omega_3 \\ 
\tilde v\, s / 2 {}&{} ~=~ {}&{} (\hat u \hat w - \hat x \hat z)^2 +\, \Omega_0 \Omega_2 
{}&{} ~=~ {}&{} (\hat u \hat w + \hat x \hat z)^2 -\, \Omega_1 \Omega_3 \\ 
\tilde w\, s / 2 {}&{} ~=~ {}&{} (\hat u \hat v - \hat y \hat z)^2 +\, \Omega_0 \Omega_3 
{}&{} ~=~  {}&{} (\hat u \hat v + \hat y \hat z)^2 -\, \Omega_1 \Omega_2 \\ 
\tilde x\, s / 2 {}&{} ~=~ {}&{} (\hat u \hat y - \hat v \hat z)^2 +\, \Omega_0 \Omega_3 
{}&{} ~=~ {}&{} (\hat u \hat y + \hat v \hat z)^2 -\, \Omega_1 \Omega_2 \\ 
\tilde y\, s / 2 {}&{} ~=~ {}&{} (\hat u \hat x - \hat w \hat z)^2 +\, \Omega_0 \Omega_2 
{}&{} ~=~ {}&{} (\hat u \hat x + \hat w \hat z)^2 -\, \Omega_1 \Omega_3 \\ 
\tilde z\, s / 2 {}&{} ~=~ {}&{} (\hat v \hat x - \hat w \hat y)^2 +\, \Omega_0 \Omega_1 
{}&{} ~=~ {}&{} (\hat v \hat x + \hat w \hat y)^2 -\, \Omega_2 \Omega_3 
\end{alignat}
\end{subequations}
where $\Omega_k$ ($k = 0,1,2,3$) are the factors of $\Omega$ regarded as a polynomial in the square-roots of the natural parameters from Eq.~(\myref[3]{eq:omegafact}).
\end{lemma}
\noindent These ways of writing the inverse natural parameters make it quite clear why the sub\-stitutions used in the proof of Proposition \myref[3]{thm:collinear} yield the results that they do.

\smallskip 
In the non-degenerate case, this lemma together with Corollary \myref[2]{thm:r2d2} shows that e.g.
\begin{equation} \begin{aligned}
{}&{} u\tilde u\, s / 2 ~=~ (\hat u \hat v \hat w - \hat u \hat x \hat y)^2 + \hat u^2\, \Omega_0 \Omega_1 ~=~ (\hat u \hat v \hat w + \hat u \hat x \hat y)^2 - \hat u^2\, \Omega_2 \Omega_3 \\
\Longleftrightarrow\quad {}&{} \MAG{AB}^2 ~ \begin{alignedat}[t]{3}
=~ {}&{} \raisebox{1pt}{\large$\big($} \MAG{AN} - \MAG{BN} \raisebox{1pt}{\large$\big)$}^{\!2} + \frac{2\, u}{r^2 s}\, \Omega_0 \Omega_1 {}&{} ~=~ {}&{} \raisebox{1pt}{\large$\big($} \MAG{AN} + \MAG{BN} \raisebox{1pt}{\large$\big)$}^{\!2} - \frac{2\, u}{r^2 s}\, \Omega_2 \Omega_3 \\
=~ {}&{} \raisebox{1pt}{\large$\big($} \MAG{AL} - \MAG{BL} \raisebox{1pt}{\large$\big)$}^{\!2} +\, \frac{2\, u}{r^2 s}\, \Omega_0 \Omega_1 {}&{} ~=~ {}&{} \raisebox{1pt}{\large$\big($} \MAG{AL} + \MAG{BL} \raisebox{1pt}{\large$\big)$}^{\!2} -\, \frac{2\, u}{r^2 s}\, \Omega_2 \Omega_3
\end{alignedat}
\end{aligned} \end{equation}
and similarly for the other complementary products and squared inter-vertex distances, where the second line was obtained from the first by dividing through by $r^2 s / 2$ then applying Proposition \myref{thm:itp2vtx} and Corollary \myref[2]{thm:r2d2}, and the last line from $\MAG{AL} = \MAG{AN}$ \& $\MAG{BL} = \MAG{BN}$. 
This observation leads to the following:
\begin{corollary}
In a non-degenerate tetrahedron $\GEO{ABCD}$ with in-radius $r$, surface area $s$ (times $2$), natural parameters $u\fit, \,\ldots\,,\, z$ and with $\alpha, \beta, \gamma, \delta$ equal to the positive square-roots of the corresponding squares in Eq.~(\myref{eq:abcd-defs}), the products $\Omega_0 \Omega_k$ \& $\Omega_k \Omega_\ell$ ($k, \ell = 1,2,3$; $k \ne \ell$) are related to the inter-vertex distances $a \coloneq \MAG{AB}$, $b \coloneq \MAG{AC}$, $c \coloneq \MAG{AD}$, $d \coloneq \MAG{BC}$, $e \coloneq \MAG{BD}$, $f \coloneq \MAG{CD}$ as:
\begin{equation} \begin{alignedat}{4}
2\, u\, \Omega_0 \Omega_1 {}&{} ~=~ {}&{} s\, (ar - \alpha + \beta) (ar + \alpha - \beta) \,,\quad 
2\, u\, \Omega_2 \Omega_3 {}&{} ~=~ {}&{} s\, (\alpha + \beta + ar) (\alpha + \beta - ar) \,; \\ 
2\, v\, \Omega_0 \Omega_2 {}&{} ~=~ {}&{} s\, (br - \alpha + \gamma) (br + \alpha - \gamma) \,,\quad 
2\, v\, \Omega_1 \Omega_3 {}&{} ~=~ {}&{} s\, (\alpha + \gamma + br) (\alpha + \gamma - br) \,; \\ 
2\, w\, \Omega_0 \Omega_3{}&{}  ~=~ {}&{} s\, (cr - \alpha + \delta) (cr + \alpha - \delta) \,,\quad 
2\, w\, \Omega_1 \Omega_2 {}&{} ~=~ {}&{} s\, (\alpha + \delta + cr) (\alpha + \delta - cr) \,; \\ 
2\, x\, \Omega_0 \Omega_3 {}&{} ~=~ {}&{} s\, (dr - \beta + \gamma) (dr + \beta - \gamma) \,,\quad 
2\, x\, \Omega_1 \Omega_2 {}&{} ~=~ {}&{} s\, (\beta + \gamma + dr) (\beta + \gamma - dr) \,; \\ 
2\, y\, \Omega_0 \Omega_2 {}&{} ~=~ {}&{} s\, (er - \beta + \delta) (er + \beta - \delta) \,,\quad 
2\, y\, \Omega_1 \Omega_3 {}&{} ~=~ {}&{} s\, (\beta + \delta + er) (\beta + \delta - er) \,; \\ 
2\, z\, \Omega_0 \Omega_1 {}&{} ~=~ {}&{} s\, (fr - \gamma + \delta) (fr + \gamma - \delta) \,,\quad 
2\, z\, \Omega_2 \Omega_3 {}&{} ~=~ {}&{} s\, (\gamma + \delta + fr) (\gamma + \delta - fr) \,. 
\end{alignedat} \end{equation}
\end{corollary}
\noindent Note that multiplying the two equations on any given line in the above yields e.g.
\begin{align}
&{} 4\, u^2 \Omega_0 \Omega_1 \Omega_2 \Omega_3 ~=~ s^2 (\alpha + \beta + ar) (\alpha + \beta - ar) (ar - \alpha + \beta) (ar + \alpha - \beta) \\
\fit[-0.5em] \Longleftrightarrow~ {}&{} 4\, u^2 t^4 / s^2 ~=~ r^4 s^2 \Big( \raisebox{0.5pt}{\large$\big($} \smash{\MAG{AB}}^2 +\, \smash{\MAG{AN}}^2 +\, \smash{\MAG{BN}}^2 \raisebox{0.5pt}{\large$\big)$}^2 -\, 2\, \raisebox{0.5pt}{\large$\big($} \smash{\MAG{AB}}^4 +\, \smash{\MAG{AN}}^4 +\, \smash{\MAG{BN}}^4 \raisebox{0.5pt}{\large$\big)$} \Big) . \nonumber
\end{align}
This last equation can be derived directly by noting that the expression in parentheses on the right is four times the Cayley-Menger determinant of the three points $\GEO{A}, \GEO{B}, \GEO{N}$, or $4\, u^2$, and that $r^4 s^2 =\, \Omega \,=\, \Omega_0 \Omega_1 \Omega_2 \Omega_3 = t^4 / s^2$.

Returning now to our generic parametrization of the zeros, we use the relations in Eq.~(\myref[3]{eq:onprats}) to eliminate $x, y, z$ from the polynomial $\Omega(u,v,w,x,y,z)$, thereby obtaining
\begin{align} \label{eq:omegauvw}
& \Omega(u,v,w) ~\coloneq~ \alpha^2\beta^2\gamma^2\delta^2\, \Omega( u,v,w, w\, \beta\gamma/\alpha\delta, v\, \beta\delta/\alpha\gamma, u\, \gamma\delta / \alpha\beta ) ~= \\ &\quad \nonumber
\begin{aligned} (\gamma\delta\, u + \beta\delta v + \beta\gamma w) (-\gamma\delta\, u + \beta\delta v + \beta\gamma w) (\gamma\delta\, u - \beta\delta v + \beta\gamma w) (\gamma\delta\, u + \beta\delta v - \beta\gamma w) & \\
\eqcolon~ \Omega_0(u,v,w)\, \Omega_1(u,v,w)\, \Omega_2(u,v,w)\, \Omega_3(u,v,w) & ~, \end{aligned}
\end{align}
so that $\tau^4 = \varsigma^2\fit \Omega(u,v,w,x,y,z) = \varsigma^2\fit \Omega(u,v,w) / (\alpha\beta\gamma\delta)^2$ where $\tau$ \& $\varsigma$ are the assumed values of $t$ \& $s$, resp.
Similarly, on using Eq.~(\myref[3]{eq:onprats}) to eliminate $x, y, z$ from the expression $2\, (u+\cdots+z)$ for the surface area $s$, we get \pagebreak[2]
\begin{align} \label{eq:sigmauvw}
\varsigma ~=~ {}&{} 2\, (u + v + w + w\, \beta\gamma/\alpha\delta + v\, \beta\delta/\alpha\gamma + u\, \gamma\delta / \alpha\beta ) ~\eqcolon~ \Sigma(u,v,w) / (\alpha\beta\gamma\delta) \\ \nonumber
=~ {}&{} \big( 2\, (\alpha\beta + \gamma\delta)\, \gamma\delta\, u + 2\, (\alpha\gamma + \beta\delta)\, \beta\delta\, v + 2\, (\alpha\delta + \beta\gamma)\, \beta\gamma\, w \big) / (\alpha\beta\gamma\delta) ~.
\end{align}
In the non-degenerate case, these two equations together with $uvw = \alpha^2 \varsigma^2 / 2$ can be inverted to obtain a multi-valued relation between the natural parameters $u, v, w$ (which determine $x, y, z$ via Eq.~(\myref[3]{eq:onprats})) and $\alpha, \beta, \gamma, \delta$ together with $\varsigma$ \& $\tau$.
Numerical experiments with the Bertini program for solving polynomial systems of equations \cite{Bates:2013} indicate that these equations (which imply all those in Eq.~(\myref{eq:abcd-defs})) are highly ill-conditioned whenever the latter six parameters are anywhere close to those in a regular tetrahedron.
For a regular tetrahedron, the equations have a single positive real solution of multiplicity $6$.

A considerable simplification is possible in the degenerate ($\tau = 0$) case, providing that one gives $\alpha, \beta, \gamma, \delta$ the signs indicated in Proposition \myref[3]{thm:collinear} according to which factor $\Omega_k$ ($k = 1,2,3$) of $\Omega(u,v,w,x,y,z)$ vanishes in Eq.~(\myref[3]{eq:collinvecs}).
When this is done, the trivariate linear polynomial $\Omega_0(u,v,w)$ in Eq.~(\myref{eq:omegauvw}) will always be zero no matter which factor of $\Omega(u,v,w,x,y,z)$ in Eq.~(\myref[3]{eq:collinvecs}) vanishes.
Numerical experiments with Bertini indicate that the resulting system of two linear and one cubic equation in $u, v, w$, namely
\begin{equation} \label{eq:eqns4uvw}
0 ~=~ \Omega_0(u,v,w) ~,\quad 0 ~=~ \Sigma(u,v,w) - \alpha\beta\gamma\delta\, \varsigma ~,\quad 0 ~=~ uvw - \varsigma\, \alpha^2 / 2 ~,
\end{equation}
generically possesses three real solutions two of which satisfy $u, v, w, u', v', w' > 0$.
Upon substituting their values into Eq.~(\myref[3]{eq:onprats}), one also obtains $x, y, z, x', y', z' > 0$ satisfying
\begin{equation} \label{eq:yet-another}
u/u' ~=~ z/z' ~,\quad v/v' ~=~ y/y' ~,\quad w/w' ~=~ x/x ~,\quad (uvw) / (u'v'w') ~=~ 1 ~.
\end{equation}
This indicates that degenerate tetrahedra generically come in pairs related by these four equations as well as $s = s'$ (other pairings applicable to degenerate tetrahedra may be found in Appendix \myref{sec:invol}).

Some additional insight into the equations (\myref{eq:eqns4uvw}) can be obtained by using $0 =\, \Omega_0(u,v,w)$ to eliminate $w$ from $0 =\, \Sigma(u,v,w) - \alpha\beta\gamma\delta\, \varsigma$, which yields
\begin{align} \label{eq:res1}
0 ~=~ &\> \mrm{Res}( \Omega_0(u,v,w),\, \Sigma(u,v,w) - \alpha\beta\gamma\delta\, \varsigma;\, w) \\ \nonumber
=~ &\> \beta\gamma\delta\, \big( 2\, \big( \gamma\, (\alpha - \gamma)(\beta - \delta)\, u \,+\, \beta\, (\alpha - \beta) (\gamma - \delta) \, v \big) \,-\, \alpha\beta\gamma\, \varsigma \big) ,
\end{align}
where \textquote{$\mrm{Res}$} is the resultant of its first two arguments with respect to the third.
The overall factor of $\beta\gamma\delta$ is non-zero by our genericity assumption and may be dropped.
On also eliminating $w$ from the equations $0 =\, \Omega_0(u,v,w)$ and $0 = uvw - \varsigma \alpha^2 / 2$, one obtains
\begin{equation} \label{eq:res2}
0 ~=~ \mrm{Res}( \varsigma \alpha^2 / 2 - uvw,\, \Omega_0(u,v,w);\, w ) ~=~ \delta\, ( \gamma\, u + \beta\, v )\, uv \,+\, \varsigma\alpha^2\beta\gamma / 2 ~.
\end{equation}
Finally, upon taking the resultant of the right-hand sides of Eqs.~(\myref{eq:res1}) \& (\myref{eq:res2}) w.r.t.~$v$ and dropping the overall non-zero factors, we obtain a cubic equation for $u$:
\begin{align}
\fit[-0.5em] 0 ~=~ \Psi(u) ~=~ a_\Psi u^3 + b_\Psi u^2 + c_\Psi u + d_\Psi ~\coloneq~  4\, \gamma\delta\, (\alpha - \gamma) (\alpha - \delta) (\beta - \gamma) {}&{} (\beta - \delta)\, u^3 \nonumber \\
\begin{aligned} -\; 2\, \alpha\beta\gamma\delta\, \big( (\alpha - \gamma)(\beta - \delta) + (\alpha - \delta)(\beta - \gamma) \big)\, \varsigma\, u^2 & \\[1pt]
+\; \alpha^2\beta^2\gamma\delta\, \varsigma^2\, u \,+\, 2\, \alpha^2\beta^2 (\alpha - \beta)^2 (\gamma - \delta)^2\, \varsigma & \end{aligned} & \label{eq:cubic} \\ \nonumber
=~ \gamma \delta\, \big( 2\, (\alpha - \delta) (\beta - \gamma)\, u - \alpha\beta\, \varsigma \big) \big( 2\, (\alpha - \gamma) (\beta - \delta)\, u - \alpha\beta\, \varsigma {}&{} \big)\, u + d_\Psi
\end{align}
Clearly the constant term is positive, but the signs of the remaining coefficients depend on the relative signs \& magnitudes of $\alpha, \beta, \gamma, \delta$.
Numerical examples nevertheless suggest that when $\varsigma$ is sufficiently large (as specified below) this cubic always has either $2$ or $3$ positive roots, depending on if $a_\Psi > 0$ or $a_\Psi < 0$, but that in the latter case putting the largest root into Eqs.~(\myref{eq:res1}) \& (\myref{eq:res2}) yields negative values for $v, w$.

The discriminant of this cubic is $\varsigma^2$ times a quadratic polynomial in $\varsigma^2$, namely:
\begin{align} 
18\, a_\Psi b_\Psi c_\Psi d_\Psi - 4\, (b_\Psi)^3 d_\Psi + (b_\Psi c_\Psi)^2 - 4\, a_\Psi (c_\Psi)^3 - 27\, (a_\Psi d_\Psi)^2 ~= & \nonumber \\
\begin{aligned}
4\, \alpha^6 \beta^6 \gamma^4 \delta^4 (\alpha \!-\! \beta)^2 (\gamma \!-\! \delta)^2\, \varsigma^6 \,+\, 32\, \alpha^5 \beta^5 \gamma^3 \delta^3 \big( (\alpha \!-\! \gamma) (\beta \!-\! \delta) + (\alpha \!-\! \delta) (\beta \!-\! \gamma) \big) \cdots & \\ \cdots
\big( (\alpha \!-\! \beta) (\gamma \!-\! \delta) + (\alpha \!-\! \gamma) (\beta \!-\! \delta) \big) \big( (\alpha \!-\! \beta) (\gamma \!-\! \delta) - (\alpha \!-\! \delta) (\beta \!-\! \gamma) \big) (\alpha \!-\! \beta)^2 (\gamma \!-\! \delta)^2\, \varsigma^4 & \end{aligned} & \nonumber \\
-\; 1728\, \alpha^4\beta^4\gamma^2\delta^2 (\alpha \!-\! \beta)^4 (\alpha \!-\! \gamma)^2 (\alpha \!-\! \delta)^2 (\beta \!-\! \gamma)^2 (\beta \!-\! \delta)^2 (\gamma \!-\! \delta)^4 \, \varsigma^2 & \nonumber \\
& 
\end{align}
Since the constant term is negative while the leading coefficient is positive, Descartes' rule of signs shows that this quadratic has a single positive root $\varrho_+( \alpha, \beta, \gamma, \delta )$, and for all $\varsigma^2$ exceeding this root the cubic in Eq.~(\myref{eq:cubic}) will have three real roots.
This suggests that the equations (\myref{eq:eqns4uvw}) are only valid when $\varsigma$ observes this lower bound and, in accord with this expectation, it has been observed that in random examples with $\varsigma^2 < \varrho_+$ Bertini returns only a single non-positive solution.

Since negating all the parameters $\alpha, \beta, \gamma, \delta$ does not change the equations, the above observations lead to the following:
\begin{conjecture}
The equations (\myref{eq:eqns4uvw}) \& (\myref[3]{eq:onprats}) generically define a $2\!:\!2$ relation between the natural parameters of degenerate tetrahedra and the semi-algebraic set $\big\{ \alpha, \beta, \gamma, \delta, \varsigma \in \mbb R \linebreak[2] \,\bigm\vert\, |\alpha| + |\beta| + |\gamma| + |\delta| > |\alpha + \beta + \gamma + \delta| ,~ \alpha\beta\gamma\delta > 0 ,~ \varsigma > \sqrt{\varrho_+(\alpha,\beta,\gamma,\delta)} \big\}$.
\end{conjecture}

The author thanks Silviana Amethyst for assistance with the Bertini calculations described herein.

\section{Two involutions on tetrahedra, non-degenerate and otherwise} \label{sec:invol}
The late Czech mathematician Miroslav Fiedler has studied an interesting involution between $n$-simplices, which yields what he called the \emph{inverse\/} simplex \cite{Fiedler:2011}.
It may be viewed as a Euclidean specialization of the more general polarity involution defined on convex polytopes in affine geometry, see e.g.~\S 2.3 in Ref.~\cite{Ziegler:2000}.
In the non-degenerate case, Fiedler's inverse simplex is most simply obtained by computing the (Moore-Penrose) pseudo-inverse of the \EMPH{lineal} Gram matrix formed from the inner products of the vectors from the centroid of the simplex to the vertices thereof, and then computing coordinates for the inverse's vertices via diagonalization, exactly as one would recover coordinates for the original simplex from its lineal Gram matrix \cite{Borg:2005,Crippen:1978,Gower:2004}.
Fiedler further showed that the pseudo-inverse of the full (hyper-)areal Gram matrix of a non-degenerate $n$-simplex, over the square of $n$ times its (hyper-)volume, is equal to this lineal Gram matrix (and vice versa), which provides an alternative geometric interpretation of the latter's pseudo-inverse.
Finally, he extended his inverse involution to degenerate $n$-simplices via a matrix construction called a \textquote{generalized biorthogonal system,} which involves a Gale transform \cite{Ziegler:2000} of the corresponding affinely dependent configuration.

For a non-degenerate tetrahedron, Fiedler's inverse may also be obtained by taking the pseudo-inverse of the full $4\times4$ areal Gram matrix $\mbf G_\msf{ext}$ as defined in the proof of Proposition \myref[3]{thm:planar}, multiplying by $16$, and then computing vertex coordinates from a $3\times3$ principal submatrix thereof as in the proof of Theorem \myref[1]{thm:bdsmc} (see Fig.~\myref{fig:fiedler}).
For a (conventional) degenerate tetrahedron (in the finite Euclidean plane), however, Fiedler's inverse involution must be based upon the lineal Gram matrix rather than the areal, since the areal Gram matrix has rank $1$ and only defines the configuration up to area-preserving affine transformations.
In the general rank $2$ case, moreover, Fiedler's definition of the inverse involution does not apply at all, because the inter-vertex distances are generally infinite so that the lineal Gram matrix does not even exist.

\begin{figure}
\input{Fig_IV-2.tex}
\label{fig:fiedler}
\end{figure}

This appendix will further explore the nature of Fiedler's inverse involution for non-degenerate tetrahedra, adapt it to the areal Gram matrices of degenerate tetrahedra with $s^2\fit\Omega = 0$, and define an apparently new involution that is specific to the $n = 3$ tetrahedral case.
For these purposes it is more convenient to focus upon the interior areal Gram matrix $\mbf G_\msf{int}$ as defined in Eq.~(\myref[2]{eq:intgramat}), rather than the exterior.
These two matrices are related by a simple equivalence transformation, namely
\begin{equation}
\mbf G_\msf{ext} ~=~ \mbf{MG}_\msf{int}\mbf M^\top \quad\text{where}\quad \mbf M ~\coloneq~ \sfrac12 \left[ \begin{smallmatrix} \hit[1.5ex] ~1 & -1 & -1~ \\[2pt] ~1 & ~1 & ~1~ \\[2pt] -1 & -1 & ~1~ \\[2pt] -1 & ~1 & -1~ \dit[1ex] \end{smallmatrix} \right] ~,
\end{equation}
which may be derived from Eq.~(\myref[1]{eq:havel0}).
Since the Moore-Penrose pseudo-inverse of $\mbf M$ is $\mbf M^\ominus = \mbf M^\top$, in the non-degenerate case the pseudo-inverse of $\mbf G_\msf{ext}$ is simply $\mbf G_\msf{ext}^\ominus = \mbf{M\fit G}_\msf{int}^{-1}\fit \mbf M^\top$.
Because the pseudo-inverse of a diagonal matrix is diagonal, this formulation of Fiedler's inverse together with Remark \myref[1]{rem:int-face-obs} shows immediately that the inverse of an equi-facial tetrahedron is again equi-facial.
It can also be shown that the inverse $\mbf G_{\smash{\msf{int}}}^{-1}$ is nothing more (or less) than the lineal Gram matrix of the \TDEF{oriented bimedians\/} times $t^{-2}$, specifically $\VEC{ZU}$, $\VEC{VY}$, $\VEC{XW}$ in the notation of Fig.~\myref[1]{fig:medoct}, or:
\begin{subequations} \begin{align}
\VEC{ZU} ~=~ \tfrac12\, \big( \GEO{A} \Vplus \GEO{B} \Vdiff \GEO{C} \Vdiff \GEO{D} \big) ~=~ \Vdiff & \tfrac12\, \big( \VEC{AC} \Vplus \VEC{BD} \big) ~=~ \Vdiff\tfrac12\, \big( \VEC{AD} \Vplus \VEC{BC} \big) \\
\VEC{VY} ~=~ \tfrac12\, \big( \GEO{B} \Vplus \GEO{D} \Vdiff \GEO{A} \Vdiff \GEO{C} \big) ~=~ \fit[0.6em]& \tfrac12\, \big( \VEC{AB} \Vplus \VEC{CD} \big) ~=~ \quad\tfrac12\, \big( \VEC{AD} \Vdiff \VEC{BC} \big) \\
\VEC{XW} ~=~ \tfrac12\, \big( \GEO{A} \Vplus \GEO{D} \Vdiff \GEO{B} \Vdiff \GEO{C} \big) ~=~ \fit[0.6em]& \tfrac12\, \big( \VEC{BD} \Vdiff \VEC{AC} \big) ~=~ \quad\tfrac12\, \big( \VEC{CD} \Vdiff \VEC{AB} \big)
\end{align} \end{subequations}
This makes it possible to compute $\mbf G_\msf{int}^{-1}$ without explicitly taking a matrix inverse.

The determinant of the Gram matrix of the oriented bimedians is easily shown to be $t^2/4$, whereas $\mrm{det}\big( \mbf G_\msf{int} \big) = 4\, t^4 \Leftrightarrow \mrm{det}\big( \mbf G_\msf{int}^{-1} \big) = 1/(4\, t^4)$.
It follows that the volume of the inverse tetrahedron is $4/t$, which in turn implies that $\VCP{AB}{CD} \,/\, t$, $\VCP{AC}{BD} \,/\, t$, $\VCP{AD}{BC} \,/\, t$ can be viewed as the oriented bimedians of the inverse.
The positions of $\GEO{B}$, $\GEO{C}$, $\GEO{D}$ relative to $\GEO{A}$ may be obtained from the original's bimedians as
\begin{equation}
\VEC{AB} ~=~ \VEC{VY} \Vdiff \VEC{XW} ~,\qquad \VEC{AC} ~=\, \Vdiff\,\VEC{ZU} \Vdiff \VEC{XW} ~,\qquad \VEC{AD} ~=~ \VEC{VY} \Vdiff \VEC{ZU} ~,
\end{equation}
and similarly the positions of the inverse's vertices $\GEO{B'}$, $\GEO{C'}$, $\GEO{D'}$ relative to $\GEO{A'}$ may be obtained as
\begin{subequations} \begin{gather}
\VEC{A'B'} ~=~ \quad\sfrac1{t}\, \big( \VCP{AC}{BD} \,\Vdiff\, \VCP{AD}{BC} \big) ~, \\
\VEC{A'C'} ~=~ \Vdiff\sfrac1{t}\, \big( \VCP{AB}{CD} \,\Vplus\, \VCP{AD}{BC} \big) ~, \\
\VEC{A'D'} ~=~ \quad\sfrac1{t}\, \big( \VCP{AC}{BD} \,\Vdiff\, \VCP{AB}{CD} \big) ~.
\end{gather} \end{subequations}
This shows that one can obtain the inverse tetrahedron not only without explicitly computing $\mbf G_{\smash{\msf{int}}}^{-1}\,$, but also without converting it to $\mbf G_\msf{ext}^\ominus$ and diagonalizing that matrix.

All of these observations, of course, only apply to the non-degenerate case; when $t = 0$ one can still compute $\mbf G_\msf{ext}^\ominus = \mbf{M\fit G}_{\smash{\msf{int}}}^\ominus \mbf M^\top$, but these pseudo-inverses are no longer $(t/4)^2$ times the corresponding lineal Gram matrices of the inverse tetrahedron (for obvious reasons).
What nevertheless \EMPH{can} be done is to take the square roots of the diagonal entries of $\mbf G_\msf{ext}^\ominus$ \& $\mbf G_{\smash{\msf{int}}}^\ominus$ to be the values of the seven areas in another degenerate ($s^2\, \Omega = 0$) tetrahedron.
These determine what will be called the \TDEF{reciprocal} tetrahedron, to distinguish it from Fiedler's inverse for (conventional) degenerate tetrahedra.
It coincides with Fiedler's inverse only in the non-degenerate case (with suitable multiples of $t^2$ or $1/t^2$ thrown in, as in Fig.~\myref{fig:fiedler}).

We now turn to another, completely different involution on the set of tetrahedra as para\-meterized by their natural parameters.
\begin{definition}
Given a tetrahedron, non-degenerate or otherwise, with natural parameters $u, v, w, x, y, z\fit$, its \TDEF{twin} is the tetrahedron obtained by swapping their values in all three opposite pairs, i.e.~$u \leftrightarrow z\fit$, $v \leftrightarrow y$, $w \leftrightarrow x$, whereas the given tetrahedron will be referred to as the \TDEF{base} tetrahedron for contrast.
\end{definition}
\begin{remark}
Note the twin as defined here differs from that in Ref.~\cite{Court:1935}; the latter is just the mirror image of the tetrahedron superimposed upon it so that both have the same circumscribing parallelopiped.
Note also that swapping the values of opposite pairs of inter-vertex \EMPH{distances} in a non-degenerate tetrahedron does not in general yield a distance matrix that can be realized in Euclidean space.
It can be shown that the tetrahedra obtained by swapping the values of any single pair of opposite natural parameters are equivalent to the above twin up to vertex relabeling, while those obtained from even numbers of such swaps are equivalent to the base up to relabeling.
Thus there is no loss of generality in restricting ourselves to the more symmetric triple swaps.
\end{remark}
This \textquote{twinning} involution is clearly quite different from Fiedler's inverse even for non-degenerate tetrahedra, and further has no obvious relation to the aforementioned polarity involution of convex geometry.
Indeed it is not immediately clear how the twin as defined here can be extended to higher (or lower) dimensions.
In this regard it bears some resemblance to the Regge symmetry of tetrahedra as described in e.g.~Ref.~\cite{Akopyan:2019}, although that is based upon the edge lengths rather than the facial areas.

The areal properties of the twin tetrahedron are summarized in the following:
\begin{proposition}
The areas $\MAG{A'B'C'}$, $\ldots\,$, $\MAG{A'D'|B'C'}$, surface area $s'$, volume $t'$, in-radius $r' = t'/s'$ and inverse natural parameters $\tilde u'$, $\ldots$, $\tilde z'$ of the twin tetrahedron with $u' \coloneq z$, $z' \coloneq u$, $v' \coloneq y$, $y' \coloneq v$, $w' \coloneq x$, $x' \coloneq w$ satisfy:
\begin{subequations} \begin{gather}
\begin{aligned}
2 \MAG{A'B'C'} ~=~ {}&{} \MAG{ABD} \,+\, \MAG{ACD} \,+\, \MAG{BCD} \,-\, \MAG{ABC} \,, \\
2 \MAG{A'B'D'} ~=~ {}&{} \MAG{ACD} \,+\, \MAG{BCD} \,+\, \MAG{ABC} \,-\, \MAG{ABD} \,, \\
2 \MAG{A'C'D'} ~=~ {}&{} \MAG{BCD} \,+\, \MAG{ABC} \,+\, \MAG{ABD} \,-\, \MAG{ACD} \,, \\
2 \MAG{B'C'D'} ~=~ {}&{} \MAG{ABC} \,+\, \MAG{ABD} \,+\, \MAG{ACD} \,-\, \MAG{BCD} \,;
\end{aligned} \label{eq:twin_a} \\[3pt]
\fit[-0.5em] \MAG{A'B'|C'D'} \,=\, \MAG{AB|CD} \,,~ \MAG{A'C'|B'D'} \,=\, \MAG{AC|BD} \,,~ \MAG{A'D'|B'C'} \,=\, \MAG{AD|BC} ; \label{eq:twin_b} \\[1pt]
s' ~=~ s \,,\quad t' ~=~ t \,,\quad r' ~=~ r \,; \label{eq:twin_c} \\
\tilde u' \:=\: \tilde u \,,~ \tilde v' \:=\: \tilde v \,,~ \tilde w' \:=\: \tilde w \,,~ \tilde x' \:=\: \tilde x \,,~ \tilde y' \:=\: \tilde y \,,~ \tilde z' \:=\: \tilde z ~. \label{eq:twin_d}
\end{gather} \end{subequations}
\end{proposition}
\begin{proof}
Equations (\myref{eq:twin_a}) \& (\myref{eq:twin_b}) follow from Eqs.~(\myref[2]{eq:idb}), (\myref[2]{eq:idc}) \& (\myref[2]{eq:ide}) in Corollary \myref[2]{thm:ids}, respectively, while the $s' = s$ in Eq.~(\myref{eq:twin_c}) follows from Eq.~(\myref[2]{eq:ida}).
The volumes of the base tetrahedron and its twin are equal since $t^4 = s^2\, \Omega(u,v,w,x,y,z) = s^2\, \Omega(z,y,x,w,v,u) = (t')^4$ by Theorem \myref[2]{thm:myform}, from which it follows that $r' = t'/s' = t/s = r$.
Equation (\myref{eq:twin_d}) is a consequence of $s' = s$ together with the formulae established in Lemma \myref[2]{thm:natinvnat}.
\end{proof}

In the non-degenerate case, it further follows from this proposition together with Corollary \myref[2]{thm:r2d2} that the squared inter-vertex distances of the twin are given by:
\begin{equation} \begin{aligned}
\MAG{A'B'}^2 ~=~ z\tilde u/r^2 \,,\quad\; & \MAG{A'C'}^2 ~=~ y\tilde v/r^2 \,,\quad \MAG{A'D'}^2 ~=~ x\tilde w/r^2 \,, \\
\MAG{B'C'}^2 ~=~ w\tilde x/r^2 \,,\quad & \MAG{B'D'}^2 ~=~ v\tilde y/r^2 \,,\quad \MAG{C'D'}^2 ~=~ u\tilde z/r^2 \,.
\end{aligned} \end{equation}
This shows that the twinning involution preserves the products of opposite pairs of distances:
\begin{equation}
\MAG{A'B'} \! \MAG{C'D'} \,=\, \MAG{AB} \! \MAG{CD} ,~ \MAG{A'C'} \! \MAG{B'D'} \,=\, \MAG{AC} \! \MAG{BD} ,~ \MAG{A'D'} \! \MAG{B'C'} \,=\, \MAG{AD} \! \MAG{BC}
\end{equation}
Since twinning also preserves the interior areas and
\begin{equation*}
16 \MAG{AB|CD}^2 ~=~ \MAG{AB}^2 \MAG{CD}^2 -\, \sfrac14\, \Big( \MAG{AD}^2 + \MAG{BC}^2 - \MAG{AC}^2 - \MAG{BD}^2 \Big)^{\!2}
\end{equation*}
etc.~by Eqs.~(\myref[1]{eq:talata}) \& (\myref[1]{eq:talata2}), it further follows that it preserves the dot products of the inter-vertex vectors between opposite pairs of vertices:
\begin{equation}
\fit[-1em] \VEC{A'B'} \fit\Dprod\, \VEC{C'D'} ~=~ \VEC{AB} \fit\Dprod\fit \VEC{CD} \,,\fit[0.5em] \VEC{A'C'} \fit\Dprod\, \VEC{B'D'} ~=~ \VEC{AC} \fit\Dprod\fit \VEC{BD} \,,\fit[0.5em] \VEC{A'D'} \fit\Dprod\, \VEC{B'C'} ~=~ \VEC{AD} \fit\Dprod\fit \VEC{BC}
\end{equation}
This shows immediately that the twin of an orthocentric tetrahedron is again orthocentric.
It also follows from Eq.~(\myref[2]{eq:idb}) of Corollary \myref[2]{thm:ids} that a tetrahedron is equi-facial if \& only if $u = z$, $v = y$ \& $w = x$, so that the fixed points of the twinning involution are exactly the equi-facial tetrahedra (non-degenerate and otherwise).

One striking difference between the twinning and reciprocal involutions is that while the latter always preserves the common rank of the areal Gram matrices, the former does not necessarily do so.
Even though twinning preserves the volume and hence always carries non-degenerate tetrahedra to the same, it can map rank $1$ tetrahedra to rank $2$ and vice versa.
An examination of the patterns of natural and inverse natural parameters vanishing in rank $1$ tetrahedra (Fig.~\myref{fig:regions-labeled}) in fact shows that whereas twinning maps cases ($\msf 4, \msf 5, \msf 6, \msf 7$), ($\msf 8, \msf 9, \msf{10}, \msf{11}$) \& ($\msf{12}, \msf{13}, \msf{14}, \msf{15}$) to themselves, applied to the cases $\msf 0$, $\msf 1$, $\msf 2$ \& $\msf 3$ it produces a rank $2$ tetrahedron.
These rank $2$ tetrahedra are those at the lowest level of the rank $2$ sub-hierarchy in Fig.~\myref[3]{fig:lattice} wherein three complementary products vanish; the other two levels in that sub-hierarchy are clearly preserved by twinning.

Via random numerical examples, it is easily shown that the reciprocal and twinning involutions do \EMPH{not} commute.
Similarly, it can be shown that their compositions (in either order) are not themselves involutions, so that the group they generate is not the Klein $4$-group.
The following conjecture is most likely wistful thinking, but if true would show that the reciprocal and twinning involutions are somehow connected.
The alternative is that alternatively applying these two involutions executes a kind of a chaotic trajectory through \textquote{tetrahedron space,} which would also be rather interesting.
\begin{conjecture}
The two compositions of the reciprocal and twinning involutions have finite orders.
\end{conjecture}
\noindent Note that since these compositions are conjugate, their orders are necessarily equal.

Yet another involution, at least on degenerate tetrahedra, is defined by the condition $s' = s$ together with Eqs.~(\myref{eq:yet-another}) of Appendix \myref{sec:2to2}; it would also be interesting to study how it relates to the two analyzed above.
With a bit of imagination \cite{Hilbert:1952}, it is not hard to come up with yet more discrete relations between tetrahedra based on the areas or areal vectors of their interior and exterior faces.
One could, for example, set the areal Gram matrix $\mbf G_\msf A$ at $\GEO{A}$ to the interior areal Gram matrix $\mbf G_\msf{int}$ of another tetrahedron (perhaps after a row/column permutation).
Clearly instead setting the areal Gram matrix at $\GEO{B}, \GEO{C}$ or $\GEO{D}$ to $\mbf G_\msf{int}$ yields the same tetrahedron up to relabeling.
The same is \EMPH{not} true, however, if one sets $\mbf G_\msf{int} \coloneq \mbf G_\msf A$, $\mbf G_\msf B$, $\mbf G_\msf C$ or $\mbf G_\msf D\fit$; in fact,
iterating on these latter identifications appears to produce a directed tree of non-equivalent tetrahedra, each node of which has out-degree four!
Is this tree infinite, or do the directed paths therein all converge, and if so, to what?
The one thing that seems certain is that there are infinitely many questions about the humble tetrahedron left to be explored \textellipsis

{\small

} 



\bigskip
\begin{center}
\medskip\rule{70pt}{2pt}\medskip
\end{center}

\medskip
\begin{displayquote}\small
Before I begin to talk to you about the sizes and shapes of things, I am going to make a request that may seem somewhat strange.
I am going to ask you to forget that you have ever lived until this moment.
It is not that I am going to tell you anything new, that you did not know before; for I am merely going to remind you of a lot of things that you have known familiarly for years.
Only I want you to observe them all quite freshly over again, as if you had not seen them before \textellipsis~for geometry, you know, is the gate of science, and the gate is so low and small that one can only enter it as a child.

\smallskip\raggedleft\textit{\textquote{Seeing and Thinking,} by William Kingdon Clifford, published 1883}
\end{displayquote}

\medskip
\begin{displayquote}\small
For my part I know nothing with any certainty, but the sight of the stars makes me dream.

\vspace{-1ex}\raggedleft\textit{Vincent Van Gogh, 1888}
\end{displayquote}



\end{document}